\documentclass[a4paper,twoside,10pt]{amsart}
\usepackage{amsmath,amsfonts,amssymb,amsthm}
\newtheorem{theorem}{Theorem}[section]
\usepackage{mathrsfs}
\usepackage{color,graphicx}
\usepackage[a4paper]{geometry}
\numberwithin{equation}{section}

\author[G. Nemes]{Gerg\H{o} Nemes}
\address{Central European University, Department of Mathematics and its Applications, H-1051 Budapest, N\'ador utca 9, Hungary}
\email{nemesgery@gmail.com}


\dedicatory{To the memory of Frank W. J. Olver (1924--2013)}

\begin{document}

\title[Resurgence of Hankel and Bessel functions]{The resurgence properties of\\ the large order asymptotics of\\ the Hankel and Bessel functions}

\begin{abstract} The aim of this paper is to derive new representations for the Hankel and Bessel functions, exploiting the reformulation of the method of steepest descents by M. V. Berry and C. J. Howls (Berry and Howls, Proc. R. Soc. Lond. A \textbf{434} (1991) 657--675). Using these representations, we obtain a number of properties of the large order asymptotic expansions of the Hankel and Bessel functions due to Debye, including explicit and numerically computable error bounds, asymptotics for the late coefficients, exponentially improved asymptotic expansions, and the smooth transition of the Stokes discontinuities.
\end{abstract}
\maketitle

\section{Introduction and main results}\label{section1}

The large $\nu$ asymptotics of the Hankel functions $H_\nu^{\left(1\right)}\left(\nu x\right)$, $H_\nu^{\left(2\right)}\left(\nu x\right)$, and the Bessel functions $J_\nu \left(\nu x\right)$, $Y_\nu \left(\nu x\right)$, beside theoretical and practical importance, have historical significance. In 1909, Debye \cite{Debye} developed the method of steepest descents to deduce asymptotic expansions for these functions with fixed $x>0$ as $\nu \to \infty$. Since then, these expansions become standard textbook examples to illustrate various techniques, such as the method of steepest descents itself or the method of stationary phase (see, for example, Copson \cite[pp. 34--35 and pp. 75--82]{Copson}; Olver \cite[pp. 133--134]{Olver}; Wong \cite[pp. 94--103]{Wong}). Meijer \cite{Meijer} and Olver \cite[p. 382]{Olver} were able to derive error bounds for Debye's expansions when $x>1$. Meijer's approach was based on the Lagrange inversion theorem while Olver used differential equation methods. Gatteschi \cite{Gatteschi} estimated the remainder of Debye's series when $x=1$. Dingle \cite[pp. 166--174]{Dingle} investigated the resurgence properties of Debye's expansions and obtained asymptotic approximations for their late terms. Nevertheless, the derivation of his results is based on interpretive, rather than rigorous, methods.

The drawback of Debye's expansions is their non-uniformity in $x$. Each expansion has different form according to whether $0<x<1$, $x=1$ or $x>1$. In the first and the third cases, it is also required that $\lim _{\nu  \to \infty } \left|\nu \right|^{\frac{2}{3}} \left| {x - 1} \right| =  + \infty$. In 1952, Olver \cite{Olver2} modified the expansions corresponding to $x=1$, obtaining asymptotic series in the so-called transition regions, i.e., when $\left|\nu\right| ^{\frac{2}{3}} \left| {x - 1} \right|$ is bounded (see also Sch\"{o}be \cite{Schobe}). Two years later, he published new asymptotic expansions for the Hankel and Bessel functions that are uniformly valid for all $x>0$ \cite{Olver3}.

It might be supposed therefore that there is no need for the non-uniform expansions of Debye and that little could be added to their well-established theory. The uniformity of the expansion is an important property, however, in general, uniform asymptotic expansions have much more complicated forms than non-uniform expansions. Olver's uniform asymptotic expansions involve the Airy functions, whereas Debye's ones build from elementary functions. In addition, we have simple recurrence formulas for the coefficients in the latter ones. The coefficients in Olver's uniform expansions are slightly more complicated, and their computation requires the ones in Debye's series. For the asymptotic properties of these coefficients, see Howls and Olde Daalhuis \cite{Howls}.

The aim of this paper is to establish new resurgence-type integral representations for the functions $H_\nu^{\left(1\right)}\left(\nu x\right)$, $H_\nu^{\left(2\right)}\left(\nu x\right)$, $J_\nu \left(\nu x\right)$ and $Y_\nu \left(\nu x\right)$ when $x\geq 1$. Our derivation is based on the reformulation of the method of steepest descents by Berry and Howls \cite{Berry} (see also Boyd \cite{Boyd} and Paris \cite[pp. 94--99]{Paris}). Here resurgence has to be understood in the sense of Berry and Howls, meaning that the function reappears in the remainder of its own asymptotic series. Using these representations, we obtain several new properties of Debye's classical expansions, including explicit and numerically computable error bounds, asymptotics for the late coefficients, exponentially improved asymptotic expansions, and the smooth transition of the Stokes discontinuities. Our analysis also provides a rigorous treatment of Dingle's formal results.

Our first theorem describes the resurgence properties of the asymptotic expansions of $H_\nu^{\left(1\right)}\left(\nu x\right)$, $H_\nu^{\left(2\right)}\left(\nu x\right)$, $J_\nu \left(\nu x\right)$ and $Y_\nu \left(\nu x\right)$ for $x > 1$. We employ the usual substitution $x = \sec \beta$ with an appropriate $0 < \beta <\frac{\pi}{2}$. The notations follow the ones given in \cite[p. 231]{NIST}. Throughout this paper, empty sums are taken to be zero.

\begin{theorem}\label{thm1} Let $0 < \beta <\frac{\pi}{2}$ be a fixed acute angle, and let $N$ be a non-negative integer. We have
\begin{equation}\label{eq21}
H_\nu ^{\left( 1 \right)} \left( {\nu \sec \beta } \right) = \frac{{e^{i\nu \left( {\tan \beta  - \beta } \right) - \frac{\pi}{4}i} }}{{\left( {\frac{1}{2}\nu \pi \tan \beta } \right)^{\frac{1}{2}} }}\left( {\sum\limits_{n = 0}^{N - 1} {\left( { - 1} \right)^n \frac{{U_n \left( {i\cot \beta } \right)}}{{\nu ^n }}}  + R_N^{\left(H\right)} \left( {\nu ,\beta } \right)} \right)
\end{equation}
for $-\frac{\pi}{2} < \arg \nu < \frac{3\pi}{2}$;
\begin{equation}\label{eq17}
H_\nu ^{\left( 2 \right)} \left( {\nu \sec \beta } \right) = \frac{{e^{ - i\nu \left( {\tan \beta  - \beta } \right) + \frac{\pi}{4}i} }}{{\left( {\frac{1}{2}\nu \pi \tan \beta } \right)^{\frac{1}{2}} }}\left( {\sum\limits_{n = 0}^{N - 1} {\frac{{U_n \left( {i\cot \beta } \right)}}{{\nu ^n }}}  + R_N^{\left(H\right)} \left( { \nu e^{\pi i}  ,\beta } \right)} \right)
\end{equation}
for $-\frac{3\pi}{2} < \arg \nu < \frac{\pi}{2}$;
\begin{equation}\label{eq22}
J_\nu  \left( {\nu \sec \beta } \right) = \left( {\frac{2}{{\nu \pi \tan \beta }}} \right)^{\frac{1}{2}} \left( {\cos \xi \sum\limits_{n = 0}^{N - 1} {\frac{{U_{2n} \left( {i\cot \beta } \right)}}{{\nu ^{2n} }}}  - i\sin \xi \sum\limits_{n = 0}^{N - 1} {\frac{{U_{2n + 1} \left( {i\cot \beta } \right)}}{{\nu ^{2n + 1} }}}  +  R_N^{\left(J\right)} \left( {\nu ,\beta } \right)} \right),
\end{equation}
\begin{equation}\label{eq23}
Y_\nu  \left( {\nu \sec \beta } \right) = \left( {\frac{2}{{\nu \pi \tan \beta }}} \right)^{\frac{1}{2}} \left( {\sin \xi \sum\limits_{n = 0}^{N - 1} {\frac{{U_{2n} \left( {i\cot \beta } \right)}}{{\nu ^{2n} }}}  + i\cos \xi \sum\limits_{n = 0}^{N - 1} {\frac{{U_{2n + 1} \left( {i\cot \beta } \right)}}{{\nu ^{2n + 1} }}}  + R_N^{\left(Y\right)} \left( {\nu ,\beta } \right)} \right)
\end{equation}
for $\left|\arg \nu\right| < \frac{\pi}{2}$ with $\xi  = \nu \left( {\tan \beta  - \beta } \right) - \frac{\pi }{4}$. The coefficients $U_n \left( {i\cot \beta } \right)$ are given by
\begin{gather}\label{eq13}
\begin{split}
U_n \left( {i\cot \beta } \right) & = \left( { - 1} \right)^n \frac{{\left( {i\cot \beta } \right)^n }}{{2^n n!}}\left[ {\frac{{d^{2n} }}{{dt^{2n} }}\left( {\frac{1}{2}\frac{{t^2 }}{{i\cot \beta \left( {t - \sinh t} \right) + \cosh t - 1}}} \right)^{n + \frac{1}{2}} } \right]_{t = 0}\\
& = \frac{i^n}{2\left( {2\pi \cot \beta } \right)^{\frac{1}{2}}}\int_0^{ + \infty } {t^{n - \frac{1}{2}} e^{ -t \left( {\tan \beta  - \beta } \right)} \left( {1 + e^{ - 2\pi t} } \right)i H_{it}^{\left( 1 \right)} \left( {it\sec \beta } \right)dt} .
\end{split}
\end{gather}
The remainder terms can be expressed as
\begin{equation}\label{eq12}
R_N^{\left(H\right)} \left( {\nu ,\beta } \right) =  \frac{1}{{2\left( {2\pi \cot \beta } \right)^{\frac{1}{2}} \left(i\nu\right)^N }}\int_0^{ + \infty } {\frac{{t^{N - \frac{1}{2}} e^{ - t\left( {\tan \beta  - \beta } \right)} }}{{1 + it/\nu }}\left( {1 + e^{ - 2\pi t} } \right)i H_{it}^{\left( 1 \right)} \left( {it\sec \beta } \right)dt} ,
\end{equation}
\begin{gather}\label{eq29}
\begin{split}
R_N^{\left(J\right)} \left( {\nu ,\beta } \right) = \; & \cos \xi \frac{{\left( { - 1} \right)^N }}{{2\left( {2\pi \cot \beta } \right)^{\frac{1}{2}} \nu ^{2N} }}\int_0^{ + \infty } {\frac{{t^{2N - \frac{1}{2}} e^{ - t\left( {\tan \beta  - \beta } \right)} }}{{1 + \left( {t/\nu } \right)^2 }}} \left( {1 + e^{ - 2\pi t} } \right)i H_{it}^{\left( 1 \right)} \left( {it\sec \beta } \right)dt\\ & + \sin \xi \frac{{\left( { - 1} \right)^N }}{{2\left( {2\pi \cot \beta } \right)^{\frac{1}{2}} \nu ^{2N + 1} }}\int_0^{ + \infty } {\frac{{t^{2N + \frac{1}{2}} e^{ - t\left( {\tan \beta  - \beta } \right)} }}{{1 + \left( {t/\nu } \right)^2 }}} \left( {1 + e^{ - 2\pi t} } \right)i H_{it}^{\left( 1 \right)} \left( {it\sec \beta } \right)dt,
\end{split}
\end{gather}
\begin{gather}\label{eq30}
\begin{split}
R_N^{\left(Y\right)} \left( {\nu ,\beta } \right) = \; & \sin \xi \frac{{\left( { - 1} \right)^N }}{{2\left( {2\pi \cot \beta } \right)^{\frac{1}{2}} \nu ^{2N} }}\int_0^{ + \infty } {\frac{{t^{2N - \frac{1}{2}} e^{ - t\left( {\tan \beta  - \beta } \right)} }}{{1 + \left( {t/\nu } \right)^2 }}} \left( {1 + e^{ - 2\pi t} } \right) i H_{it}^{\left( 1 \right)} \left( {it\sec \beta } \right)dt\\ & - \cos \xi \frac{{\left( { - 1} \right)^N }}{{2\left( {2\pi \cot \beta } \right)^{\frac{1}{2}} \nu ^{2N + 1} }}\int_0^{ + \infty } {\frac{{t^{2N + \frac{1}{2}} e^{ - t\left( {\tan \beta  - \beta } \right)} }}{{1 + \left( {t/\nu } \right)^2 }}} \left( {1 + e^{ - 2\pi t} } \right) i H_{it}^{\left( 1 \right)} \left( {it\sec \beta } \right)dt .
\end{split}
\end{gather}
The square roots are defined to be positive on the positive real line and are defined by analytic continuation elsewhere.
\end{theorem}

Using the continuation formulas (see, e.g., \cite[p. 226]{NIST})
\begin{align*}
\sin \left( {\pi \nu } \right)H_{\nu e^{2\pi im} }^{\left( 1 \right)} \left( {\nu e^{2\pi im} \sec \beta } \right) & = \sin \left( {\pi \nu } \right)H_\nu ^{\left( 1 \right)} \left( {\nu e^{2\pi im} \sec \beta } \right)\\ & =  - \sin \left( {\left( {2m - 1} \right)\pi \nu } \right)H_\nu ^{\left( 1 \right)} \left( {\nu \sec \beta } \right) - e^{ - \pi i\nu } \sin \left( {2\pi m\nu } \right)H_\nu ^{\left( 2 \right)} \left( {\nu \sec \beta } \right),
\end{align*}
\begin{align*}
\sin \left( {\pi \nu } \right)H_{\nu e^{2\pi im} }^{\left( 2 \right)} \left( {\nu e^{2\pi im} \sec \beta } \right) & = \sin \left( {\pi \nu } \right)H_\nu ^{\left( 2 \right)} \left( {\nu e^{2\pi im} \sec \beta } \right)\\ & = \sin \left( {\left( {2m + 1} \right)\pi \nu } \right)H_\nu ^{\left( 2 \right)} \left( {\nu \sec \beta } \right) + e^{\pi i\nu } \sin \left( {2\pi m\nu } \right)H_\nu ^{\left( 1 \right)} \left( {\nu \sec \beta } \right),
\end{align*}
and the resurgence formulas \eqref{eq21}, \eqref{eq17} and \eqref{eq12}, we can derive analogous representations in sectors of the form
\[
\left( {2m - \frac{1}{2}} \right)\pi  < \arg \nu  < \left( {2m + \frac{3}{2}} \right)\pi \; \text{ and } \; \left( {2m - \frac{3}{2}} \right)\pi  < \arg \nu  < \left( {2m + \frac{1}{2}} \right)\pi, \; m \in \mathbb{Z},
\]
respectively. The lines $\arg \nu  = \left( {2m - \frac{1}{2}} \right)\pi$ are the Stokes lines for the function $H_\nu ^{\left( 1 \right)} \left( {\nu \sec \beta } \right)$, and the lines $\arg \nu  = \left( {2m + \frac{1}{2}} \right)\pi$ are the Stokes lines for the function $H_\nu ^{\left( 2 \right)} \left( {\nu \sec \beta } \right)$.

Similarly, applying the continuation formulas
\begin{align*}
& J_{\nu e^{\left( {2m + 1} \right)\pi i} } \left( {\nu e^{\left( {2m + 1} \right)\pi i} \sec \beta } \right) = J_{ - \nu } \left( {\nu e^{\left( {2m + 1} \right)\pi i} \sec \beta } \right)\\
& =  e^{2\pi im\nu } J_\nu  \left( {\nu \sec \beta } \right) - i\sin \left( {2\pi m\nu } \right)H_\nu ^{\left( 1 \right)} \left( {\nu \sec \beta } \right)  - ie^{ - \pi i\nu } \sin \left( {\left( {2m + 1} \right)\nu \pi } \right)H_\nu ^{\left( 2 \right)} \left( {\nu \sec \beta } \right),
\end{align*}
\begin{align*}
& Y_{\nu e^{\left( {2m + 1} \right)\pi i} } \left( {\nu e^{\left( {2m + 1} \right)\pi i} \sec \beta } \right) = Y_{ - \nu } \left( {\nu e^{\left( {2m + 1} \right)\pi i} \sec \beta } \right) \\
& = e^{ - 2\left( {m + 1} \right)\pi i\nu } Y_\nu  \left( {\nu \sec \beta } \right) + 2ie^{ - \pi i\nu } \sin \left( {\left( {2m + 1} \right)\pi \nu } \right)\cot \left( {\pi \nu } \right)J_\nu  \left( {\nu \sec \beta } \right) \\
 & \;\; - \sin \left( {2\pi m\nu } \right)H_\nu ^{\left( 1 \right)} \left( {\nu \sec \beta } \right) - e^{ - \pi i\nu } \sin \left( {\left( {2m + 1} \right)\pi \nu } \right)H_\nu ^{\left( 2 \right)} \left( {\nu \sec \beta } \right)
\end{align*}
and the representations \eqref{eq21}--\eqref{eq23}, we can obtain analogous formulas in any sector of the form
\[
\left( {2m + \frac{1}{2}} \right)\pi  < \arg \nu  < \left( {2m + \frac{3}{2}} \right)\pi , \; m \in \mathbb{Z}.
\]
Finally, from
\[
J_{\nu e^{2\pi im} } \left( {\nu e^{2\pi im} \sec \beta } \right) = J_\nu  \left( {\nu e^{2\pi im} \sec \beta } \right) = e^{2\pi im\nu } J_\nu  \left( {\nu \sec \beta } \right),
\]
\[
Y_{\nu e^{2\pi im} } \left( {\nu e^{2\pi im} \sec \beta } \right) = Y_\nu  \left( {\nu e^{2\pi im} \sec \beta } \right) = e^{ - 2\pi im\nu } Y_\nu  \left( {\nu \sec \beta } \right) + 2i\sin \left( {2\pi m\nu } \right)\cot \left( {\pi \nu } \right)J_\nu  \left( {\nu \sec \beta } \right)
\]
and the resurgence formulas \eqref{eq22}--\eqref{eq23}, we can derive the corresponding representations in sectors of the form
\[
\left( {2m - \frac{1}{2}} \right)\pi  < \arg \nu  < \left( {2m + \frac{1}{2}} \right)\pi , \; m \in \mathbb{Z}.
\]
The lines $\arg \nu  = \left( {2m \pm \frac{1}{2}} \right)\pi$ are the Stokes lines for the functions $J_\nu\left( {\nu \sec \beta } \right)$ and $Y_\nu \left( {\nu \sec \beta } \right)$.

When $\nu$ is an integer, the limiting values have to be taken in these continuation formulas.

The second theorem provides resurgence formulas for $H_\nu^{\left(1\right)}\left(\nu\right)$, $H_\nu^{\left(2\right)}\left(\nu\right)$, $J_\nu \left(\nu\right)$ and $Y_\nu \left(\nu\right)$. The notations follow the ones used by Wong \cite[p. 102]{Wong}.

\begin{theorem}\label{thm2} For any non-negative integer $N$, we have
\begin{equation}\label{eq44}
H_\nu ^{\left( 1 \right)} \left( \nu  \right) =  - \frac{2}{{3\pi }}\sum\limits_{n = 0}^{N - 1} {d_{2n} e^{\frac{{2\left( {2n + 1} \right)\pi i}}{3}} \sin \left( {\frac{{\left( {2n + 1} \right)\pi }}{3}} \right)\frac{{\Gamma \left( {\frac{{2n + 1}}{3}} \right)}}{{\nu ^{\frac{{2n + 1}}{3}} }}}  + R_N^{\left( H \right)} \left( \nu  \right)
\end{equation}
for $-\frac{\pi}{2} < \arg \nu < \frac{3\pi}{2}$;
\begin{equation}\label{eq18}
H_\nu ^{\left( 2 \right)} \left( \nu  \right) =  - \frac{2}{{3\pi }}\sum\limits_{n = 0}^{N - 1} {d_{2n} e^{ - \frac{{2\left( {2n + 1} \right)\pi i}}{3}} \sin \left( {\frac{{\left( {2n + 1} \right)\pi }}{3}} \right)\frac{{\Gamma \left( {\frac{{2n + 1}}{3}} \right)}}{{\nu ^{\frac{{2n + 1}}{3}} }}}  - R_N^{\left( H \right)} \left( {\nu e^{\pi i} } \right)
\end{equation}
for $-\frac{3\pi}{2} < \arg \nu < \frac{\pi}{2}$;
\begin{equation}\label{eq40}
J_\nu  \left( \nu  \right) = \frac{1}{{3\pi }}\sum\limits_{n = 0}^{N - 1} {d_{2n} \sin \left( {\frac{{\left( {2n + 1} \right)\pi }}{3}} \right)\frac{{\Gamma \left( {\frac{{2n + 1}}{3}} \right)}}{{\nu ^{\frac{{2n + 1}}{3}} }}}  + R_N^{\left( J \right)} \left( \nu  \right),
\end{equation}
\begin{equation}\label{eq41}
Y_\nu  \left( \nu  \right) =  - \frac{2}{{3\pi }}\sum\limits_{n = 0}^{N - 1} {d_{2n} \sin ^2 \left( {\frac{{\left( {2n + 1} \right)\pi }}{3}} \right)\frac{{\Gamma \left( {\frac{{2n + 1}}{3}} \right)}}{{\nu ^{\frac{{2n + 1}}{3}} }}}  + R_N^{\left( Y \right)} \left( \nu  \right)
\end{equation}
for $\left|\arg \nu\right| < \frac{\pi}{2}$. The coefficients $d_{2n}$ are given by
\begin{equation}\label{eq15}
d_{2n} = \frac{1}{{\left( {2n} \right)!}}\left[ {\frac{{d^{2n} }}{{dt^{2n} }}\left( {\frac{{t^3 }}{{\sinh t - t}}} \right)^{\frac{{2n + 1}}{3}} } \right]_{t = 0}  = \frac{{\left( { - 1} \right)^n }}{{\Gamma \left( {\frac{{2n + 1}}{3}} \right)}}\int_0^{ + \infty } {t^{\frac{{2n - 2}}{3}} e^{ - 2\pi t} i H_{it}^{\left( 1 \right)} \left( {it} \right)dt} .
\end{equation}
The remainder terms can be expressed as
\begin{equation}\label{eq14}
R_N^{\left( H \right)} \left( \nu  \right) = \frac{{\left( { - 1} \right)^N }}{{3\pi \nu ^{\frac{{2N + 1}}{3}} }}\int_0^{ + \infty } {t^{\frac{{2N - 2}}{3}} e^{ - 2\pi t} \left( {\frac{{e^{\frac{{\left( {2N + 1} \right)\pi i}}{3}} }}{{1 + \left( {t/\nu } \right)^{\frac{2}{3}} e^{\frac{{2\pi i}}{3}} }} + \frac{1}{{1 + \left( {t/\nu } \right)^{\frac{2}{3}} }}} \right)H_{it}^{\left( 1 \right)} \left( {it} \right)dt} ,
\end{equation}
\begin{equation}\label{eq42}
R_N^{\left( J \right)} \left( \nu  \right) = \frac{{\left( { - 1} \right)^N }}{{6\pi \nu ^{\frac{{2N + 1}}{3}} }}\int_0^{ + \infty } {t^{\frac{{2N - 2}}{3}} e^{ - 2\pi t} \left( {\frac{{e^{\frac{{\left( {2N + 1} \right)\pi i}}{3}} }}{{1 + \left( {t/\nu } \right)^{\frac{2}{3}} e^{\frac{{2\pi i}}{3}} }} - \frac{{e^{ - \frac{{\left( {2N + 1} \right)\pi i}}{3}} }}{{1 + \left( {t/\nu } \right)^{\frac{2}{3}} e^{ - \frac{{2\pi i}}{3}} }}} \right)H_{it}^{\left( 1 \right)} \left( {it} \right)dt} ,
\end{equation}
\begin{equation}\label{eq43}
R_N^{\left( Y \right)} \left( \nu  \right) = \frac{{\left( { - 1} \right)^N }}{{6\pi i\nu ^{\frac{{2N + 1}}{3}} }}\int_0^{ + \infty } {t^{\frac{{2N - 2}}{3}} e^{ - 2\pi t} \left( {\frac{{e^{\frac{{\left( {2N + 1} \right)\pi i}}{3}} }}{{1 + \left( {t/\nu } \right)^{\frac{2}{3}} e^{\frac{{2\pi i}}{3}} }} + \frac{{e^{ - \frac{{\left( {2N + 1} \right)\pi i}}{3}} }}{{1 + \left( {t/\nu } \right)^{\frac{2}{3}} e^{ - \frac{{2\pi i}}{3}} }} + \frac{2}{{1 + \left( {t/\nu } \right)^{\frac{2}{3}} }}} \right)H_{it}^{\left( 1 \right)} \left( {it} \right)dt} .
\end{equation}
The cube roots are defined to be positive on the positive real line and are defined by analytic continuation elsewhere.
\end{theorem}

Again, these formulas can be extended to other sectors of the complex plane like the ones in Theorem \ref{thm1}. (One has to replace the factor $\sec \beta$ by $1$ in the continuation formulas given above.)

In Section \ref{section3}, we will show how to obtain numerically computable bounds for the remainder terms using their explicit form given in Theorems \ref{thm1} and \ref{thm2}. If we neglect these remainder terms and extend the sums to $N=\infty$ in Theorems \ref{thm1} and \ref{thm2}, we recover Debye's asymptotic series. Some other formulas for the coefficients $U_n\left(i \cot \beta\right)$ and $d_{2n}$ can be found in Appendix \ref{appendixa}.

In the following two theorems, we give exponentially improved asymptotic expansions for the function $H_\nu ^{\left( 1 \right)} \left( \nu x \right)$ when $x>1$ and $x=1$, respectively. The related expansions for the functions $H_\nu^{\left(2\right)}\left(\nu x\right)$, $J_\nu \left(\nu x\right)$ and $Y_\nu \left(\nu x\right)$ may be derived from the corresponding connection formulas. The resulting expansions for $J_\nu \left(\nu x\right)$ and $Y_\nu \left(\nu x\right)$ can be viewed as the mathematically rigorous forms of the terminated expansions of Dingle \cite[pp. 469--472]{Dingle}. In these theorems we truncate the asymptotic series of $H_\nu ^{\left( 1 \right)} \left( \nu x \right)$ at about the least term and re-expand the remainders into new asymptotic expansions. The resulting exponentially improved asymptotic series are valid in larger regions than the original Debye expansions. The terms in these new series involve the Terminant function $\widehat T_p\left(z\right)$, which allows the smooth transition through the Stokes line $\arg \nu = -\frac{\pi}{2}$. For the definition and basic properties of the Terminant function, see Section \ref{section5}. It is interesting to note that for the case $x>1$, an additional sum is needed for the re-expansion beside the truncated Debye series. The appearance of this second sum is due to the presence of the exponential quantity $e^{-2\pi t}$ in \eqref{eq12}.

For exponentially improved expansions using Hadamard series, see Paris \cite{Paris3}, \cite[pp. 198--206]{Paris}.

\begin{theorem}\label{thm3} Let $0 < \beta < \frac{\pi}{2}$ be a fixed acute angle. Define $R_{N,M}^{\left( H \right)} \left( {\nu ,\beta } \right)$ by
\[
H_\nu ^{\left( 1 \right)} \left( {\nu \sec \beta } \right) = \frac{{e^{i\nu \left( {\tan \beta  - \beta } \right) - \frac{\pi }{4}i} }}{{\left( {\frac{1}{2}\nu \pi \tan \beta } \right)^{\frac{1}{2}} }}\left(\sum\limits_{n = 0}^{N - 1} {\left( { - 1} \right)^n \frac{{U_n \left( {i\cot \beta } \right)}}{{\nu ^n }}}  + \sum\limits_{m = N}^{M - 1} {\left( { - 1} \right)^m \frac{{\widetilde U_m \left( {i\cot \beta } \right)}}{{\nu ^m }}}  + R_{N,M}^{\left( H \right)} \left( {\nu ,\beta } \right)\right)
\]
where
\[
\widetilde U_m \left( {i\cot \beta } \right) = \frac{{i^m }}{{2\left( {2\pi \cot \beta } \right)^{\frac{1}{2}} }}\int_0^{ + \infty } {t^{m - \frac{1}{2}} e^{ - t\left( {\tan \beta  - \beta  + 2\pi } \right)} iH_{it}^{\left( 1 \right)} \left( {it\sec \beta } \right)dt} ,
\]
\[
N = 2\left| \nu  \right|\left( {\tan \beta  - \beta } \right) + \rho \; \text{ and } \; M = 2\left| \nu  \right|\left( {\tan \beta  - \beta  + \pi } \right) + \sigma,
\]
$\left|\nu\right|$ being large, $\rho$ and $\sigma$ being bounded quantities such that $1 \leq N \leq M$. Then
\begin{align*}
R_{N,M}^{\left( H \right)} \left( {\nu ,\beta } \right) = \; & ie^{ - 2i\nu\left( {\tan \beta  - \beta } \right) } \sum\limits_{k = 0}^{K - 1} {\frac{{U_k \left( {i\cot \beta } \right)}}{{\nu ^k }}\widehat T_{N - k} \left( { - 2i\nu\left( {\tan \beta  - \beta } \right) } \right)} \\ & + ie^{ - 2i \nu \left( {\tan \beta  - \beta  + \pi } \right) } \sum\limits_{\ell  = 0}^{L - 1} {\frac{{U_\ell  \left( {i\cot \beta } \right)}}{{\nu ^\ell  }}\widehat T_{M - \ell } \left( { - 2i\nu\left( {\tan \beta  - \beta  + \pi } \right)} \right)}  + R_{N,M,K,L}^{\left( H \right)} \left( {\nu ,\beta } \right),
\end{align*}
where $K$ and $L$ are arbitrary fixed non-negative integers, and
\[
R_{N,M,K,L}^{\left( H \right)} \left( {\nu ,\beta } \right) = \mathcal{O}_{K,\rho } \left( {e^{ - 2\left| \nu  \right|\left( {\tan \beta  - \beta } \right)} \frac{{\left|U_K \left( {i\cot \beta } \right)\right|}}{{\left| \nu  \right|^K }}} \right) + \mathcal{O}_{L,\sigma } \left( {e^{ - 2\left| \nu  \right|\left( {\tan \beta  - \beta  + \pi } \right)} \frac{{\left|U_L \left( {i\cot \beta } \right)\right|}}{{\left| \nu  \right|^L }}} \right)
\]
for $-\frac{\pi}{2} \le \arg \nu \le \frac{3\pi}{2}$,
\[
R_{N,M,K,L}^{\left( H \right)} \left( {\nu ,\beta } \right) = \mathcal{O}_{K,\rho} \left( {e^{2\Im \left( \nu  \right)\left( {\tan \beta  - \beta } \right)} \frac{{\left| {U_K \left( {i\cot \beta } \right)} \right|}}{{\left| \nu  \right|^K }}} \right) + \mathcal{O}_{L,\sigma} \left( {e^{2\Im \left( \nu  \right)\left( {\tan \beta  - \beta  + \pi } \right)} \frac{{\left| {U_L \left( {i\cot \beta } \right)} \right|}}{{\left| \nu  \right|^L }}} \right)
\]
for $-\frac{3\pi}{2} \leq \arg \nu  \leq  - \frac{\pi }{2}$.
\end{theorem}

\begin{theorem}\label{thm4} Define $R_{N,M}^{\left( H \right)} \left( \nu \right)$ by
\[
H_\nu ^{\left( 1 \right)} \left( \nu  \right) = \frac{{e^{ - \frac{\pi}{3}i} }}{{\sqrt 3 \pi \nu ^{\frac{1}{3}} }}\sum\limits_{n = 0}^{N - 1} {d_{6n} \frac{{\Gamma \left( {2n+\frac{1}{3}} \right)}}{{\nu ^{2n} }}}  - \frac{{e^{\frac{\pi}{3}i} }}{{\sqrt 3 \pi \nu ^{\frac{5}{3}} }}\sum\limits_{m = 0}^{M - 1} {d_{6m + 4} \frac{{\Gamma \left( {2m+\frac{5}{3}} \right)}}{{\nu ^{2m} }}}  + R_{N,M}^{\left( H \right)} \left( \nu  \right),
\]
where
\[
N = \pi \left| \nu \right| + \rho \; \text{ and } \; M = \pi \left| \nu \right| + \sigma,
\]
$\left|\nu\right|$ being large, $\rho$ and $\sigma$ being bounded quantities such that $M,N \geq 1$. Then
\begin{align*}
R_{N,M}^{\left( H \right)} \left( \nu  \right) = \; & ie^{ - \frac{\pi }{3}i} \frac{{e^{ - 2\pi i\nu } }}{{\sqrt 3 }}\frac{2}{{3\pi }}\sum\limits_{k = 0}^{K - 1} {d_{2k} \sin \left( {\frac{{\left( {2k + 1} \right)\pi }}{3}} \right)\frac{{\Gamma \left( {\frac{{2k + 1}}{3}} \right)}}{{\nu ^{\frac{{2k + 1}}{3}} }}\widehat T_{2N - \frac{{2k}}{3}} \left( { - 2\pi i\nu } \right)} \\
& - i\frac{{e^{2\pi i\nu } }}{{\sqrt 3 }}\frac{2}{{3\pi }}\sum\limits_{k = 0}^{K - 1} {d_{2k} e^{\frac{{2\left( {2k + 1} \right)\pi i}}{3}} \sin \left( {\frac{{\left( {2k + 1} \right)\pi }}{3}} \right)\frac{{\Gamma \left( {\frac{{2k + 1}}{3}} \right)}}{{\nu ^{\frac{{2k + 1}}{3}} }}\widehat T_{2N - \frac{{2k}}{3}} \left( {2\pi i\nu } \right)} \\
& - ie^{\frac{\pi }{3}i} \frac{{e^{ - 2\pi i\nu } }}{{\sqrt 3 }}\frac{2}{{3\pi }}\sum\limits_{\ell  = 0}^{L - 1} {d_{2\ell } \sin \left( {\frac{{\left( {2\ell  + 1} \right)\pi }}{3}} \right)\frac{{\Gamma \left( {\frac{{2\ell  + 1}}{3}} \right)}}{{\nu ^{\frac{{2\ell  + 1}}{3}} }}\widehat T_{2M - \frac{{2\ell  - 4}}{3}} \left( { - 2\pi i\nu } \right)} \\
& + i\frac{{e^{2\pi i\nu } }}{{\sqrt 3 }}\frac{2}{{3\pi }}\sum\limits_{\ell  = 0}^{L - 1} {d_{2\ell } e^{\frac{{2\left( {2\ell  + 1} \right)\pi i}}{3}} \sin \left( {\frac{{\left( {2\ell  + 1} \right)\pi }}{3}} \right)\frac{{\Gamma \left( {\frac{{2\ell  + 1}}{3}} \right)}}{{\nu ^{\frac{{2\ell  + 1}}{3}} }} \widehat T_{2M - \frac{{2\ell  - 4}}{3}} \left( {2\pi i\nu } \right)} \\
& + R_{N,M,K,L}^{\left( H \right)} \left( \nu  \right) ,
\end{align*}
where $K$ and $L$ are arbitrary fixed non-negative integers satisfying $K,L \equiv 0 \mod 3$, and
\begin{equation}\label{eq80}
R_{N,M,K,L}^{\left( H \right)} \left( \nu  \right) = \mathcal{O}_{K,\rho } \left( {e^{ - 2\pi \left| \nu  \right|} \left| {d_{2K} } \right|\frac{{\Gamma \left( {\frac{{2K + 1}}{3}} \right)}}{{\left| \nu  \right|^{\frac{{2K + 1}}{3}} }}} \right) + \mathcal{O}_{L,\sigma } \left( {e^{ - 2\pi \left| \nu  \right|} \left| {d_{2L} } \right|\frac{{\Gamma \left( {\frac{{2L + 1}}{3}} \right)}}{{\left| \nu  \right|^{\frac{{2L + 1}}{3}} }}} \right)
\end{equation}
for $- \frac{\pi }{2} \le \arg \nu \le \frac{\pi }{2}$,
\[
R_{N,M,K,L}^{\left( H \right)} \left( \nu  \right) = \mathcal{O}_{K,\rho} \left( {e^{ \mp 2\pi \Im \left( \nu  \right)} \left| {d_{2K} } \right|\frac{{\Gamma \left( {\frac{{2K + 1}}{3}} \right)}}{{\left| \nu  \right|^{\frac{{2K + 1}}{3}} }}} \right) + \mathcal{O}_{L,\sigma} \left( {e^{ \mp 2\pi \Im \left( \nu  \right)} \left| {d_{2L} } \right|\frac{{\Gamma \left( {\frac{{2L + 1}}{3}} \right)}}{{\left| \nu  \right|^{\frac{{2L + 1}}{3}} }}} \right)
\]
for $\frac{\pi}{2} \le \pm \arg \nu  \le \frac{3\pi}{2}$. Moreover, if $K=L$ then the bound \eqref{eq80} remains valid in the larger sector $- \frac{\pi }{2} \le \arg \nu \le \frac{3\pi}{2}$.
\end{theorem}

The assumption that $K,L \equiv 0 \mod 3$ is only for simplicity. Estimations for $R_{N,M,K,L}^{\left( H \right)} \left(\nu\right)$ when $K$ or $L$ may not be divisible by $3$ can be obtained similarly.

The rest of the paper is organised as follows. In Section \ref{section2}, we prove the resurgence formulas stated in Theorems \ref{thm1} and \ref{thm2}. In Section \ref{section3}, we give explicit and numerically computable error bounds for Debye's expansions when $x\geq 1$ using the results of Section \ref{section2}. In Section \ref{section4}, asymptotic approximations for the late terms in Debye's expansions are given. In Section \ref{section5}, we prove the exponentially improved expansions presented in Theorems \ref{thm3} and \ref{thm4}, and provide a detailed discussion of the Stokes phenomenon related to the expansions of $H_\nu ^{\left( 1 \right)} \left( \nu x \right)$. The paper concludes with a discussion in Section \ref{section6}.

\section{Proofs of the resurgence formulas}\label{section2}

In this section we prove the resurgence formulas given in Theorems \ref{thm1} and \ref{thm2}. In the first subsection, we prove the results related to the Hankel functions. In the second subsection we show how the corresponding formulas for the Bessel functions can be derived using the results for the Hankel functions.

\subsection{The Hankel functions $H_\nu^{\left(1\right)}\left(\nu x\right)$ and $H_\nu^{\left(2\right)}\left(\nu x\right)$} Our analysis is based on the Schl\"afli--Sommerfeld integral representation
\[
H_\nu ^{\left( 1 \right)} \left( z \right) = \frac{1}{{\pi i}}\int_{ - \infty }^{\infty  + \pi i} {e^{z\sinh t - \nu t} dt} \quad \left|\arg z\right| <\frac{\pi}{2}
\]
and the connection formula $H_\nu ^{\left( 2 \right)} \left( z \right) = \overline {H_{ \bar{\nu} }^{\left( 1 \right)} \left( {\bar z} \right)}$ \cite[pp. 224 and 226]{NIST}. If $z = \nu x$, where $x$ is a positive constant, then
\begin{equation}\label{eq3}
H_\nu ^{\left( 1 \right)} \left( {\nu x} \right) = \frac{1}{{\pi i}}\int_{ - \infty }^{\infty  + \pi i} {e^{\nu \left( {x\sinh t - t} \right)} dt} \quad \left|\arg \nu \right| <\frac{\pi}{2}.
\end{equation}
The analysis is significantly different according to whether $x>1$ or $x=1$. The saddle points of the integrand are the roots of the equation $x\cosh t = 1$. Hence, the saddle points are given by $t_{\pm}^{\left( k \right)} = \pm \mathrm{sech}^{ - 1} x + 2\pi ik$ where $k$ is an arbitrary integer. When $x=1$, we shall use the simpler notation $t^{\left( k \right)} = 2\pi ik$. We denote by $\mathscr{C}_{\pm}^{\left( k \right)}\left(\theta\right)$ the portion of the steepest descent paths that pass through the saddle point $t_{\pm}^{\left( k \right)}$. Here, and subsequently, we write $\theta = \arg \nu$. Similarly, $\mathscr{C}^{\left( k \right)}\left(\theta\right)$ denotes the steepest descent paths through the saddle point $t^{\left( k \right)}$.

\begin{figure}[!t]
\def\svgwidth{0.8\columnwidth}
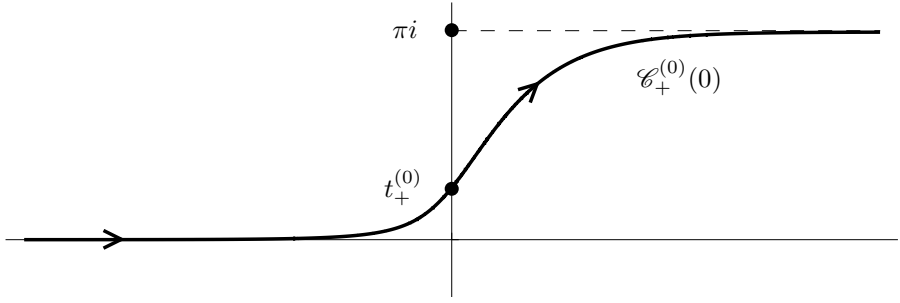
\caption{The steepest descent path for $H_\nu ^{\left( 1 \right)} \left( {\nu \sec \beta } \right)$ through the saddle point $t_+^{\left( 0 \right)} = i \beta$ when $\arg \nu =0$.}
\label{fig1}
\end{figure}

\subsubsection{Case (i): $x>1$} 
\begin{figure}[!t]
\def\svgwidth{0.8\columnwidth}
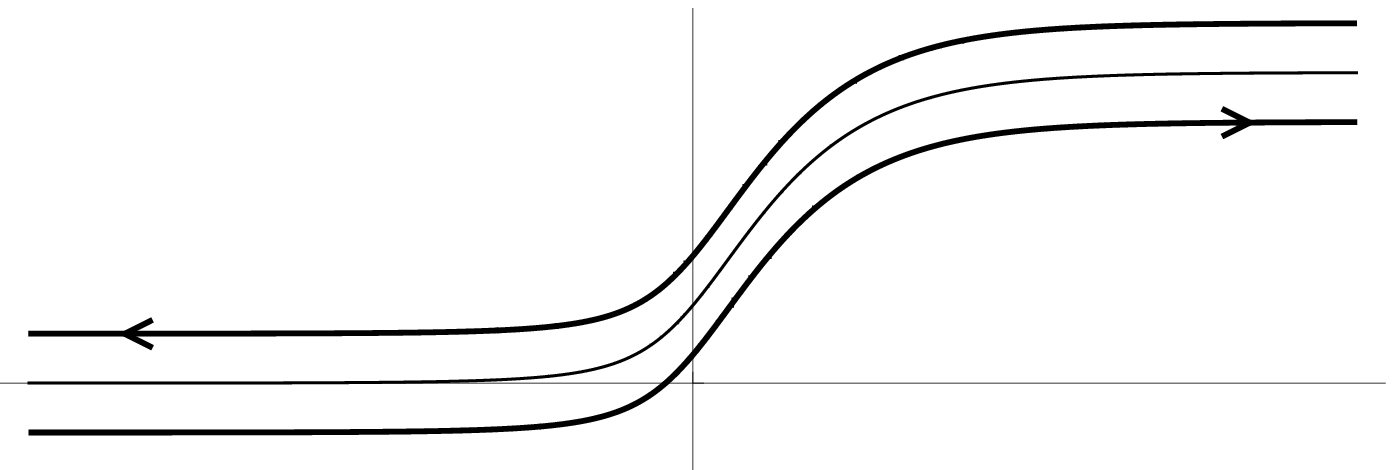
\caption{The contour $\Gamma_+^{\left( 0 \right)}$ encircling the path $\mathscr{C}_+^{\left( 0 \right)}\left(0\right)$.}
\label{fig2}
\end{figure}

Let $0<\beta<\frac{\pi}{2}$ be defined by $\sec \beta = x$. The saddle points in $\left|\Im \left(t\right)\right| < \pi$ are located at $t =t_\pm^{\left(0\right)}= \pm i\beta$. For the integral \eqref{eq3}, we need to consider only the saddle point $t_+^{\left( 0 \right)} = i \beta$. (The saddle point $t_-^{\left( 0 \right)} = -i \beta$ is suitable for the corresponding integral representation for $H_\nu ^{\left( 2 \right)} \left( {\nu \sec \beta } \right)$.) Let $f\left(t,\beta\right) = \sec \beta \sinh t -t$. The steepest descent path in this case is given by
\[
\arg \left[ {e^{i\theta } \left( {f\left( {t_ + ^{\left( 0 \right)} ,\beta } \right) - f\left( {t,\beta } \right)} \right)} \right] = 0 .
\]
The path $\mathscr{C}_+^{\left( 0 \right)}\left(\theta\right)$ for $\theta = 0$, with an appropriate orientation, is shown in Figure \ref{fig1}. For simplicity, we assume that $\theta = 0$. In due course, we shall appeal to an analytic continuation argument to extend our results to complex $\nu$. By Cauchy's theorem, we can deform the path of integration in \eqref{eq3} to $\mathscr{C}_+^{\left( 0 \right)}\left(0\right)$. If
\begin{equation}\label{eq2}
\tau  =  f\left( {t_+^{\left( 0 \right)} ,\beta } \right) - f\left( {t,\beta } \right) = f\left( {i \beta ,\beta } \right) - f\left( {t,\beta } \right) =  i\left(\tan\beta -\beta\right)  - \left(\sec \beta \sinh t - t\right),
\end{equation}
then $\tau$ is real on the curve $\mathscr{C}_+^{\left( 0 \right)}\left(0\right)$, and, as $t$ travels along this curve from $-\infty$ to $\infty + \pi i$, $\tau$ decreases from $+\infty$ to $0$ and then increases to $+ \infty$. Therefore, corresponding to each positive value of $\tau$, there are two values of $t$, say $t_1$ and $t_2$, satisfying \eqref{eq2} with $\Re \left(t_1\right) > 0$ and $\Re \left(t_2\right) < 0$. In terms of $\tau$, we have
\begin{align*}
H_\nu ^{\left( 1 \right)} \left( {\nu \sec \beta } \right) & = \frac{{e^{i\nu \left( {\tan \beta  - \beta } \right)} }}{{\pi i}}\int_0^{ + \infty } {e^{ - \nu \tau } \left( {\frac{{dt_1 }}{{d\tau }} - \frac{{dt_2 }}{{d\tau }}} \right)d\tau } \\
& = \frac{{e^{i\nu \left( {\tan \beta  - \beta } \right)} }}{{\pi i}}\int_0^{ + \infty } {e^{ - \nu \tau } \left( {\frac{1}{{1 - \sec \beta \cosh t_1 \left( \tau  \right)}} - \frac{1}{{1 - \sec \beta \cosh t_2 \left( \tau  \right)}}} \right)d\tau } .
\end{align*}
Following Berry and Howls, we express the quantity in the large parentheses as a contour integral using the residue theorem, to find
\[
H_\nu ^{\left( 1 \right)} \left( {\nu \sec \beta } \right) = \frac{{e^{i\nu \left( {\tan \beta  - \beta } \right)} }}{{\pi i}}\int_0^{ + \infty } {\tau ^{ - \frac{1}{2}} e^{ - \nu \tau } \frac{1}{{2\pi i}}\oint_{\Gamma_+^{\left( 0 \right)}}  {\frac{{\left( {f\left( {i\beta ,\beta } \right) - f\left( {u,\beta } \right)} \right)^{ - \frac{1}{2}} }}{{1 - \tau \left( {f\left( {i\beta ,\beta } \right) - f\left( {u,\beta } \right)} \right)^{ - 1} }}du} d\tau } ,
\]
where the contour $\Gamma_+^{\left( 0 \right)}$ encircles the path $\mathscr{C}_+^{\left( 0 \right)}\left(0\right)$ in the positive direction and does not enclose any of the saddle points $t_ \pm ^{\left( k \right)} \neq t_ + ^{\left( 0 \right)}$ (see Figure \ref{fig2}). The square root is defined so that $\left( {f\left( {i\beta ,\beta } \right) - f\left( {t,\beta } \right)} \right)^{\frac{1}{2}}$ is positive on the portion of $\mathscr{C}_+^{\left( 0 \right)}\left(0\right)$ that starts at $t= t_+^{\left( 0 \right)}$ and ends at $t=\infty + \pi i$. Now, we employ the well-known expression for non-negative integer $N$
\begin{equation}\label{eq7}
\frac{1}{1 - z} = \sum\limits_{n = 0}^{N-1} {z^n}  + \frac{z^N}{1 - z},\; z \neq 1,
\end{equation}
to expand the function under the contour integral in powers of ${\tau \left( {f\left( {i\beta ,\beta } \right) - f\left( {u,\beta } \right)} \right)^{ - 1} }$. The result is
\begin{align*}
H_\nu ^{\left( 1 \right)} \left( {\nu \sec \beta } \right) = \; & \frac{{e^{i\nu \left( {\tan \beta  - \beta } \right)} }}{{\pi i}}\sum\limits_{n = 0}^{N - 1} {\int_0^{ + \infty } {\tau ^{n - \frac{1}{2}} e^{ - \nu \tau } \frac{1}{{2\pi i}}\oint_{\Gamma_+^{\left( 0 \right)}}  {\frac{{du}}{{\left( {f\left( {i\beta ,\beta } \right) - f\left( {u,\beta } \right)} \right)^{n + \frac{1}{2}} }}} d\tau } } \\ & + \frac{{e^{i\nu \left( {\tan \beta  - \beta } \right) - \frac{\pi}{4}i} }}{{\left( {\frac{1}{2}\nu \pi \tan \beta } \right)^{\frac{1}{2}} }}R_N^{\left(H\right)} \left( {\nu ,\beta } \right),
\end{align*}
where
\begin{equation}\label{eq1}
R_N^{\left(H\right)} \left( {\nu ,\beta } \right) = \frac{{e^{\frac{\pi }{4}i} \nu ^{\frac{1}{2}} }}{{i\left( {2\pi \cot \beta } \right)^{\frac{1}{2}} }}\int_0^{ + \infty } {\tau ^{N - \frac{1}{2}} e^{ - \nu \tau } \frac{1}{{2\pi i}}\oint_{\Gamma _ + ^{\left( 0 \right)} } {\frac{{\left( {f\left( {i\beta ,\beta } \right) - f\left( {u,\beta } \right)} \right)^{ - N - \frac{1}{2}} }}{{1 - \tau \left( {f\left( {i\beta ,\beta } \right) - f\left( {u,\beta } \right)} \right)^{ - 1} }}du} d\tau } .
\end{equation}
The path $\Gamma_+^{\left( 0 \right)}$ in the sum can be shrunk into a small circle around $i \beta$, and we arrive at
\begin{equation}\label{eq10}
H_\nu ^{\left( 1 \right)} \left( {\nu \sec \beta } \right) = \frac{{e^{i\nu \left( {\tan \beta  - \beta } \right) - \frac{\pi}{4}i} }}{{\left( {\frac{1}{2}\nu \pi \tan \beta } \right)^{\frac{1}{2}} }}\left( {\sum\limits_{n = 0}^{N - 1} {\left( { - 1} \right)^n \frac{{U_n \left( {i\cot \beta } \right)}}{{\nu ^n }}}  + R_N^{\left(H\right)} \left( {\nu ,\beta } \right)} \right),
\end{equation}
where
\begin{align*}
U_n \left( {i\cot \beta } \right) & = \frac{{\left( { - 1} \right)^n e^{\frac{\pi }{4}i} }}{{i\left( {2\pi \cot \beta } \right)^{\frac{1}{2}} }} \frac{{\Gamma \left( {n + \frac{1}{2}} \right)}}{{2\pi i}}\oint_{\left( {i\beta ^ +  } \right)} {\frac{{du}}{{\left( {f\left( {i\beta ,\beta } \right) - f\left( {u,\beta } \right)} \right)^{n + \frac{1}{2}} }}} \\
& = \frac{{\left( { - 1} \right)^n i^n }}{{\left( {2\pi \cot \beta } \right)^{\frac{1}{2}} }}\frac{{\Gamma \left( {n + \frac{1}{2}} \right)}}{{2\pi i}}\oint_{\left( {0^ +  } \right)} {\left( {\frac{{t^2 }}{{i\left( {f\left( {i\beta ,\beta } \right) - f\left( {t + i\beta ,\beta } \right)} \right)}}} \right)^{n + \frac{1}{2}} \frac{{dt}}{{t^{2n + 1} }}}\\
& = \left( { - 1} \right)^n \frac{{\left( {i\cot \beta } \right)^n }}{{2^n n!}}\left[ {\frac{{d^{2n} }}{{dt^{2n} }}\left( {\frac{1}{2}\frac{{t^2 }}{{i\cot \beta \left( {t - \sinh t} \right) + \cosh t - 1}}} \right)^{n + \frac{1}{2}} } \right]_{t = 0} .
\end{align*}

\begin{figure}[!t]
\def\svgwidth{0.66\columnwidth}
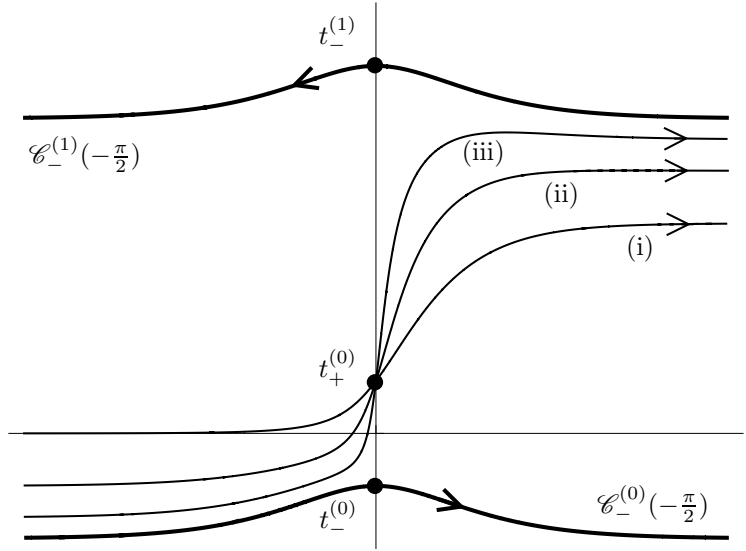
\caption{The steepest descent path for $H_\nu ^{\left( 1 \right)} \left( {\nu \sec \beta } \right)$ through the saddle point $t_+^{\left( 0 \right)}$ when (i) $\arg \nu =0$, (ii) $\arg \nu =-\frac{\pi}{4}$, (iii) $\arg \nu =-\frac{2\pi}{5}$. The paths $\mathscr{C}_ - ^{\left( 1 \right)} \left( { - \frac{\pi }{2}} \right)$ and $\mathscr{C}_ - ^{\left( 0 \right)} \left( { - \frac{\pi }{2}} \right)$ are the adjacent contours to $t_+^{\left( 0 \right)}$. The domain $\Delta^{\left(0\right)}$ comprises all points between these two paths.}
\label{fig3}
\end{figure}

Performing the change of variable $\nu \tau = s$ in \eqref{eq1} yields
\begin{equation}\label{eq4}
R_N^{\left(H\right)} \left( {\nu ,\beta } \right) = \frac{{e^{\frac{\pi }{4}i} }}{{i\left( {2\pi \cot \beta } \right)^{\frac{1}{2}} \nu ^N }}\int_0^{ + \infty } {s^{N - \frac{1}{2}} e^{ - s} \frac{1}{{2\pi i}}\oint_{\Gamma _ + ^{\left( 0 \right)} } {\frac{{\left( {f\left( {i\beta ,\beta } \right) - f\left( {u,\beta } \right)} \right)^{ - N - \frac{1}{2}} }}{{1 - \left( {s/\nu } \right)\left( {f\left( {i\beta ,\beta } \right) - f\left( {u,\beta } \right)} \right)^{ - 1} }}du} ds} .
\end{equation}
This representation of $R_N^{\left(H\right)} \left( {\nu ,\beta } \right)$ and the formula \eqref{eq10} can be continued analytically if we choose $\Gamma _ + ^{\left( 0 \right)} = \Gamma _ + ^{\left( 0 \right)}\left(\theta\right)$ to be an infinite contour that surrounds the steepest descent path $\mathscr{C}_{+}^{\left(0\right)}\left(\theta\right)$ in the anti-clockwise direction and that does not encircle any of the saddle points $t_ \pm ^{\left( k \right)} \neq t_ + ^{\left( 0 \right)}$. This continuation argument works until the path $\mathscr{C}_{+}^{\left(0\right)}\left(\theta\right)$ runs into an other saddle point. In the terminology of Berry and Howls, such saddle points are called adjacent to $t_+^{\left(0\right)}$. As
\[
\arg \left( {f\left( {t_ + ^{\left( 0 \right)} ,\beta } \right) - f\left( {t_ \pm ^{\left( k \right)} ,\beta } \right)} \right) = \frac{\pi }{2}
\]
for any saddle point $t_ \pm ^{\left( k \right)} \neq t_ + ^{\left( 0 \right)}$, we infer that \eqref{eq4} is valid as long as $-\frac{\pi}{2} < \theta < \frac{3\pi}{2}$ with a contour $\Gamma _ + ^{\left( 0 \right)}\left(\theta\right)$ specified above. When $\theta = -\frac{\pi}{2}$ or $\frac{3\pi}{2}$, the path $\mathscr{C}_{+}^{\left(0\right)}\left(\theta\right)$ connects to the saddle points $t_ - ^{\left( 0 \right)} = -i\beta$ and $t_ - ^{\left( 1 \right)} = -i\beta +2 \pi i$. These are the adjacent saddles. The set
\[
\Delta ^{\left( 0 \right)}  = \left\{ {u \in \mathscr{C}_ + ^{\left( 0 \right)} \left( \theta  \right) : - \frac{\pi }{2} < \theta  < \frac{{3\pi }}{2}} \right\}
\]
forms a domain in the complex plane whose boundaries are themselves steepest descent paths through the adjacent saddles (see Figure \ref{fig3}). These paths are $\mathscr{C}_ - ^{\left( 1 \right)} \left( { - \frac{\pi }{2}} \right)$ and $\mathscr{C}_ - ^{\left( 0 \right)} \left( { - \frac{\pi }{2}} \right)$, and they are called the adjacent contours to $t_ + ^{\left( 0 \right)}$. The function under the contour integral in \eqref{eq4} is an analytic function of $u$ in the domain $\Delta ^{\left( 0 \right)}$, therefore we can deform $\Gamma _ + ^{\left( 0 \right)}$ over the adjacent contours. We thus find that for $-\frac{\pi}{2} < \theta < \frac{3\pi}{2}$ and $N \geq 0$, \eqref{eq4} may be written
\begin{gather}\label{eq5}
\begin{split}
R_N^{\left(H\right)} \left( {\nu ,\beta } \right) = \; & \frac{e^{\frac{\pi}{4}i}}{i\left( {2\pi \cot \beta } \right)^{\frac{1}{2}} \nu ^N }\int_0^{ + \infty } {s^{N - \frac{1}{2}} e^{ - s} \frac{1}{2\pi i}\int_{\mathscr{C}_ - ^{\left( 1 \right)} \left( { - \frac{\pi }{2}} \right)} {\frac{{\left( {f\left( {i\beta ,\beta } \right) - f\left( {u,\beta } \right)} \right)^{ - N - \frac{1}{2}} }}{{1 - \left( {s/\nu } \right)\left( {f\left( {i\beta ,\beta } \right) - f\left( {u,\beta } \right)} \right)^{ - 1} }}du} ds} \\
& + \frac{e^{\frac{\pi}{4}i}}{i\left( {2\pi \cot \beta } \right)^{\frac{1}{2}} \nu ^N}\int_0^{ + \infty } {s^{N - \frac{1}{2}} e^{ - s} \frac{1}{2\pi i}\int_{\mathscr{C}_ - ^{\left( 0 \right)} \left( { - \frac{\pi }{2}} \right)} {\frac{{\left( {f\left( {i\beta ,\beta } \right) - f\left( {u,\beta } \right)} \right)^{ - N - \frac{1}{2}} }}{{1 - \left( {s/\nu } \right)\left( {f\left( {i\beta ,\beta } \right) - f\left( {u,\beta } \right)} \right)^{ - 1} }}du} ds} .
\end{split}
\end{gather}
Now we make the changes of variable
\[
s = t\frac{{\left| {f\left( { - i\beta  + 2\pi i,\beta } \right) - f\left( {i\beta ,\beta } \right)} \right|}}{{f\left( { - i\beta  + 2\pi i,\beta } \right) - f\left( {i\beta ,\beta } \right)}}\left( {f\left( {u,\beta } \right) - f\left( {i\beta ,\beta } \right)} \right) = it\left( {f\left( {u,\beta } \right) - f\left( {i\beta ,\beta } \right)} \right)
\]
in the first, and
\[
s = t\frac{{\left| {f\left( { - i\beta ,\beta } \right) - f\left( {i\beta ,\beta } \right)} \right|}}{{f\left( { - i\beta ,\beta } \right) - f\left( {i\beta ,\beta } \right)}}\left( {f\left( {u,\beta } \right) - f\left( {i\beta ,\beta } \right)} \right) = it\left( {f\left( {u,\beta } \right) - f\left( {i\beta ,\beta } \right)} \right)
\]
in the second double integral. Clearly, by the definition of the adjacent contours, $t$ is positive. The quantities $f\left( { - i\beta  + 2\pi i,\beta } \right) - f\left( {i\beta ,\beta } \right)= -2i\left( {\tan \beta  - \beta  + \pi } \right)$ and $f\left( { - i\beta ,\beta } \right) - f\left( {i\beta ,\beta } \right) = -2i\left( {\tan \beta  - \beta } \right)$ were essentially called the ``singulants" by Dingle \cite[p. 147]{Dingle}. When using these changes of variable, we should take $\left( { - i} \right)^{\frac{1}{2}}  = e^{ - \frac{\pi }{4}i}$ in the first, and $\left( { - i} \right)^{\frac{1}{2}}  = e^{\frac{3\pi }{4}i}$ in the second double integral. With these changes of variable, the representation \eqref{eq5} for $R_N^{\left(H\right)} \left( {\nu ,\beta } \right)$ becomes
\begin{gather}\label{eq11}
\begin{split}
R_N^{\left(H\right)} \left( {\nu ,\beta } \right) = \; & \frac{-i}{\left( {2\pi \cot \beta } \right)^{\frac{1}{2}} \left(i\nu \right)^N }\int_0^{ + \infty } {\frac{{t^{N - \frac{1}{2}} e^{ - t\left( {\tan \beta  - \beta } \right)} }}{{1 + it/\nu }}\frac{1}{2\pi i}\int_{\mathscr{C}_ - ^{\left( 1 \right)} \left( { - \frac{\pi }{2}} \right)} {e^{ - itf\left( {u,\beta } \right)} du} dt} \\
& +\frac{-i}{\left( {2\pi \cot \beta } \right)^{\frac{1}{2}} \left(i\nu \right)^N}\int_0^{ + \infty } {\frac{{t^{N - \frac{1}{2}} e^{ - t\left( {\tan \beta  - \beta } \right)} }}{{1 + it/\nu }}\frac{- 1}{2\pi i}\int_{\mathscr{C}_ - ^{\left( 0 \right)} \left( { - \frac{\pi }{2}} \right)} {e^{ - itf\left( {u,\beta } \right)} du} dt},
\end{split}
\end{gather}
for $-\frac{\pi}{2} < \theta < \frac{3\pi}{2}$ and $N \geq 0$. Finally, the contour integrals can themselves be represented in terms of the Hankel functions since
\begin{align*}
\frac{1}{{\pi i}}\int_{\mathscr{C}_ - ^{\left( 1 \right)} \left( { - \frac{\pi }{2}} \right)} {e^{ - itf\left( {u,\beta } \right)} du} & = \frac{-1}{{\pi i}}\int_{\mathscr{C}_ - ^{\left( 0 \right)} \left( { - \frac{\pi }{2}} \right)} {e^{ - itf\left( {u + 2\pi i,\beta } \right)} du}  =  - \frac{{e^{ - 2\pi t} }}{{\pi i}}\int_{\mathscr{C}_ - ^{\left( 0 \right)} \left( { - \frac{\pi }{2}} \right)} {e^{ - itf\left( {u,\beta } \right)} du}\\ & = e^{ - 2\pi t} H_{ - it}^{\left( 2 \right)} \left( { - it\sec \beta } \right) =  - e^{ - 2\pi t} H_{it}^{\left( 1 \right)} \left( {it\sec \beta } \right),
\end{align*}
and
\[
\frac{{ - 1}}{{\pi i}}\int_{\mathscr{C}_ - ^{\left( 0 \right)} \left( { - \frac{\pi }{2}} \right)} {e^{ - itf\left( {u,\beta } \right)} du}  = H_{ - it}^{\left( 2 \right)} \left( { - it\sec \beta } \right) =  - H_{it}^{\left( 1 \right)} \left( {it\sec \beta } \right).
\]
Substituting these into \eqref{eq11} gives \eqref{eq12}. Formula \eqref{eq17} follows from the connection formula $H_\nu ^{\left( 2 \right)} \left( \nu \sec \beta \right) = \overline {H_{ \bar{\nu} }^{\left( 1 \right)} \left( {\bar \nu \sec \beta} \right)}$. To prove the second representation in \eqref{eq13}, we apply \eqref{eq12} for the right-hand side of $U_n \left( {i\cot \beta } \right) = \left( { - 1} \right)^n \nu ^n \left( {R_n^{\left( H \right)} \left( {\nu ,\beta } \right) - R_{n + 1}^{\left( H \right)} \left( {\nu ,\beta } \right)} \right)$. 

\begin{figure}[!t]
\def\svgwidth{0.8\columnwidth}
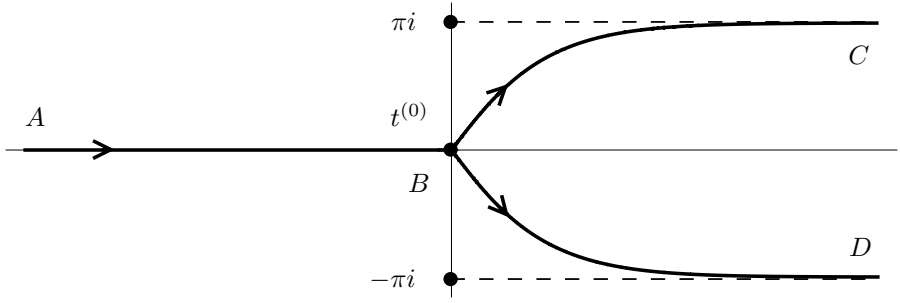
\caption{The steepest descent paths for $H_\nu ^{\left( 1 \right)} \left(\nu\right)$ through the saddle point $t^{\left( 0 \right)} = 0$ when $\arg \nu =0$.}
\label{fig4}
\end{figure}

\subsubsection{Case (ii): $x=1$} For the integral \eqref{eq3}, we need to consider only the saddle point $t^{\left( 0 \right)} = 0$. Let $f\left(t\right) = \sinh t-t$. The steepest descent paths in this case are given by
\[
\arg \left[ {e^{i\theta } \left( {f\left( t^{\left( 0 \right)} \right) - f\left( t \right)} \right)} \right] = 0 .
\]
The path $\mathscr{C}^{\left( 0 \right)}\left(\theta\right)$ for $\theta=0$, with an appropriate orientation, is shown in Figure \ref{fig4}. We assume that $\theta = 0$ and later we shall use an analytic continuation argument to extend the results to complex $\nu$. Denote by $\mathscr{P}^{\left( 0 \right)}\left(0\right)$ and $\mathscr{L}^{\left( 0 \right)}\left(0\right)$ the curves $ABC$ and $ABD$, respectively (see Figure \ref{fig4}). By Cauchy's theorem, we can deform the path of integration in \eqref{eq3} to $\mathscr{P}^{\left( 0 \right)}\left(0\right)$. (The path $\mathscr{L}^{\left( 0 \right)}\left(0\right)$ can be taken as the integration path in the corresponding integral representation for $H_\nu^{\left(2\right)}\left(\nu\right)$.) Let
\[
\tau = f(t^{\left( 0 \right)} ) - f\left( t \right)=  - f\left( t \right).
\]
By the definition of $\mathscr{C}^{\left( 0 \right)} \left( 0  \right)$, for $t \in \mathscr{P}^{\left( 0 \right)} \left( 0  \right)$, $\tau$ is real and positive. Then, corresponding to each positive $\tau$, there are two values of $t$, called $t_1$ and $t_2$, such that $t_1$ is a complex number with a positive real part and $t_2$ is a negative real number. In terms of $\tau$, \eqref{eq3} becomes
\begin{align*}
H_\nu ^{\left( 1 \right)} \left( \nu  \right) & = \frac{1}{{\pi i}}\int_0^{ + \infty } {e^{ - \nu \tau } \left( {\frac{{dt_1 }}{{d\tau }} - \frac{{dt_2 }}{{d\tau }}} \right)d\tau } \\ &  = \frac{1}{{\pi i}}\int_0^{ + \infty } {e^{ - \nu \tau } \left( {\frac{1}{{1 - \cosh t_1 \left( \tau  \right)}} - \frac{1}{{1 - \cosh t_2 \left( \tau  \right)}}} \right)d\tau } .
\end{align*}
The function in the large parentheses is an even function of $\tau$ and can be represented as a contour integral using the residue theorem, to yield
\[
H_\nu ^{\left( 1 \right)} \left( \nu  \right) = \frac{1}{{3\pi i}}\int_0^{ + \infty } {\tau ^{ - \frac{2}{3}} e^{ - \nu \tau } \frac{1}{{2\pi i}}\oint_{\Gamma ^{\left( 0 \right)} } {\left( {\frac{{e^{\frac{{\pi i}}{3}} }}{{1 - \tau ^{\frac{2}{3}} f^{ - \frac{2}{3}} \left( u \right)e^{\frac{{2\pi i}}{3}} }} + \frac{1}{{1 - \tau ^{\frac{2}{3}} f^{ - \frac{2}{3}} \left( u \right)}}} \right)\frac{{du}}{{f^{\frac{1}{3}} \left( u \right)}}} d\tau } ,
\]
where the contour $\Gamma ^{\left( 0 \right)}$ encloses $\mathscr{P}^{\left( 0 \right)} \left( 0  \right)$ in an anti-clockwise loop, and does not enclose any of the saddle points $t ^{\left( k \right)} \neq t^{\left( 0 \right)}$. The cube root is defined so that $\left(-f\left( t \right) \right)^{\frac{1}{3}}$ is positive on the portion of $\mathscr{P}^{\left( 0 \right)}\left(0\right)$ that starts at $t=-\infty$ and ends at $t= t^{\left( 0 \right)}$. Next we apply the expression \eqref{eq7} to expand the function under the contour integral in powers of $\tau ^{\frac{2}{3}} f^{ - \frac{2}{3}} \left( u \right)$. The result is
\[
H_\nu ^{\left( 1 \right)} \left( \nu  \right) = \frac{1}{{3\pi i}}\sum\limits_{n = 0}^{N - 1} {\left( {1 + e^{\frac{{\left( {2n + 1} \right)\pi i}}{3}} } \right)\int_0^{ + \infty } {\tau ^{\frac{{2n - 2}}{3}} e^{ - \nu \tau } \frac{1}{{2\pi i}}\oint_{\Gamma ^{\left( 0 \right)} } {\frac{{du}}{{f^{\frac{{2n + 1}}{3}} \left( u \right)}}} d\tau } }  + R_N^{\left( H \right)} \left( \nu  \right),
\]
where
\begin{equation}\label{eq6}
R_N^{\left( H \right)} \left( \nu  \right) = \frac{1}{{3\pi i}}\int_0^{ + \infty } {\tau ^{\frac{{2N - 2}}{3}} e^{ - \nu \tau } \frac{1}{{2\pi i}}\oint_{\Gamma ^{\left( 0 \right)} } {\left( {\frac{{e^{\frac{{\left( {2N + 1} \right)\pi i}}{3}} }}{{1 - \tau ^{\frac{2}{3}} f^{ - \frac{2}{3}} \left( u \right)e^{\frac{{2\pi i}}{3}} }} + \frac{1}{{1 - \tau ^{\frac{2}{3}} f^{ - \frac{2}{3}} \left( u \right)}}} \right)\frac{{du}}{{f^{\frac{{2N + 1}}{3}} \left( u \right)}}} d\tau } .
\end{equation}
The path $\Gamma ^{\left( 0 \right)}$ in the sum can be shrunk into a small circle around $t^{\left( 0 \right)} = 0$, and we arrive at
\[
H_\nu ^{\left( 1 \right)} \left( \nu  \right) =  - \frac{2}{{3\pi }}\sum\limits_{n = 0}^{N - 1} {d_{2n} e^{\frac{{2\left( {2n + 1} \right)\pi i}}{3}} \sin \left( {\frac{{\left( {2n + 1} \right)\pi }}{3}} \right)\frac{{\Gamma \left( {\frac{{2n + 1}}{3}} \right)}}{{\nu ^{\frac{{2n + 1}}{3}} }}}  + R_N^{\left( H \right)} \left( \nu  \right),
\]
where
\[
d_{2n}  = \frac{1}{{2\pi i}}\oint_{\left( {0^ +  } \right)} {\frac{{du}}{{f^{\frac{{2n + 1}}{3}} \left( u \right)}}}  = \frac{1}{{2\pi i}}\oint_{\left( {0^ +  } \right)} {\left( {\frac{{u^3 }}{{f\left( u \right)}}} \right)^{\frac{{2n + 1}}{3}} \frac{{du}}{{u^{2n + 1} }}}  = \frac{1}{{\left( {2n} \right)!}}\left[ {\frac{{d^{2n} }}{{dt^{2n} }}\left( {\frac{{t^3 }}{{\sinh t - t}}} \right)^{\frac{{2n + 1}}{3}} } \right]_{t = 0} .
\]

\begin{figure}[!t]
\def\svgwidth{\columnwidth}
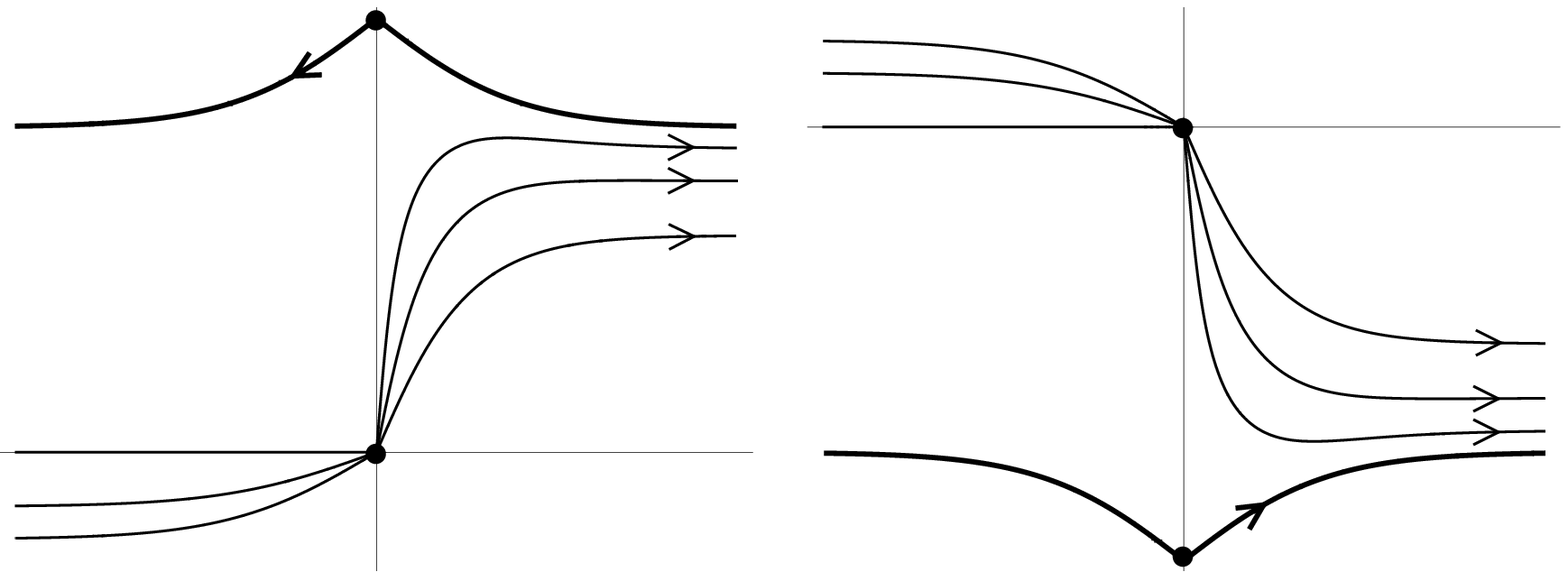
\caption{(a) The steepest descent path for $H_\nu ^{\left( 1 \right)} \left(\nu\right)$ through the saddle point $t^{\left( 0 \right)}$ when (i) $\arg \nu =0$, (ii) $\arg \nu =-\frac{\pi}{4}$, (iii) $\arg \nu =-\frac{2\pi}{5}$. The path $\mathscr{L}^{\left( 1 \right)} \left( { - \frac{\pi }{2}} \right)$ is an adjacent contour to $t^{\left( 0 \right)}$. (b) The steepest descent path $\mathscr{L}^{\left( 0 \right)} \left( \theta \right)$ through the saddle point $t^{\left( 0 \right)}$ when (i) $\arg \nu =0$, (ii) $\arg \nu = \frac{\pi}{4}$, (iii) $\arg \nu = \frac{2\pi}{5}$. The path $\mathscr{P}^{\left(-1\right)}\left({\frac{\pi}{2}}\right)$ is an adjacent contour to $t^{\left( 0 \right)}$.}
\label{fig5}
\end{figure}

Applying the change of variable $\nu \tau = s$ in \eqref{eq6} gives
\begin{equation}\label{eq8}
R_N^{\left( H \right)} \left( \nu  \right) = \frac{1}{{3\pi i\nu ^{\frac{{2N + 1}}{3}} }}\int_0^{ + \infty } {s^{\frac{{2N - 2}}{3}} e^{ - s} \frac{1}{{2\pi i}}\oint_{\Gamma ^{\left( 0 \right)} } {\left( {\frac{{e^{\frac{{\left( {2N + 1} \right)\pi i}}{3}} }}{{1 - \left( {s/\nu f\left( u \right)} \right)^{\frac{2}{3}} e^{\frac{{2\pi i}}{3}} }} + \frac{1}{{1 - \left( {s/\nu f\left( u \right)} \right)^{\frac{2}{3}} }}} \right)\frac{{du}}{{f^{\frac{{2N + 1}}{3}} \left( u \right)}}} ds} .
\end{equation}
As in the first case, we need to locate the adjacent saddle points. As $\theta$ varies, the path $\mathscr{C}^{\left(0\right)}\left(\theta\right)$ varies smoothly, so we can define the paths $\mathscr{P}^{\left(0\right)}\left(\theta\right)$ and $\mathscr{L}^{\left(0\right)}\left(\theta\right)$. In our case, the relevant path is $\mathscr{P}^{\left(0\right)}\left(\theta\right)$, for the adjacency problem, however, we also need to consider $\mathscr{L}^{\left(0\right)}\left(\theta\right)$. When $\theta = -\frac{\pi}{2}$ or $\frac{3\pi}{2}$, the path $\mathscr{P}^{\left(0\right)}\left(\theta\right)$ connects to the saddle point $t^{\left( 1 \right)} = 2\pi i$. Similarly, when $\theta = \frac{\pi}{2}$ or $-\frac{3\pi}{2}$, the path $\mathscr{L}^{\left(0\right)}\left(\theta\right)$ connects to the saddle point $t^{\left( -1 \right)} = -2\pi i$. Therefore, the adjacent saddles are $t^{\left( \pm 1 \right)}$ (see Figure \ref{fig5}). The set
\[
\Delta ^{\left( 0 \right)}  = \left\{ {u \in \mathscr{P}^{\left( 0 \right)} \left( \theta  \right) : - \frac{\pi }{2} < \theta  < \frac{{3\pi }}{2}} \right\} \cup \left\{ {u \in \mathscr{L}^{\left( 0 \right)} \left( \theta  \right) : - \frac{3\pi }{2} < \theta  < \frac{{\pi }}{2}} \right\}
\]
forms a domain in the complex plane whose boundaries are the steepest descent paths $\mathscr{L}^{\left( 1 \right)} \left( { - \frac{\pi }{2}} \right)$ and $\mathscr{P}^{\left( -1 \right)} \left( { \frac{\pi }{2}} \right)$, the adjacent contours to $t^{\left( 0 \right)}$ (they are defined analogously to $\mathscr{P}^{\left(0\right)}$ and $\mathscr{L}^{\left(0\right)}$). The function under the contour integral in \eqref{eq8} is an analytic function of $u$ in the domain $\Delta ^{\left( 0 \right)}$, therefore we can deform $\Gamma^{\left( 0 \right)}$ over the adjacent contours. We thus find that for $-\frac{\pi}{2} < \theta < \frac{3\pi}{2}$ and $N \geq 0$, \eqref{eq8} can be written
\begin{gather}\label{eq9}
\begin{split}
R_N^{\left( H \right)} \left( \nu  \right) = \; & \frac{1}{{3\pi i\nu ^{\frac{{2N + 1}}{3}} }}\int_0^{ + \infty } {s^{\frac{{2N - 2}}{3}} e^{ - s} \frac{1}{{2\pi i}}\int_{\mathscr{L}^{\left( 1 \right)} \left( { - \frac{\pi }{2}} \right)} {\left( {\frac{{e^{\frac{{\left( {2N + 1} \right)\pi i}}{3}} }}{{1 - \left( {s/\nu f\left( u \right)} \right)^{\frac{2}{3}} e^{\frac{{2\pi i}}{3}} }} + \frac{1}{{1 - \left( {s/\nu f\left( u \right)} \right)^{\frac{2}{3}} }}} \right)\frac{{du}}{{f^{\frac{{2N + 1}}{3}} \left( u \right)}}} ds}\\
& + \frac{1}{{3\pi i\nu ^{\frac{{2N + 1}}{3}} }}\int_0^{ + \infty } {s^{\frac{{2N - 2}}{3}} e^{ - s} \frac{1}{{2\pi i}}\int_{\mathscr{P}^{\left( -1 \right)} \left( { \frac{\pi }{2}} \right)} {\left( {\frac{{e^{\frac{{\left( {2N + 1} \right)\pi i}}{3}} }}{{1 - \left( {s/\nu f\left( u \right)} \right)^{\frac{2}{3}} e^{\frac{{2\pi i}}{3}} }} + \frac{1}{{1 - \left( {s/\nu f\left( u \right)} \right)^{\frac{2}{3}} }}} \right)\frac{{du}}{{f^{\frac{{2N + 1}}{3}} \left( u \right)}}} ds}.
\end{split}
\end{gather}
Now we perform the changes of variable
\[
s = t\frac{{\left| {f\left( {2\pi i} \right) - f\left( 0 \right)} \right|}}{{f\left( {2\pi i} \right) - f\left( 0 \right)}}\left( {f\left( u \right) - f\left( 0 \right)} \right) = itf\left( u \right)
\]
in the first, and
\[
s = t\frac{{\left| {f\left( { - 2\pi i} \right) - f\left( 0 \right)} \right|}}{{f\left( { - 2\pi i} \right) - f\left( 0 \right)}}\left( {f\left( u \right) - f\left( 0 \right)} \right) =  - itf\left( u \right)
\]
in the second double integral. In this case, Dingle's singulants are $f\left( { \pm 2\pi i} \right) - f\left( 0 \right)=  \mp 2\pi i$. When using these changes of variable, we should take $i^{\frac{2}{3}}  = -1$ in the first, and $\left( { - i} \right)^{\frac{2}{3}}  = -1$ in the second double integral. With these changes of variable, the representation \eqref{eq9} for $R_N^{\left(H\right)} \left(\nu\right)$ becomes
\begin{gather}\label{eq16}
\begin{split}
R_N^{\left( H \right)} \left( \nu  \right) = \; & \frac{{\left( { - 1} \right)^N }}{{3\pi \nu ^{\frac{{2N + 1}}{3}} }}\int_0^{ + \infty } {t^{\frac{{2N - 2}}{3}} \left( {\frac{{e^{\frac{{\left( {2N + 1} \right)\pi i}}{3}} }}{{1 + \left( {t/\nu } \right)^{\frac{2}{3}} e^{\frac{{2\pi i}}{3}} }} + \frac{1}{{1 + \left( {t/\nu } \right)^{\frac{2}{3}} }}} \right)\frac{{ - 1}}{{2\pi i}}\int_{\mathscr{L}^{\left( 1 \right)} \left( { - \frac{\pi }{2}} \right)} {e^{ - itf\left( u \right)} du} dt} \\
& + \frac{{\left( { - 1} \right)^N }}{{3\pi \nu ^{\frac{{2N + 1}}{3}} }}\int_0^{ + \infty } {t^{\frac{{2N - 2}}{3}} \left( {\frac{{e^{\frac{{\left( {2N + 1} \right)\pi i}}{3}} }}{{1 + \left( {t/\nu } \right)^{\frac{2}{3}} e^{\frac{{2\pi i}}{3}} }} + \frac{1}{{1 + \left( {t/\nu } \right)^{\frac{2}{3}} }}} \right)\frac{1}{{2\pi i}}\int_{\mathscr{P}^{\left( { - 1} \right)} \left( {\frac{\pi }{2}} \right)} {e^{itf\left( u \right)} du} dt} ,
\end{split}
\end{gather}
for $-\frac{\pi}{2} < \theta < \frac{3\pi}{2}$ and $N \geq 0$. Finally, the contour integrals can themselves be represented in terms of the Hankel functions since
\begin{align*}
\frac{{ - 1}}{{\pi i}}\int_{\mathscr{L}^{\left( 1 \right)} \left( { - \frac{\pi }{2}} \right)} {e^{ - itf\left( u \right)} du} & = \frac{1}{\pi i}\int_{\mathscr{L}^{\left( 0 \right)} \left( { - \frac{\pi }{2}} \right)} {e^{ - itf\left( {u + 2\pi i} \right)} du}  =  \frac{{e^{ - 2\pi t} }}{{\pi i}}\int_{\mathscr{L}^{\left( 0 \right)} \left( { - \frac{\pi }{2}} \right)} {e^{ - itf\left( u \right)} du}  \\ & =  - e^{ - 2\pi t} H_{ - it}^{\left( 2 \right)} \left( { - it} \right) = e^{ - 2\pi t} H_{it}^{\left( 1 \right)} \left( {it} \right),
\end{align*}
and
\[
\frac{1}{{\pi i}}\int_{\mathscr{P}^{\left( { - 1} \right)} \left( {\frac{\pi }{2}} \right)} {e^{itf\left( u \right)} du}  = \frac{1}{{\pi i}}\int_{\mathscr{P}^{\left( 0 \right)} \left( {\frac{\pi }{2}} \right)} {e^{itf\left( {u - 2\pi i} \right)} du}  = \frac{{e^{ - 2\pi t} }}{{\pi i}}\int_{\mathscr{P}^{\left( 0 \right)} \left( {\frac{\pi }{2}} \right)} {e^{itf\left( u \right)} du}  = e^{ - 2\pi t} H_{it}^{\left( 1 \right)} \left( {it} \right).
\]
Substituting these into \eqref{eq16} gives \eqref{eq14}. Formula \eqref{eq18} follows from the connection formula $H_\nu ^{\left( 2 \right)} \left( \nu \right) = \overline {H_{ \bar{\nu} }^{\left( 1 \right)} \left( {\bar \nu } \right)}$. Finally, we prove the second representation in \eqref{eq15}. To avoid complications caused by the zeros of the sine function, we proceed in a different way than in the case of $U_n\left(i \cot \beta\right)$. We note that
\begin{align*}
\Gamma \left( {\frac{{2n + 1}}{3}} \right)d_{2n} & = \int_0^{ + \infty } {s^{\frac{{2n - 2}}{3}} e^{ - s} \frac{1}{2\pi i} \oint_{\left( {0^ +  } \right)} {\frac{{du}}{{f^{\frac{{2n + 1}}{3}} \left( u \right)}}} ds} \\
& = \int_0^{ + \infty } {s^{\frac{{2n - 2}}{3}} e^{ - s} \frac{1}{2\pi i} \int_{\mathscr{L}^{\left( 1 \right)} \left( { - \frac{\pi }{2}} \right)} {\frac{{du}}{{f^{\frac{{2n + 1}}{3}} \left( u \right)}}} ds}  + \int_0^{ + \infty } {s^{\frac{{2n - 2}}{3}} e^{ - s} \frac{1}{2\pi i}\int_{\mathscr{P}^{\left( { - 1} \right)} \left( {\frac{\pi }{2}} \right)} {\frac{{du}}{{f^{\frac{{2n + 1}}{3}} \left( u \right)}}} ds} .
\end{align*}
For these integrals we can apply the same changes of variable as above, then using the same argument yields the second representation in \eqref{eq15}.

\subsection{The Bessel functions $J_\nu\left(\nu x\right)$ and $Y_\nu\left(\nu x\right)$} To prove the resurgence representations for the Bessel functions, we apply the well-known connection formulas
\begin{equation}\label{eq39}
J_\nu  \left( z \right) = \frac{1}{2}\left( {H_\nu ^{\left( 1 \right)} \left( z \right) + H_\nu ^{\left( 2 \right)} \left( z \right)} \right),
\end{equation}
\[
Y_\nu  \left( z \right) = \frac{1}{{2i}}\left( {H_\nu ^{\left( 1 \right)} \left( z \right) - H_\nu ^{\left( 2 \right)} \left( z \right)} \right).
\]
To show \eqref{eq22} and \eqref{eq29}, we apply \eqref{eq21} and \eqref{eq17} with $2N$ in place of $N$ and substitute them into \eqref{eq39}. Making use of the identity $e^{i \xi} = \cos \xi +i \sin \xi$ then gives the desired result. The proofs of \eqref{eq23} and \eqref{eq30} are analogous.

To prove \eqref{eq40} and \eqref{eq42}, we substitute \eqref{eq44} and \eqref{eq18} into \eqref{eq39} and note that
\[
\frac{{e^{\frac{{2\left( {2n + 1} \right)\pi i}}{3}}  + e^{ - \frac{{2\left( {2n + 1} \right)\pi i}}{3}} }}{2}\sin \left( {\frac{{\left( {2n + 1} \right)\pi }}{3}} \right) = \cos \left( {\frac{{2\left( {2n + 1} \right)\pi }}{3}} \right)\sin \left( {\frac{{\left( {2n + 1} \right)\pi }}{3}} \right) =  - \frac{1}{2}\sin \left( {\frac{{\left( {2n + 1} \right)\pi }}{3}} \right).
\]
The proofs of \eqref{eq41} and \eqref{eq43} are similar.

\section{Error bounds for Debye's expansions}\label{section3}

In this section we derive explicit and numerically computable error bounds for the large order asymptotic series of the Hankel and Bessel functions. The proofs are based on the resurgence formulas given in Theorems \ref{thm1} and \ref{thm2}. We consider the expansions of the Hankel functions in the first, and the expansions of the Bessel functions in the second, subsection.

We comment on the relation between Meijer's work on the Debye expansions \cite{Meijer} and ours for the case $x>1$. Some of the estimates in \cite{Meijer} coincide with ours and are valid in wider sectors of the complex $\nu$-plane. However, it should be noted that those bounds become less effective outside the sectors of validity of the representations \eqref{eq12}--\eqref{eq30} due to the Stokes phenomenon. For those sectors we recommend the use of the continuation formulas given in Section \ref{section1}.

To estimate the remainder terms, we shall use frequently the elementary result that
\begin{equation}\label{eq20}
\frac{1}{{\left| {1 - re^{i\varphi } } \right|}} \le \begin{cases} \left|\csc \varphi \right| & \; \text{ if } \; 0 < \left|\varphi \text{ mod } 2\pi\right| <\frac{\pi}{2} \\ 1 & \; \text{ if } \; \frac{\pi}{2} \leq \left|\varphi \text{ mod } 2\pi\right| \leq \pi \end{cases}
\end{equation}
holds for any $r>0$. We will also need the fact that
\begin{equation}\label{eq19}
iH_{it}^{\left(1\right)} \left( {itx} \right) \ge 0
\end{equation}
for any $t>0$ and $x\geq 1$. This follows from the connection formula with the modified Bessel function of the third kind of purely imaginary order
\[
iH_{it}^{\left( 1 \right)} \left( {itx} \right) = \frac{2}{\pi }e^{\frac{\pi}{2}t} K_{it} \left( {tx} \right) .
\]
It is known that $K_{it} \left( {tx} \right)$ is real and non-zero for $t>0$ and $x\geq 1$, and, as its uniform asymptotic expansion shows, it is positive for large $t$ and fixed $x\geq 1$ (see, e.g., Dunster \cite{Dunster}).

\subsection{Error bounds for the expansions of $H_\nu^{\left(1\right)}\left(\nu x\right)$ and $H_\nu^{\left(2\right)}\left(\nu x\right)$}

\subsubsection{Case (i): $x>1$} As usual, let $0<\beta<\frac{\pi}{2}$ be defined by $\sec \beta = x$. We observe that from \eqref{eq13} and \eqref{eq19} it follows that
\begin{equation}\label{eq28}
\left| {U_n \left( {i\cot \beta } \right)} \right| = \frac{1}{{2\left( {2\pi \cot \beta } \right)^{\frac{1}{2}} }}\int_0^{ + \infty } {t^{n - \frac{1}{2}} e^{ - t\left( {\tan \beta  - \beta } \right)} \left( {1 + e^{ - 2\pi t} } \right)iH_{it}^{\left( 1 \right)} \left( {it\sec \beta } \right)dt} .
\end{equation}
Using this formula, together with the estimate \eqref{eq20} and the representation \eqref{eq12}, we obtain the error bound
\begin{equation}\label{eq24}
\left| {R_N^{\left( H \right)} \left( {\nu ,\beta } \right)} \right|  \le \frac{{\left| {U_N \left( {i\cot \beta } \right)} \right|}}{{\left| \nu  \right|^N }} \begin{cases} \left|\sec \theta \right| & \; \text{ if } \; { - \frac{\pi }{2} < \theta  < 0 \; \text{ or } \; \pi  < \theta  < \frac{{3\pi }}{2}} \\ 1 & \; \text{ if } \; {0 \le \theta  \le \pi }. \end{cases}
\end{equation}
The remainder term of the expansion \eqref{eq17} can be estimated in exactly the same way, the result is
\begin{equation}\label{eq73}
\left| {R_N^{\left( H \right)} \left( {\nu e^{\pi i} ,\beta } \right)} \right|  \le \frac{{\left| {U_N \left( {i\cot \beta } \right)} \right|}}{{\left| \nu  \right|^N }} \begin{cases} \left|\sec \theta \right| & \; \text{ if } \; { 0 < \theta  < \frac{\pi}{2} \; \text{ or } \; -\frac{3\pi}{2}}  < \theta  < -\pi \\ 1 & \; \text{ if } \; {-\pi \le \theta \le 0 }. \end{cases}
\end{equation}
These error bounds become singular as $\theta \to -\frac{\pi }{2}, \frac{3\pi }{2}$ and $\theta \to -\frac{3\pi }{2}, \frac{\pi}{2}$, and therefore unrealistic near the Stokes lines. A better bound for $R_N^{\left( H \right)} \left( {\nu ,\beta } \right)$ near the Stokes lines $\theta = -\frac{\pi }{2}, \frac{3\pi }{2}$ can be derived as follows. Let $0 < \varphi  < \frac{\pi }{2}$ be an acute angle that may depend on $N$. Suppose that $ - \frac{\pi }{2} \le \theta  <  - \varphi$. An analytic continuation of the representation \eqref{eq21} to this sector can be found by rotating the path of integration in \eqref{eq12} by $-\varphi$:
\[
R_N^{\left( H \right)} \left( {\nu ,\beta } \right)  =  \frac{1}{{2\left( {2\pi \cot \beta } \right)^{\frac{1}{2}} \left(i\nu\right)^N }}\int_0^{ + \infty e^{ - i\varphi } } {\frac{{t^{N - \frac{1}{2}} e^{ - t\left( {\tan \beta  - \beta } \right)} }}{{1 + it/\nu }}\left( {1 + e^{ - 2\pi t} } \right)i H_{it}^{\left( 1 \right)} \left( {it\sec \beta } \right)dt} .
\]
Substituting $t = \frac{s}{{e^{i\varphi } \cos \varphi }}$ and applying the estimation \eqref{eq20}, we obtain
\begin{equation}\label{eq25}
\left| {R_N^{\left( H \right)} \left( {\nu ,\beta } \right)} \right| \le \frac{{\sec \left( {\theta  + \varphi } \right)}}{{2\left( {2\pi \cot \beta } \right)^{\frac{1}{2}} \cos^{N + \frac{1}{2}} \varphi \left| \nu  \right|^N }}\int_0^{ + \infty } {s^{N - \frac{1}{2}} e^{ - s\left( {\tan \beta  - \beta } \right)} \left| {1 + e^{ - \frac{{2\pi s}}{{e^{i\varphi } \cos \varphi }}} } \right|\left| {H_{\frac{{is}}{{e^{i\varphi } \cos \varphi }}}^{\left( 1 \right)} \left( {\frac{{is}}{{e^{i\varphi } \cos \varphi }}\sec \beta } \right)} \right|ds} .
\end{equation}
We would like to have $\left| {U_n \left( {i\cot \beta } \right)} \right|$ in our error bound. First, we observe that
\begin{equation}\label{eq32}
\left| {1 + e^{ - \frac{{2\pi s}}{{e^{i\varphi } \cos \varphi }}} } \right| \le 1 + e^{ - 2\pi s} .
\end{equation}
To make the formula \eqref{eq28} applicable, we need to ``replace" the modulus of the Hankel function inside the integral by $iH_{is}^{\left( 1 \right)} \left( {is\sec \beta } \right)$. To this end, we employ the formula \eqref{eq12} with $N=0$, to find
\begin{gather}\label{eq26}
\begin{split}
\left| {H_{\frac{{is}}{{e^{i\varphi } \cos \varphi }}}^{\left( 1 \right)} \left( {\frac{{is}}{{e^{i\varphi } \cos \varphi }}\sec \beta } \right)} \right| & = \left| {\frac{{e^{i\frac{{is}}{{e^{i\varphi } \cos \varphi }}\left( {\tan \beta  - \beta } \right) + \frac{\pi }{4}i} }}{{2\pi \left( {\frac{{is}}{{e^{i\varphi } \cos \varphi }}} \right)^{\frac{1}{2}} }}\int_0^{ + \infty } {\frac{{t^{ - \frac{1}{2}} e^{ - t\left( {\tan \beta  - \beta } \right)} }}{{1 + e^{i\varphi } \cos \varphi t/s}}\left( {1 + e^{ - 2\pi t} } \right)H_{it}^{\left( 1 \right)} \left( {it\sec \beta } \right)dt} } \right| \\
& \le \sqrt {\cos \varphi } \frac{{e^{ - s\left( {\tan \beta  - \beta } \right)} }}{{2\pi s^{\frac{1}{2}} }}\int_0^{ + \infty } {\frac{{t^{ - \frac{1}{2}} e^{ - t\left( {\tan \beta  - \beta } \right)} }}{{\left| {1 + e^{i\varphi } \cos \varphi t/s} \right|}}\left( {1 + e^{ - 2\pi t} } \right)iH_{it}^{\left( 1 \right)} \left( {it\sec \beta } \right)dt} \\
& \le \frac{1}{{\sqrt {\cos \varphi } }}\frac{{e^{ - s\left( {\tan \beta  - \beta } \right)} }}{{2\pi s^{\frac{1}{2}} }}\int_0^{ + \infty } {\frac{{t^{ - \frac{1}{2}} e^{ - t\left( {\tan \beta  - \beta } \right)} }}{{1 + t/s}}\left( {1 + e^{ - 2\pi t} } \right)iH_{it}^{\left( 1 \right)} \left( {it\sec \beta } \right)dt} \\
& = \frac{1}{{\sqrt {\cos \varphi } }}\left| {H_{is}^{\left( 1 \right)} \left( {is\sec \beta } \right)} \right| = \frac{1}{{\sqrt {\cos \varphi } }}iH_{is}^{\left( 1 \right)} \left( {is\sec \beta } \right) .
\end{split}
\end{gather}
Finally, we estimate the ratio $\sec \left(\theta + \varphi \right)\cos^{-N-1} \varphi$. It is easy to show that the value $\varphi = \arctan \left( {\left(N+1\right)^{ - \frac{1}{2}} } \right)$ minimises the function $\sec \left( { - \frac{\pi }{2} + \varphi } \right)\cos^{ - N-1} \varphi$, and
\begin{align*}
\frac{{\sec \left( {\theta  + \arctan \left( {\left(N+1\right)^{ - \frac{1}{2}} } \right)} \right)}}{{\cos ^{N+1} \left( {\arctan \left( {\left(N+1\right)^{ - \frac{1}{2}} } \right)} \right)}} \le \frac{{\sec \left( { - \frac{\pi }{2} + \arctan \left( {\left(N+1\right)^{ - \frac{1}{2}} } \right)} \right)}}{{\cos^{N+1} \left( {\arctan \left( {\left(N+1\right)^{ - \frac{1}{2}} } \right)} \right)}} & = \left( {1 + \frac{1}{N+1}} \right)^{\frac{N}{2}+1} \sqrt{N+1} \\ & \le \sqrt {e\left( {N + \frac{3}{2}} \right)} 
\end{align*}
for all $ - \frac{\pi }{2} \le \theta  <  - \varphi = - \arctan \left( {\left(N+1\right)^{ - \frac{1}{2}} } \right)$ with $N \geq 0$. Here we used that the sequence $\left( {1 + \frac{1}{{N + 1}}} \right)^{N/2 + 1} \sqrt {\frac{{N + 1}}{{N + 3/2}}}$ is increasing and its limit is $\sqrt{e}$. Substituting all this estimations into \eqref{eq25} yields the error bound
\begin{equation}\label{eq27}
\left| {R_N^{\left( H \right)} \left( {\nu ,\beta } \right)} \right|  \le \sqrt {e\left( {N + \frac{3}{2}} \right)}  \frac{{\left| {U_N \left( {i\cot \beta } \right)} \right|}}{{\left| \nu  \right|^N }},
\end{equation}
which is valid for $ - \frac{\pi }{2} \le \theta  <  - \arctan \left( {\left(N+1\right)^{ - \frac{1}{2}} } \right)$ with $N \geq 0$. A similar argument shows that this bound is also valid when $\pi  + \arctan \left( {\left(N+1\right)^{ - \frac{1}{2}} } \right) < \theta  \le \frac{3\pi}{2}$. In the ranges $-\frac{\pi}{4} \le \theta  < 0$ and $\pi < \theta  \le \frac{5\pi}{4}$, it holds that $\sqrt {e\left( {N + \frac{3}{2}} \right)}  \ge \sqrt {\frac{3}{2}e}  \ge \left| {\sec \theta } \right|$, whence the estimate \eqref{eq27} remains valid for $ - \frac{\pi}{2} \le \theta  <  0$ and $\pi < \theta  \le \frac{3\pi}{2}$.

The remainder term $R_N^{\left( H \right)} \left( {\nu e^{\pi i},\beta } \right)$ of the expansion \eqref{eq17} can be estimated similarly and we obtain the same upper bound as in \eqref{eq27} but with the modified conditions $0 < \theta  \le \frac{\pi}{2}$ and $ - \frac{3\pi }{2} \le \theta  <  -\pi$.

\subsubsection{Case (ii): $x=1$}\label{bounds} We note that from \eqref{eq15} and \eqref{eq19} it follows that
\begin{equation}\label{eq37}
\left| {d_{2n} } \right| = \frac{1}{{\Gamma \left( {\frac{{2n + 1}}{3}} \right)}}\int_0^{ + \infty } {t^{\frac{{2n - 2}}{3}} e^{ - 2\pi t} iH_{it}^{\left( 1 \right)} \left( {it} \right)dt} .
\end{equation}
A little manipulation of \eqref{eq14} shows that
\begin{align}
R_N^{\left( H \right)} \left( \nu  \right) & = \frac{{\left( { - 1} \right)^{N+1} }}{{\sqrt 3 \pi \nu ^{\frac{{2N + 1}}{3}} }}\int_0^{ + \infty } {t^{\frac{{2N - 2}}{3}} e^{ - 2\pi t} e^{\frac{{2\left( {2N + 1} \right)\pi i}}{3}} \frac{{1 + \left( {t/\nu } \right)^{\frac{2}{3}} e^{\frac{\pi}{3}i} }}{{\left( {1 + \left( {t/\nu } \right)^{\frac{2}{3}} e^{\frac{{2\pi i}}{3}} } \right)\left( {1 + \left( {t/\nu } \right)^{\frac{2}{3}} } \right)}}i H_{it}^{\left( 1 \right)} \left( {it} \right)dt} ,\nonumber \\
R_N^{\left( H \right)} \left( \nu  \right) & = \frac{{\left( { - 1} \right)^N }}{{\sqrt 3 \pi \nu ^{\frac{{2N + 3}}{3}} }}\int_0^{ + \infty } {t^{\frac{{2N}}{3}} e^{ - 2\pi t} e^{\frac{{2\left( {2N + 1} \right)\pi i}}{3}} \frac{{e^{\frac{\pi}{3}i} }}{{\left( {1 + \left( {t/\nu } \right)^{\frac{2}{3}} e^{\frac{{2\pi i}}{3}} } \right)\left( {1 + \left( {t/\nu } \right)^{\frac{2}{3}} } \right)}}i H_{it}^{\left( 1 \right)} \left( {it} \right)dt} ,\label{eq35} \\
R_N^{\left( H \right)} \left( \nu  \right) & = \frac{{\left( { - 1} \right)^N }}{{\sqrt 3 \pi \nu ^{\frac{{2N + 1}}{3}} }}\int_0^{ + \infty } {t^{\frac{{2N - 2}}{3}} e^{ - 2\pi t} e^{\frac{{2\left( {2N + 1} \right)\pi i}}{3}} \frac{1}{{\left( {1 + \left( {t/\nu } \right)^{\frac{2}{3}} e^{\frac{{2\pi i}}{3}} } \right)\left( {1 + \left( {t/\nu } \right)^{\frac{2}{3}} } \right)}}i H_{it}^{\left( 1 \right)} \left( {it} \right)dt} \nonumber
\end{align}
according to whether $N\equiv 0 \mod 3$, $N\equiv 1 \mod 3$ or $N\equiv 2 \mod 3$, respectively. We show in Appendix \ref{appendixb} that
\begin{equation}\label{eq47}
\frac{1}{{\left| {1 + \left( {t/\nu } \right)^{\frac{2}{3}} e^{\frac{{2\pi i}}{3}} } \right|\left| {1 + \left( {t/\nu } \right)^{\frac{2}{3}} } \right|}} \le \begin{cases} \left|\sec \theta \right| & \; \text{ if } \; { - \frac{\pi }{2} < \theta  < 0 \; \text{ or } \; \pi  < \theta  < \frac{{3\pi }}{2}} \\ 1 & \; \text{ if } \; {0 \le \theta  \le \pi }. \end{cases}
\end{equation}
Applying this estimate together with \eqref{eq37} and \eqref{eq35} yields the error bound
\begin{equation}\label{eq51}
\left| {R_N^{\left( H \right)} \left( \nu  \right)} \right| \le \frac{2}{{3\pi }}\left| {d_{2N + 2} } \right|\frac{{\sqrt 3 }}{2}\frac{{\Gamma \left( {\frac{{2N + 3}}{3}} \right)}}{{\left| \nu  \right|^{\frac{{2N + 3}}{3}} }} \begin{cases} \left|\sec \theta \right| & \; \text{ if } \; { - \frac{\pi }{2} < \theta  < 0 \; \text{ or } \; \pi  < \theta  < \frac{{3\pi }}{2}} \\ 1 & \; \text{ if } \; {0 \le \theta  \le \pi }, \end{cases}
\end{equation}
when $N\equiv 1 \mod 3$. Similarly, we have
\begin{equation}\label{eq58}
\left| {R_N^{\left( H \right)} \left( \nu  \right)} \right| \le \frac{2}{{3\pi }}\left| {d_{2N} } \right|\frac{{\sqrt 3 }}{2}\frac{{\Gamma \left( {\frac{{2N + 1}}{3}} \right)}}{{\left| \nu  \right|^{\frac{{2N + 1}}{3}} }} \begin{cases} \left|\sec \theta \right| & \; \text{ if } \; { - \frac{\pi }{2} < \theta  < 0 \; \text{ or } \; \pi  < \theta  < \frac{{3\pi }}{2}} \\ 1 & \; \text{ if } \; {0 \le \theta  \le \pi }, \end{cases}
\end{equation}
if $N\equiv 2 \mod 3$. Finally, the corresponding estimate for the case $N\equiv 0 \mod 3$ is
\begin{equation}\label{eq38}
\left| {R_N^{\left( H \right)} \left( \nu  \right)} \right| \le \left(\frac{2}{{3\pi }}\left| {d_{2N} } \right|\frac{{\sqrt 3 }}{2}\frac{{\Gamma \left( {\frac{{2N + 1}}{3}} \right)}}{{\left| \nu  \right|^{\frac{{2N + 1}}{3}} }} + \frac{2}{{3\pi }}\left| {d_{2N + 2} } \right|\frac{{\sqrt 3 }}{2}\frac{{\Gamma \left( {\frac{{2N + 3}}{3}} \right)}}{{\left| \nu  \right|^{\frac{{2N + 3}}{3}} }}\right) \begin{cases} \left|\sec \theta \right| & \; \text{ if } \; { - \frac{\pi }{2} < \theta  < 0 \; \text{ or } \; \pi  < \theta  < \frac{{3\pi }}{2}} \\ 1 & \; \text{ if } \; {0 \le \theta  \le \pi }. \end{cases}
\end{equation}
To obtain the analogous bounds for the remainder term $R_N^{\left( H \right)} \left( \nu e^{\pi i} \right)$ of the expansion \eqref{eq18} one has to simply replace $\theta$ by $\theta+\pi$ in each of the three estimates above.

Our bounds for $R_N^{\left( H \right)} \left( \nu \right)$ are unrealistic near the Stokes lines $\theta = -\frac{\pi }{2}, \frac{3\pi }{2}$ due to the presence of the factor $\sec \theta$. We shall derive better bounds for $R_N^{\left( H \right)} \left( \nu \right)$ near these lines using the method we applied in the previous case. Let $0 < \varphi  < \frac{\pi }{2}$ be an acute angle that may depend on $N$ and suppose that $ - \frac{\pi }{2} \le \theta  <  - \varphi$. We rotate the path of integration in \eqref{eq35} by $-\varphi$, and apply the inequality \eqref{eq47} to obtain
\begin{equation}\label{eq49}
\left| {R_N^{\left( H \right)} \left( \nu  \right)} \right| \le \frac{\sec \left( {\theta + \varphi} \right)}{{\sqrt 3 \pi \cos^{\frac{{2N + 3}}{3}} \varphi \left| \nu  \right|^{\frac{{2N + 3}}{3}} }}\int_0^{ + \infty } {t^{\frac{{2N}}{3}} e^{ - 2\pi t} \left| {H_{\frac{{it}}{{e^{i\varphi} \cos \varphi }}}^{\left( 1 \right)} \left( {\frac{{it}}{{e^{i\varphi} \cos \varphi }}} \right)} \right|dt} 
\end{equation}
for $ - \frac{\pi }{2} \le \theta  <  - \varphi$ with $N\equiv 1 \mod 3$. Using a continuity argument for the inequality \eqref{eq26}, yields
\begin{equation}\label{eq55}
\left| {H_{\frac{{it}}{{e^{i\varphi } \cos \varphi }}}^{\left( 1 \right)} \left( {\frac{{it}}{{e^{i\varphi } \cos \varphi }}} \right)} \right| \le \frac{1}{{\sqrt {\cos \varphi } }}iH_{it}^{\left( 1 \right)} \left( {it} \right) \le \frac{1}{{\cos ^{\frac{2}{3}} \varphi }}iH_{it}^{\left( 1 \right)} \left( {it} \right).
\end{equation}
The angle $\varphi = \arctan \left( {\left(\frac{N}{3}\right)^{ - \frac{1}{2}} } \right)$ minimises the function $\sec \left( { - \frac{\pi }{2} + \varphi } \right)\cos^{-\frac{N}{3}} \varphi$, and
\[
\frac{{\sec \left( {\theta  + \arctan \left( {\left(\frac{N}{3}\right)^{ - \frac{1}{2}} } \right)} \right)}}{{\cos^{\frac{N}{3}} \left( {\arctan \left( {\left(\frac{N}{3}\right)^{ - \frac{1}{2}} } \right)} \right)}} \le \frac{{\sec \left( { - \frac{\pi }{2} + \arctan \left( {\left(\frac{N}{3}\right)^{ - \frac{1}{2}} } \right)} \right)}}{{\cos^{\frac{N}{3}} \left( {\arctan \left( {\left(\frac{N}{3}\right)^{ - \frac{1}{2}} } \right)} \right)}} = 
\frac{1}{{\sqrt 3 }}\left( {1 + \frac{3}{N}} \right)^{\frac{{N + 3}}{6}} \sqrt{N} \le \sqrt {\frac{e}{3}\left( {N + \frac{3}{2}} \right)} 
\]
for all $ - \frac{\pi }{2} \le \theta  <  - \varphi = -\arctan \left( {\left(\frac{N}{3}\right)^{ - \frac{1}{2}} } \right)$ with $N \geq 1$. Substituting these estimations into \eqref{eq49} yields the bound
\begin{equation}\label{eq50}
\left| {R_N^{\left( H \right)} \left( \nu  \right)} \right| \le \sqrt {\frac{e}{3}\left( {2N + \frac{{13}}{2}} \right)} \frac{2}{{3\pi }}\left| {d_{2N + 2} } \right|\frac{{\sqrt 3 }}{2}\frac{{\Gamma \left( {\frac{{2N + 3}}{3}} \right)}}{{\left| \nu  \right|^{\frac{{2N + 3}}{3}} }},
\end{equation}
which is valid for $ - \frac{\pi}{2} \le \theta  <  - \arctan \left( {\left(\frac{2N+5}{3}\right)^{ - \frac{1}{2}} } \right)$ with $N\equiv 1 \mod 3$. A similar argument shows that this bound is also valid when $\pi + \arctan \left( {\left(\frac{2N+5}{3}\right)^{ - \frac{1}{2}} } \right) < \theta \le \frac{3\pi}{2}$. Since $\sqrt {\frac{e}{3}\left( {2N + \frac{{13}}{2}} \right)}  \ge \sqrt {\frac{{17}}{6}e}  \ge \left| {\sec \theta } \right|$ for $ - \arctan \left( {\left( {\frac{7}{3}} \right)^{ - \frac{1}{2}} } \right) \le \theta  < 0$ and $\pi  < \theta  \le \pi  + \arctan \left( {\left( {\frac{7}{3}} \right)^{ - \frac{1}{2}} } \right)$, the estimate \eqref{eq50} remains valid for $ - \frac{\pi}{2} \le \theta  <  0$ and $\pi < \theta  \le \frac{3\pi}{2}$.

The corresponding bound for the case $N\equiv 2 \mod 3$ takes the form
\begin{equation}\label{eq53}
\left| {R_N^{\left( H \right)} \left( \nu  \right)} \right| \le \sqrt {\frac{e}{3}\left( {2N + \frac{9}{2}} \right)} \frac{2}{3\pi}\left| {d_{2N} } \right|\frac{{\sqrt 3 }}{2}\frac{{\Gamma \left( {\frac{{2N + 1}}{3}} \right)}}{{\left| \nu  \right|^{\frac{{2N + 1}}{3}} }}
\end{equation}
for $ - \frac{\pi }{2} \le \theta  <  0$ and $\pi < \theta \le \frac{3\pi}{2}$.

For the third case, $N\equiv 0 \mod 3$, the upper bound for $\left| {R_N^{\left( H \right)} \left( \nu  \right)} \right|$ is the sum of the right-hand sides of \eqref{eq50} and \eqref{eq53}.

The remainder term $R_N^{\left( H \right)} \left( \nu e^{\pi i} \right)$ of the expansion \eqref{eq18} can be estimated similarly and we obtain the same upper bounds but with the modified conditions $0 < \theta  \le \frac{\pi}{2}$ and $ - \frac{3\pi }{2} \le \theta  <  -\pi$.

\subsection{Error bounds for the expansions of $J_\nu \left(\nu x\right)$ and $Y_\nu \left(\nu x\right)$}

\subsubsection{Case (i): $x>1$} Using formula \eqref{eq28}, together with the estimate \eqref{eq20} and the representation \eqref{eq29}, we obtain the error bound
\begin{equation}\label{eq31}
\left| {R_N^{\left( J \right)} \left( {\nu ,\beta } \right)} \right| \le \left( {\left|\cos \xi\right|\frac{{\left| {U_{2N} \left( {i\cot \beta } \right)} \right|}}{{\left| \nu  \right|^{2N} }} + \left|\sin \xi\right|\frac{{\left| {U_{2N + 1} \left( {i\cot \beta } \right)} \right|}}{{\left| \nu  \right|^{2N + 1} }}} \right) \begin{cases} \left|\csc\left(2\theta\right)\right| & \; \text{ if } \; \frac{\pi}{4} < \left|\theta\right| <\frac{\pi}{2} \\ 1 & \; \text{ if } \; \left|\theta\right| \leq \frac{\pi}{4}. \end{cases}
\end{equation}
The same bound holds for the remainder $R_N^{\left( Y \right)} \left( {\nu ,\beta } \right)$ with $\cos \xi$ and $\sin \xi$ being interchanged. When $\nu$ is real and positive, we can obtain more precise estimates. Indeed, as $0 < \frac{1}{{1 + \left( {t/\nu } \right)^2 }} < 1 $ for $t,\nu>0$, from \eqref{eq29}, \eqref{eq30} and \eqref{eq28} we find
\[
\left( { - 1} \right)^N R_N^{\left( J \right)} \left( {\nu ,\beta } \right) = \cos \xi \frac{{\left| {U_{2N} \left( {i\cot \beta } \right)} \right|}}{{\nu ^{2N} }}\Theta_1  + \sin \xi \frac{{\left| {U_{2N + 1} \left( {i\cot \beta } \right)} \right|}}{{\nu ^{2N + 1} }}\Theta_2 ,
\]
\[
\left( { - 1} \right)^N R_N^{\left( Y \right)} \left( {\nu ,\beta } \right) = \sin \xi \frac{{\left| {U_{2N} \left( {i\cot \beta } \right)} \right|}}{{\nu ^{2N} }}\Theta_1  - \cos \xi \frac{{\left| {U_{2N + 1} \left( {i\cot \beta } \right)} \right|}}{{\nu ^{2N + 1} }}\Theta_2 ,
\]
where $0 < \Theta_i < 1$ ($i=1,2$) is an appropriate number depending on $\nu,\beta$ and $N$.

The error bound \eqref{eq31} becomes singular as $\theta \to \pm\frac{\pi }{2}$, and hence unrealistic near the Stokes lines. A better bound for $R_N^{\left( J \right)} \left( {\nu ,\beta } \right)$ near $\theta = \pm\frac{\pi }{2}$ can be derived using the method we applied in the case of the Hankel functions. Let $0 < \varphi  < \frac{\pi }{2}$ be an acute angle that may depend on $N$. First, suppose that $\frac{\pi }{4} + \varphi  < \theta  \le \frac{\pi }{2}$. We rotate the path of integration in \eqref{eq29} through the angle $\varphi$, and employ the estimates \eqref{eq20}, \eqref{eq32} and \eqref{eq26} (with $\varphi$ in place of $-\varphi$) to obtain
\[
\left| {R_N^{\left( J \right)} \left( {\nu ,\beta } \right)} \right| \le \frac{{\csc \left( {2\left( {\theta  - \varphi } \right)} \right)}}{{\cos ^{2N + 1} \varphi }}\left| {\cos \xi } \right|\frac{{\left| {U_{2N} \left( {i\cot \beta } \right)} \right|}}{{\left| \nu  \right|^{2N} }} + \frac{{\csc \left( {2\left( {\theta  - \varphi } \right)} \right)}}{{\cos ^{2N + 2} \varphi }}\left| {\sin \xi } \right|\frac{{\left| {U_{2N + 1} \left( {i\cot \beta } \right)} \right|}}{{\left| \nu  \right|^{2N + 1} }}
\]
for $\frac{\pi }{4} + \varphi  < \theta  \le \frac{\pi }{2}$ and $N\geq 0$. The value $\varphi = \arctan \left( {\left(N+2\right)^{ - \frac{1}{2}} } \right)$ minimises the function $\csc \left( { 2\left(\frac{\pi }{2}-\varphi\right) } \right)\cos^{ - N-1} \varphi$, and
\begin{align*}
 \frac{{\csc \left( {2\left( {\theta - \arctan \left( {\left( {N + 2} \right)^{ - \frac{1}{2}} } \right) } \right)} \right)}}{{\cos ^{N+1} \left( {\arctan \left( {\left( {N + 2} \right)^{ - \frac{1}{2}} } \right)} \right)}} & \le \frac{{\csc \left( {2\left( {\frac{\pi }{2}- \arctan \left( {\left( {N + 2} \right)^{ - \frac{1}{2}} } \right)} \right)} \right)}}{{\cos ^{N+1} \left( {\arctan \left( {\left( {N + 2} \right)^{ - \frac{1}{2}} } \right)} \right)}}\\ & = \frac{1}{2}\left( {1 + \frac{1}{{N + 2}}} \right)^{\frac{N+3}{2}} \sqrt{N+2}  \le \frac{1}{2}\sqrt {e\left( {N + \frac{5}{2}} \right)}  
\end{align*}
for all $\frac{\pi }{4} + \varphi = \frac{\pi }{4} + \arctan \left( {\left(N+2\right)^{ - \frac{1}{2}} } \right) < \theta  \le \frac{\pi }{2}$ with $N \geq 0$. Therefore, we have
\begin{equation}\label{eq45}
\left| {R_N^{\left( J \right)} \left( {\nu ,\beta } \right)} \right| \le \frac{1}{2}\sqrt {e\left( {2N + \frac{5}{2}} \right)} \left| {\cos \xi } \right|\frac{{\left| {U_{2N} \left( {i\cot \beta } \right)} \right|}}{{\left| \nu  \right|^{2N} }} + \frac{1}{2}\sqrt {e\left( {2N + \frac{7}{2}} \right)} \left| {\sin \xi } \right|\frac{{\left| {U_{2N + 1} \left( {i\cot \beta } \right)} \right|}}{{\left| \nu  \right|^{2N + 1} }}
\end{equation}
for $\frac{\pi }{4} + \arctan \left( {\left(2N+2\right)^{ - \frac{1}{2}} } \right) < \theta  \le \frac{\pi }{2}$ and $N\geq 0$. Since $\left| {R_N^{\left( J \right)} \left( {\bar \nu ,\beta } \right)} \right| = \left| {\overline {R_N^{\left( J \right)} \left( {\nu ,\beta } \right)} } \right| = \left| {R_N^{\left( J \right)} \left( {\nu ,\beta } \right)} \right|$, this bound also holds when $ - \frac{\pi }{2} \le \theta  <  - \frac{\pi }{4} - \arctan \left( {\left( {2N + 2} \right)^{ - \frac{1}{2}} } \right)$. In the ranges $\frac{\pi }{4} < \theta  \le \frac{\pi }{4} + \arctan \left( {\frac{1}{2}} \right)$ and $- \frac{\pi }{4} - \arctan \left( {\frac{1}{2}} \right) \le \theta  <  - \frac{\pi }{4}$, it holds that $\frac{1}{2}\sqrt {\frac{9}{2}e}  \ge \left| {\csc \left( {2\theta } \right)} \right|$, whence the estimate \eqref{eq45} is valid in the wider sectors $\frac{\pi }{4} < \left| \theta  \right| \le \frac{\pi }{2}$ as long as $N \geq 1$.

The same bound holds for the remainder $R_N^{\left( Y \right)} \left( {\nu ,\beta } \right)$ with $\cos \xi$ and $\sin \xi$ being interchanged.

\subsubsection{Case (ii): $x=1$} First, we consider the error bounds for the expansion of $J_\nu\left(\nu\right)$. Formula \eqref{eq42} can be simplified to
\begin{align}
R_N^{\left( J \right)} \left( \nu  \right) & = \frac{{\left( { - 1} \right)^N }}{{2\sqrt 3 \pi \nu ^{\frac{{2N + 1}}{3}} }}\int_0^{ + \infty } {t^{\frac{{2N - 2}}{3}} e^{ - 2\pi t} \frac{{1 - \left( {t/\nu } \right)^{\frac{4}{3}} }}{{1 + \left( {t/\nu } \right)^2 }}iH_{it}^{\left( 1 \right)} \left( {it} \right)dt} , \label{eq54} \\
R_N^{\left( J \right)} \left( \nu  \right) & = \frac{{\left( { - 1} \right)^N }}{{2\sqrt 3 \pi \nu ^{\frac{{2N + 3}}{3}} }}\int_0^{ + \infty } {t^{\frac{{2N}}{3}} e^{ - 2\pi t} \frac{{1 + \left( {t/\nu } \right)^{\frac{2}{3}} }}{{1 + \left( {t/\nu } \right)^2 }}iH_{it}^{\left( 1 \right)} \left( {it} \right)dt} , \nonumber \\
R_N^{\left( J \right)} \left( \nu  \right) & = \frac{{\left( { - 1} \right)^{N + 1} }}{{2\sqrt 3 \pi \nu ^{\frac{{2N + 1}}{3}} }}\int_0^{ + \infty } {t^{\frac{{2N - 2}}{3}} e^{ - 2\pi t} \frac{{1 + \left( {t/\nu } \right)^{\frac{2}{3}} }}{{1 + \left( {t/\nu } \right)^2 }}iH_{it}^{\left( 1 \right)} \left( {it} \right)dt} \nonumber
\end{align}
according to whether $N\equiv 0 \mod 3$, $N\equiv 1 \mod 3$ or $N\equiv 2 \mod 3$, respectively. Hence, by \eqref{eq20} and \eqref{eq37} we get the bounds
\begin{equation}\label{eq56}
\left|R_N^{\left( J \right)} \left( \nu  \right)\right| \leq \left( {\frac{1}{{3\pi }}\left| {d_{2N} } \right|\frac{{\sqrt 3 }}{2}\frac{{\Gamma \left( {\frac{{2N + 1}}{3}} \right)}}{{\left| \nu  \right|^{\frac{{2N + 1}}{3}} }} + \frac{1}{{3\pi }}\left| {d_{2N + 4} } \right|\frac{{\sqrt 3 }}{2}\frac{{\Gamma \left( {\frac{{2N + 5}}{3}} \right)}}{{\left| \nu  \right|^{\frac{{2N + 5}}{3}} }}} \right) \begin{cases} \left|\csc\left(2\theta\right)\right| & \; \text{ if } \; \frac{\pi}{4} < \left|\theta\right| <\frac{\pi}{2} \\ 1 & \; \text{ if } \; \left|\theta\right| \leq \frac{\pi}{4}, \end{cases}
\end{equation}
when $N\equiv 0 \mod 3$;
\[
\left|R_N^{\left( J \right)} \left( \nu  \right)\right| \leq \left( {\frac{1}{{3\pi }}\left| {d_{2N + 2} } \right|\frac{{\sqrt 3 }}{2}\frac{{\Gamma \left( {\frac{{2N + 3}}{3}} \right)}}{{\left| \nu  \right|^{\frac{{2N + 3}}{3}} }} + \frac{1}{{3\pi }}\left| {d_{2N + 4} } \right|\frac{{\sqrt 3 }}{2}\frac{{\Gamma \left( {\frac{{2N + 5}}{3}} \right)}}{{\left| \nu  \right|^{\frac{{2N + 5}}{3}} }}} \right) \begin{cases} \left|\csc\left(2\theta\right)\right| & \; \text{ if } \; \frac{\pi}{4} < \left|\theta\right| <\frac{\pi}{2} \\ 1 & \; \text{ if } \; \left|\theta\right| \leq \frac{\pi}{4}, \end{cases}
\]
when $N\equiv 1 \mod 3$;
\[
\left|R_N^{\left( J \right)} \left( \nu  \right)\right| \leq \left( {\frac{1}{{3\pi }}\left| {d_{2N} } \right|\frac{{\sqrt 3 }}{2}\frac{{\Gamma \left( {\frac{{2N + 1}}{3}} \right)}}{{\left| \nu  \right|^{\frac{{2N + 1}}{3}} }} + \frac{1}{{3\pi }}\left| {d_{2N + 2} } \right|\frac{{\sqrt 3 }}{2}\frac{{\Gamma \left( {\frac{{2N + 3}}{3}} \right)}}{{\left| \nu  \right|^{\frac{{2N + 3}}{3}} }}} \right) \begin{cases} \left|\csc\left(2\theta\right)\right| & \; \text{ if } \; \frac{\pi}{4} < \left|\theta\right| <\frac{\pi}{2} \\ 1 & \; \text{ if } \; \left|\theta\right| \leq \frac{\pi}{4}, \end{cases}
\]
when $N\equiv 2 \mod 3$. Again, when $\nu$ is real and positive, we can deduce better estimates:
\[
\left( { - 1} \right)^N R_N^{\left( J \right)} \left( \nu  \right) = \frac{1}{{3\pi }}\left| {d_{2N} } \right|\frac{{\sqrt 3 }}{2}\frac{{\Gamma \left( {\frac{{2N + 1}}{3}} \right)}}{{\nu ^{\frac{{2N + 1}}{3}} }}\Xi _1  - \frac{1}{{3\pi }}\left| {d_{2N + 4} } \right|\frac{{\sqrt 3 }}{2}\frac{{\Gamma \left( {\frac{{2N + 5}}{3}} \right)}}{{\nu ^{\frac{{2N + 5}}{3}} }}\Xi _3 ,
\]
when $N\equiv 0 \mod 3$;
\[
\left( { - 1} \right)^N R_N^{\left( J \right)} \left( \nu  \right) = \frac{1}{{3\pi }}\left| {d_{2N + 2} } \right|\frac{{\sqrt 3 }}{2}\frac{{\Gamma \left( {\frac{{2N + 3}}{3}} \right)}}{{\nu ^{\frac{{2N + 3}}{3}} }}\Xi _2  + \frac{1}{{3\pi }}\left| {d_{2N + 4} } \right|\frac{{\sqrt 3 }}{2}\frac{{\Gamma \left( {\frac{{2N + 5}}{3}} \right)}}{{\nu ^{\frac{{2N + 5}}{3}} }}\Xi _3 ,
\]
when $N\equiv 1 \mod 3$;
\[
\left( { - 1} \right)^{N + 1} R_N^{\left( J \right)} \left( \nu  \right) = \frac{1}{{3\pi }}\left| {d_{2N} } \right|\frac{{\sqrt 3 }}{2}\frac{{\Gamma \left( {\frac{{2N + 1}}{3}} \right)}}{{\nu ^{\frac{{2N + 1}}{3}} }}\Xi _1  + \frac{1}{{3\pi }}\left| {d_{2N + 2} } \right|\frac{{\sqrt 3 }}{2}\frac{{\Gamma \left( {\frac{{2N + 3}}{3}} \right)}}{{\nu ^{\frac{{2N + 3}}{3}} }}\Xi _2 ,
\]
when $N\equiv 2 \mod 3$. Here $0 < \Xi_i < 1$ ($i=1,2,3$) is an appropriate number depending on $\nu$ and $N$. Substituting $N=0$ into the first one yields
\[
J_\nu  \left( \nu  \right) < \frac{1}{{3\pi }}\left| {d_0 } \right|\frac{{\sqrt 3 }}{2}\frac{{\Gamma \left( {\frac{1}{3}} \right)}}{{\nu ^{\frac{1}{3}} }} = \frac{{\Gamma \left( {\frac{1}{3}} \right)}}{{2^{\frac{2}{3}} 3^{\frac{1}{6}} \pi \nu ^{\frac{1}{3}} }} \; \text{ for } \; \nu >0.
\]
This upper bound was established by Watson \cite[pp. 258--259]{Watson} using a method different from ours. We remark that it can be shown easily that
\[
\left| {\frac{{1 - \left( {t/\nu } \right)^{\frac{4}{3}} }}{{1 + \left( {t/\nu } \right)^2 }}} \right| \le 1,\; \text{ for } \left| \theta  \right| \le \frac{\pi }{4},
\]
whence
\[
\left| {R_N^{\left( J \right)} \left( \nu  \right)} \right| \le \frac{1}{{3\pi }}\left| {d_{2N} } \right|\frac{{\sqrt 3 }}{2}\frac{{\Gamma \left( {\frac{{2N + 1}}{3}} \right)}}{{\left| \nu  \right|^{\frac{{2N + 1}}{3}} }}
\]
holds when $N\equiv 0 \mod 3$ and $\left| \theta  \right| \le \frac{\pi}{4}$.

As before, it is possible to derive better bounds near the Stokes lines $\theta = \pm\frac{\pi }{2}$. Let $0 < \varphi  < \frac{\pi}{2}$ be an acute angle that may depend on $N$. First, suppose that $\frac{\pi }{4} + \varphi  < \theta  \le \frac{\pi }{2}$. We rotate the path of integration in \eqref{eq54} by $\varphi$, and apply the inequality \eqref{eq20} to obtain
\begin{align*}
\left| {R_N^{\left( J \right)} \left( \nu  \right)} \right|  \le \; & \frac{{\csc \left( {2\left( {\theta  - \varphi } \right)} \right)}}{{2\sqrt 3 \pi \cos ^{\frac{{2N + 1}}{3}} \varphi \left| \nu  \right|^{\frac{{2N + 1}}{3}} }}\int_0^{ + \infty } {t^{\frac{{2N - 2}}{3}} e^{ - 2\pi t} \left| {H_{\frac{{ite^{i\varphi } }}{{\cos \varphi }}}^{\left( 1 \right)} \left( {\frac{{ite^{i\varphi } }}{{\cos \varphi }}} \right)} \right|dt}\\ &   + \frac{{\csc \left( {2\left( {\theta  - \varphi } \right)} \right)}}{{2\sqrt 3 \pi \cos ^{\frac{{2N + 5}}{3}} \varphi \left| \nu  \right|^{\frac{{2N + 5}}{3}} }}\int_0^{ + \infty } {t^{\frac{{2N + 2}}{3}} e^{ - 2\pi t} \left| {H_{\frac{{ite^{i\varphi } }}{{\cos \varphi }}}^{\left( 1 \right)} \left( {\frac{{ite^{i\varphi } }}{{\cos \varphi }}} \right)} \right|dt} 
\end{align*}
for $\frac{\pi }{4} + \varphi  < \theta  \le \frac{\pi }{2}$ and $N\equiv 0 \mod 3$, $N\geq 3$. The modulus of the Hankel function can be estimated via \eqref{eq55} by replacing $-\varphi$ with $\varphi$. The angle $\varphi  = \arctan \left( {\left( {\frac{1}{3}N + 1} \right)^{ - \frac{1}{2}} } \right)$ minimises the function $\csc \left( { 2\left(\frac{\pi }{2}-\varphi\right) } \right)\cos^{ - \frac{N}{3}} \varphi$, and
\begin{align*}
\frac{{\csc \left( {2\left( {\theta  - \arctan \left( {\left( {\frac{1}{3}N + 1} \right)^{ - \frac{1}{2}} } \right)} \right)} \right)}}{{\cos ^{\frac{N}{3}} \left( {\arctan \left( {\left( {\frac{1}{3}N + 1} \right)^{ - \frac{1}{2}} } \right)} \right)}} & \le \frac{{\csc \left( {2\left( {\frac{\pi }{2} - \arctan \left( {\left( {\frac{1}{3}N + 1} \right)^{ - \frac{1}{2}} } \right)} \right)} \right)}}{{\cos ^{\frac{N}{3}} \left( {\arctan \left( {\left( {\frac{1}{3}N + 1} \right)^{ - \frac{1}{2}} } \right)} \right)}}\\ & = \frac{{\sqrt 3 }}{6}\left( {1 + \frac{3}{{N + 3}}} \right)^{\frac{N}{6} + 1} \sqrt {N + 3}  \le \frac{1}{2}\sqrt {\frac{e}{3}\left( {N + \frac{9}{2}} \right)}
\end{align*}
for all $\frac{\pi}{4} + \varphi = \frac{\pi}{4} + \arctan \left( {\left( {\frac{1}{3}N + 1} \right)^{ - \frac{1}{2}} } \right) < \theta  \le \frac{\pi}{2}$ with $N \geq 3$. It follows that
\[
\left| {R_N^{\left( J \right)} \left( \nu  \right)} \right| \le \frac{1}{2}\sqrt {\frac{e}{3}\left( {2N + \frac{{15}}{2}} \right)} \frac{1}{{3\pi }}\left| {d_{2N} } \right|\frac{{\sqrt 3 }}{2}\frac{{\Gamma \left( {\frac{{2N + 1}}{3}} \right)}}{{\left| \nu  \right|^{\frac{{2N + 1}}{3}} }} + \frac{1}{2}\sqrt {\frac{e}{3}\left( {2N + \frac{{23}}{2}} \right)} \frac{1}{{3\pi }}\left| {d_{2N + 4} } \right|\frac{{\sqrt 3 }}{2}\frac{{\Gamma \left( {\frac{{2N + 5}}{3}} \right)}}{{\left| \nu  \right|^{\frac{{2N + 5}}{3}} }}
\]
for $\frac{\pi}{4} + \arctan \left( {\left( \frac{2N + 6}{3} \right)^{ - \frac{1}{2}} } \right) < \theta  \le \frac{\pi}{2}$ and $N\equiv 0 \mod 3$. As $\left| {R_N^{\left( J \right)} \left( {\bar \nu } \right)} \right| = \left| {\overline {R_N^{\left( J \right)} \left( \nu  \right)} } \right| = \left| {R_N^{\left( J \right)} \left( \nu  \right)} \right|$, this bound also holds when $-\frac{\pi}{2} \le \theta < -\frac{\pi}{4} - \arctan \left( {\left( \frac{2N + 6}{3} \right)^{ - \frac{1}{2}} } \right)$. Since $\frac{1}{2}\sqrt {\frac{e}{3}\left( {2N + \frac{{15}}{2}} \right)}  \ge \frac{1}{2}\sqrt {\frac{9}{2}e}  \ge \left| {\csc \left( {2\theta } \right)} \right|$ for $\frac{\pi }{4} < \theta  \le \frac{\pi }{4} + \arctan \left( {\left( {\frac{{12}}{3}} \right)^{ - \frac{1}{2}} } \right)$ and $- \frac{\pi }{4} - \arctan \left( {\left( {\frac{{12}}{3}} \right)^{ - \frac{1}{2}} } \right) < \theta  \le  - \frac{\pi }{4}$, the estimate \eqref{eq50} remains valid for $\frac{\pi }{4} < \left| \theta  \right| \le \frac{\pi }{2}$ as long as $N \geq 3$.

Similarly, we have
\[
\left| {R_N^{\left( J \right)} \left( \nu  \right)} \right| \le \frac{1}{2}\sqrt {\frac{e}{3}\left( {2N + \frac{{19}}{2}} \right)} \frac{1}{{3\pi }}\left| {d_{2N + 2} } \right|\frac{{\sqrt 3 }}{2}\frac{{\Gamma \left( {\frac{{2N + 3}}{3}} \right)}}{{\left| \nu  \right|^{\frac{{2N + 3}}{3}} }} + \frac{1}{2}\sqrt {\frac{e}{3}\left( {2N + \frac{{23}}{2}} \right)} \frac{1}{{3\pi }}\left| {d_{2N + 4} } \right|\frac{{\sqrt 3 }}{2}\frac{{\Gamma \left( {\frac{{2N + 5}}{3}} \right)}}{{\left| \nu  \right|^{\frac{{2N + 5}}{3}} }}
\]
for $\frac{\pi}{4} + \arctan \left( {\left( {\frac{2N+8}{3}} \right)^{ - \frac{1}{2}} } \right)  \le \theta  \le \frac{\pi}{2}$, $-\frac{\pi}{2} \le \theta \le -\frac{\pi}{4} - \arctan \left( {\left( {\frac{2N+8}{3}} \right)^{ - \frac{1}{2}} } \right)$ and $N\equiv 1 \mod 3$. If $N \geq 4$, this bound becomes valid in the wider range $\frac{\pi }{4} < \left| \theta  \right| \le \frac{\pi }{2}$.

Finally, it holds that
\[
\left| {R_N^{\left( J \right)} \left( \nu  \right)} \right| \le \frac{1}{2}\sqrt {\frac{e}{3}\left( {2N + \frac{{15}}{2}} \right)} \frac{1}{{3\pi }}\left| {d_{2N} } \right|\frac{{\sqrt 3 }}{2}\frac{{\Gamma \left( {\frac{{2N + 1}}{3}} \right)}}{{\left| \nu  \right|^{\frac{{2N + 1}}{3}} }} + \frac{1}{2}\sqrt {\frac{e}{3}\left( {2N + \frac{{19}}{2}} \right)} \frac{1}{{3\pi }}\left| {d_{2N + 2} } \right|\frac{{\sqrt 3 }}{2}\frac{{\Gamma \left( {\frac{{2N + 3}}{3}} \right)}}{{\left| \nu  \right|^{\frac{{2N + 3}}{3}} }}
\]
for $\frac{\pi}{4} + \arctan \left( {\left( \frac{2N + 6}{3} \right)^{ - \frac{1}{2}} } \right)  \le \theta  \le \frac{\pi}{2}$, $-\frac{\pi}{2} \le \theta \le -\frac{\pi}{4} - \arctan \left( {\left( \frac{2N + 6}{3} \right)^{ - \frac{1}{2}} } \right)$ and $N\equiv 2 \mod 3$. For $N \geq 5$, this bound is valid in the wider sector $\frac{\pi}{4} < \left| \theta  \right| \le \frac{\pi}{2}$.

Lastly, we consider the error bounds for the expansion of $Y_\nu\left(\nu\right)$. The formula \eqref{eq43} can be simplified to
\begin{align*}
R_N^{\left( Y \right)} \left( \nu  \right) & = \frac{{\left( { - 1} \right)^{N + 1} }}{{2\pi \nu ^{\frac{{2N + 1}}{3}} }}\int_0^{ + \infty } {t^{\frac{{2N - 2}}{3}} e^{ - 2\pi t} \frac{{1 + \left( {t/\nu } \right)^{\frac{4}{3}} }}{{1 + \left( {t/\nu } \right)^2 }}iH_{it}^{\left( 1 \right)} \left( {it} \right)dt} ,\\
R_N^{\left( Y \right)} \left( \nu  \right) & = \frac{{\left( { - 1} \right)^N }}{{2\pi \nu ^{\frac{{2N + 3}}{3}} }}\int_0^{ + \infty } {t^{\frac{{2N}}{3}} e^{ - 2\pi t} \frac{{1 - \left( {t/\nu } \right)^{\frac{2}{3}} }}{{1 + \left( {t/\nu } \right)^2 }}iH_{it}^{\left( 1 \right)} \left( {it} \right)dt} ,\\
R_N^{\left( Y \right)} \left( \nu  \right) & = \frac{{\left( { - 1} \right)^{N + 1} }}{{2\pi \nu ^{\frac{{2N + 1}}{3}} }}\int_0^{ + \infty } {t^{\frac{{2N - 2}}{3}} e^{ - 2\pi t} \frac{{1 - \left( {t/\nu } \right)^{\frac{2}{3}} }}{{1 + \left( {t/\nu } \right)^2 }}iH_{it}^{\left( 1 \right)} \left( {it} \right)dt}
\end{align*}
according to whether $N\equiv 0 \mod 3$, $N\equiv 1 \mod 3$ or $N\equiv 2 \mod 3$, respectively. We show in Appendix \ref{appendixb} that
\begin{equation}\label{eq48}
\left| {\frac{{1 - \left( {t/\nu } \right)^{\frac{2}{3}} }}{{1 + \left( {t/\nu } \right)^2 }}} \right| \le \begin{cases} \left|\csc\left(2\theta\right)\right| & \; \text{ if } \; \frac{\pi}{4} < \left|\theta\right| <\frac{\pi}{2} \\ 1 & \; \text{ if } \; \left|\theta\right| \leq \frac{\pi}{4}. \end{cases}
\end{equation}
Applying formula \eqref{eq37} together with the inequality \eqref{eq20} in the first and with the inequality \eqref{eq48} in the second and the third case, we deduce
\[
\left| {R_N^{\left( Y \right)} \left( \nu  \right)} \right| \le \left( \frac{2}{{3\pi }}\left| {d_{2N} } \right|\frac{3}{4}\frac{{\Gamma \left( {\frac{{2N + 1}}{3}} \right)}}{{\left| \nu  \right|^{\frac{{2N + 1}}{3}} }} + \frac{2}{{3\pi }}\left| {d_{2N + 4} } \right|\frac{3}{4}\frac{{\Gamma \left( {\frac{{2N + 5}}{3}} \right)}}{{\left| \nu  \right|^{\frac{{2N + 5}}{3}} }}\right) \begin{cases} \csc\left(2\theta\right) & \; \text{ if } \; \frac{\pi}{4} < \left|\theta\right| <\frac{\pi}{2} \\ 1 & \; \text{ if } \; \left|\theta\right| \leq \frac{\pi}{4}, \end{cases}
\]
when $N\equiv 0 \mod 3$;
\[
\left| {R_N^{\left( Y \right)} \left( \nu  \right)} \right| \le \frac{2}{{3\pi }}\left| {d_{2N + 2} } \right|\frac{3}{4}\frac{{\Gamma \left( {\frac{{2N + 3}}{3}} \right)}}{{\left| \nu  \right|^{\frac{{2N + 3}}{3}} }}  \begin{cases} \csc\left(2\theta\right) & \; \text{ if } \; \frac{\pi}{4} < \left|\theta\right| <\frac{\pi}{2} \\ 1 & \; \text{ if } \; \left|\theta\right| \leq \frac{\pi}{4}, \end{cases}
\]
when $N\equiv 1 \mod 3$;
\[
\left| {R_N^{\left( Y \right)} \left( \nu  \right)} \right| \le \frac{2}{{3\pi }}\left| {d_{2N} } \right|\frac{3}{4}\frac{{\Gamma \left( {\frac{{2N + 1}}{3}} \right)}}{{\left| \nu  \right|^{\frac{{2N + 1}}{3}} }} \begin{cases} \csc\left(2\theta\right) & \; \text{ if } \; \frac{\pi}{4} < \left|\theta\right| <\frac{\pi}{2} \\ 1 & \; \text{ if } \; \left|\theta\right| \leq \frac{\pi}{4}, \end{cases}
\]
when $N\equiv 2 \mod 3$, respectively. When $\nu$ is real and positive, we have the following better estimates:
\[
\left( { - 1} \right)^{N + 1} R_N^{\left( Y \right)} \left( \nu  \right) = \frac{2}{{3\pi }}\left| {d_{2N} } \right|\frac{3}{4}\frac{{\Gamma \left( {\frac{{2N + 1}}{3}} \right)}}{{\nu ^{\frac{{2N + 1}}{3}} }}\Phi _1  + \frac{2}{{3\pi }}\left| {d_{2N + 4} } \right|\frac{3}{4}\frac{{\Gamma \left( {\frac{{2N + 5}}{3}} \right)}}{{\nu ^{\frac{{2N + 5}}{3}} }}\Phi _3 ,
\]
when $N\equiv 0 \mod 3$;
\[
\left( { - 1} \right)^N R_N^{\left( Y \right)} \left( \nu  \right) = \frac{2}{{3\pi }}\left| {d_{2N + 2} } \right|\frac{3}{4}\frac{{\Gamma \left( {\frac{{2N + 3}}{3}} \right)}}{{\nu ^{\frac{{2N + 3}}{3}} }}\Phi _2  - \frac{2}{{3\pi }}\left| {d_{2N + 4} } \right|\frac{3}{4}\frac{{\Gamma \left( {\frac{{2N + 5}}{3}} \right)}}{{\nu ^{\frac{{2N + 5}}{3}} }}\Phi _3 ,
\]
when $N\equiv 1 \mod 3$;
\[
\left( { - 1} \right)^{N + 1} R_N^{\left( Y \right)} \left( \nu  \right) = \frac{2}{{3\pi }}\left| {d_{2N} } \right|\frac{3}{4}\frac{{\Gamma \left( {\frac{{2N + 1}}{3}} \right)}}{{\nu ^{\frac{{2N + 1}}{3}} }}\Phi _1  - \frac{2}{{3\pi }}\left| {d_{2N + 2} } \right|\frac{3}{4}\frac{{\Gamma \left( {\frac{{2N + 3}}{3}} \right)}}{{\nu ^{\frac{{2N + 3}}{3}} }}\Phi _2 ,
\]
when $N\equiv 2 \mod 3$. Here $0 < \Phi_i < 1$ ($i=1,2,3$) is an appropriate number depending on $\nu$ and $N$.

To bound $R_N^{\left( Y \right)} \left( \nu  \right)$ near the Stokes lines, we use the same argument as for $R_N^{\left( J \right)} \left( \nu  \right)$, and find that
\[
\left| {R_N^{\left( Y \right)} \left( \nu  \right)} \right| \le \frac{1}{2}\sqrt {\frac{e}{3}\left( {2N + \frac{{15}}{2}} \right)} \frac{2}{{3\pi }}\left| {d_{2N} } \right|\frac{3}{4}\frac{{\Gamma \left( {\frac{{2N + 1}}{3}} \right)}}{{\left| \nu  \right|^{\frac{{2N + 1}}{3}} }} + \frac{1}{2}\sqrt {\frac{e}{3}\left( {2N + \frac{{23}}{2}} \right)} \frac{2}{{3\pi }}\left| {d_{2N + 4} } \right|\frac{3}{4}\frac{{\Gamma \left( {\frac{{2N + 5}}{3}} \right)}}{{\left| \nu  \right|^{\frac{{2N + 5}}{3}} }}
\]
holds for $\frac{\pi}{4} + \arctan \left( {\left( \frac{2N+6}{3} \right)^{ - \frac{1}{2}} } \right)  \le \theta  \le \frac{\pi}{2}$, $-\frac{\pi}{2} \le \theta \le -\frac{\pi}{4} - \arctan \left( {\left( \frac{2N+6}{3} \right)^{ - \frac{1}{2}} } \right)$ and $N\equiv 0 \mod 3$. Similarly,
\[
\left| {R_N^{\left( Y \right)} \left( \nu  \right)} \right| \le \frac{1}{2}\sqrt {\frac{e}{3}\left( {2N + \frac{{19}}{2}} \right)} \frac{2}{{3\pi }}\left| {d_{2N + 2} } \right|\frac{3}{4}\frac{{\Gamma \left( {\frac{{2N + 3}}{3}} \right)}}{{\left| \nu  \right|^{\frac{{2N + 3}}{3}} }}
\]
holds for $\frac{\pi}{4} + \arctan \left( {\left( {\frac{2N+8}{3}} \right)^{ - \frac{1}{2}} } \right)  \le \theta  \le \frac{\pi}{2}$, $-\frac{\pi}{2} \le \theta \le -\frac{\pi}{4} - \arctan \left( {\left( {\frac{2N+8}{3}} \right)^{ - \frac{1}{2}} } \right)$ and $N\equiv 1 \mod 3$. Finally, we have that
\[
\left| {R_N^{\left( Y \right)} \left( \nu  \right)} \right| \le \frac{1}{2}\sqrt {\frac{e}{3}\left( {2N + \frac{{15}}{2}} \right)} \frac{2}{{3\pi }}\left| {d_{2N} } \right|\frac{3}{4}\frac{{\Gamma \left( {\frac{{2N + 1}}{3}} \right)}}{{\left| \nu  \right|^{\frac{{2N + 1}}{3}} }}
\]
for $\frac{\pi}{4} + \arctan \left( {\left( \frac{2N+6}{3} \right)^{ - \frac{1}{2}} } \right)  \le \theta  \le \frac{\pi}{2}$, $-\frac{\pi}{2} \le \theta \le -\frac{\pi}{4} - \arctan \left( {\left( \frac{2N+6}{3} \right)^{ - \frac{1}{2}} } \right)$ and $N\equiv 2 \mod 3$.

If $N \geq 4$, these bounds become valid in the wider range $\frac{\pi }{4} < \left| \theta  \right| \le \frac{\pi }{2}$.

\section{Asymptotics for the late coefficients}\label{section4}

In this section, we investigate the asymptotic nature of the coefficients $U_n\left(i \cot \beta\right)$ and $d_{2n}$ as $n \to +\infty$. First, we consider the coefficients $U_n\left(i \cot \beta\right)$. For our purposes, the most appropriate representation of these coefficients is the second integral formula in \eqref{eq13}. Upon replacing $H_{it}^{\left(1\right)}\left(it \sec \beta\right)$ by its representation \eqref{eq21} in this integral, we obtain
\begin{gather}\label{eq60}
\begin{split}
U_n \left( {i\cot \beta } \right) = \; & \frac{{\left(-1\right)^n \Gamma \left( n \right)}}{{2\pi \left( {2i\left( {\tan \beta  - \beta } \right)} \right)^n }}\left( {\sum\limits_{m = 0}^{M - 1} {\left( {2i\left( {\tan \beta  - \beta } \right)} \right)^m U_m \left( {i\cot \beta } \right)\frac{{\Gamma \left( {n - m} \right)}}{{\Gamma \left( n \right)}}}  + A_M^{\left( 1 \right)} \left( {n,\beta } \right)} \right) \\ & + \frac{{\left(-1\right)^n  \Gamma\left( n \right)}}{{2\pi \left( {2i\left( {\tan \beta  - \beta  + \pi } \right)} \right)^n }}\left( {\sum\limits_{m = 0}^{M - 1} {\left( {2i\left( {\tan \beta  - \beta  + \pi } \right)} \right)^m U_m \left( {i\cot \beta } \right)\frac{{\Gamma \left( {n - m} \right)}}{{\Gamma \left( n \right)}}}  + A_M^{\left( 2 \right)} \left( {n,\beta } \right)} \right),
\end{split}
\end{gather}
for any fixed $0 \le M \le n - 1$, provided that $n\geq 1$. Note that the second sum is the contribution of the exponential quantity $e^{ - 2\pi t}$ in \eqref{eq13}. The remainder terms $A_M^{\left( 1 \right)} \left( {n,\beta } \right)$ and $A_M^{\left( 2 \right)} \left( {n,\beta } \right)$ are given by the integral formulas
\[
A_M^{\left( 1 \right)} \left( {n,\beta } \right) = \frac{{\left( {2\left( {\tan \beta  - \beta } \right)} \right)^n }}{{\Gamma\left( n \right)}}\int_0^{+\infty} {t^{n - 1} e^{ - 2t\left( {\tan \beta  - \beta } \right)} R_M^{\left( H \right)} \left( {it,\beta } \right)dt} 
\]
and
\[
A_M^{\left( 2 \right)} \left( {n,\beta } \right) = \frac{{\left( {2\left( {\tan \beta  - \beta  + \pi } \right)} \right)^n }}{{\Gamma\left( n\right)}}\int_0^{+\infty} {t^{n - 1} e^{ - 2t\left( {\tan \beta  - \beta  + \pi } \right)} R_M^{\left( H \right)} \left( {it,\beta } \right)dt} ,
\]
respectively. To bound these error terms, we apply the estimate \eqref{eq24} to find
\begin{equation}\label{eq59}
\left| {A_M^{\left( 1 \right)} \left( {n,\beta } \right)} \right| \le \left( {2\left( {\tan \beta  - \beta } \right)} \right)^M \left| {U_M \left( {i\cot \beta } \right)} \right|\frac{{\Gamma \left( {n - M} \right)}}{{\Gamma \left( n \right)}}
\end{equation}
and
\begin{equation}\label{eq63}
\left| {A_M^{\left( 2 \right)} \left( {n,\beta } \right)} \right| \le \left( {2\left( {\tan \beta  - \beta  + \pi } \right)} \right)^M \left| {U_M \left( {i\cot \beta } \right)} \right|\frac{{\Gamma \left( {n - M} \right)}}{{\Gamma \left( n \right)}} .
\end{equation}
Expansions of type \eqref{eq60} are called inverse factorial series in the literature. Numerically, their character is similar to the character of asymptotic power series, because the consecutive Gamma functions decrease asymptotically by a factor $n$.

\begin{table*}[!ht]
\begin{center}
\begin{tabular}
[c]{ l r @{\,}c@{\,} l}\hline
 & \\ [-1ex]
 values of $\beta$ and $M$ &  $\beta=\frac{\pi}{6}$, $M=25$ & &\\ [1ex]
 exact numerical value of $U_{50}\left(i\cot\beta\right)$ & $-0.25922998993906050847874$ & $\times$ & $10^{111}$ \\ [1ex]
 Dingle's approximation \eqref{eq62} to $U_{50}\left(i\cot\beta\right)$ & $-0.25922998993906052149604$ & $\times$ & $10^{111}$  \\ [1ex]
 error & $0.1301729$ & $\times$ & $10^{95}$ \\ [1ex]
 approximation \eqref{eq60} to $U_{50}\left(i\cot\beta\right)$ & $-0.25922998993906052149604$ & $\times$ & $10^{111}$  \\ [1ex]
 error & $0.1301729$ & $\times$ & $10^{95}$\\ [1ex] 
  error bound using \eqref{eq59} and \eqref{eq63} & $0.2603745$ & $\times$ & $10^{95}$\\ [1ex] \hline
 & \\ [-1ex]
 values of $\beta$ and $M$ & $\beta=\frac{\pi}{3}$, $M=25$ & & \\ [1ex]
 exact numerical value of $U_{50}\left(i\cot\beta\right)$ & $-0.14230192249287421747599$ & $\times$ & $10^{56}$ \\ [1ex]
 Dingle's approximation \eqref{eq62} to $U_{50}\left(i\cot\beta\right)$ & $-0.14230192249287422461949$ & $\times$ & $10^{56}$  \\ [1ex]
 error & $0.714350$ & $\times$ & $10^{39}$ \\ [1ex]
 approximation \eqref{eq60} to $U_{50}\left(i\cot\beta\right)$ & $-0.14230192249287422461949$ & $\times$ & $10^{56}$  \\ [1ex]
 error & $0.714350$ & $\times$ & $10^{39}$\\ [1ex]
 error bound using \eqref{eq59} and \eqref{eq63} & $0.1428874$ & $\times$ & $10^{40}$\\ [1ex] \hline
 & \\ [-1ex]
 values of $\beta$ and $M$ & $\beta=\frac{6\pi}{13}$, $M=25$ & & \\ [1ex]
 exact numerical value of $U_{50}\left(i\cot\beta\right)$ & $-0.22522390129012627337081$ & $\times$ & $10^6$ \\ [1ex]
 Dingle's approximation \eqref{eq62} to $U_{50}\left(i\cot\beta\right)$ & $-0.22522390005970895996288$ & $\times$ & $10^6$  \\ [1ex]
 error & $-0.123041731340794$ & $\times$ & $10^{-2}$ \\ [1ex]
 approximation \eqref{eq60} to $U_{50}\left(i\cot\beta\right)$ & $-0.22522390129012628466504$ & $\times$ & $10^6$ \\ [1ex]
 error & $0.1129422$ & $\times$ & $10^{-10}$\\ [1ex]
 error bound using \eqref{eq59} and \eqref{eq63} & $0.2259222$ & $\times$ & $10^{-10}$\\ [1ex] \hline
 & \\ [-1ex]
 values of $\beta$ and $M$ & $\beta=\frac{7\pi}{15}$, $M=25$ & & \\ [1ex]
 exact numerical value of $U_{50}\left(i\cot\beta\right)$ & $-0.44399210443101419183462$ & $\times$ & $10^2$ \\ [1ex]
 Dingle's approximation \eqref{eq62} to $U_{50}\left(i\cot\beta\right)$ & $-0.44399207348793668328156$ & $\times$ & $10^2$ \\ [1ex]
 error & $-0.3094307750855306$ & $\times$ & $10^{-5}$ \\ [1ex]
 approximation \eqref{eq60} to $U_{50}\left(i\cot\beta\right)$ & $-0.44399210443101421410584$ & $\times$ & $10^2$ \\ [1ex]
 error & $0.2227122$ & $\times$ & $10^{-14}$\\ [1ex] 
 error bound using \eqref{eq59} and \eqref{eq63} & $0.4455002$ & $\times$ & $10^{-14}$\\ [-1ex]
 & \\\hline
\end{tabular}
\end{center}
\caption{Approximations for $U_{50}\left(i\cot\beta\right)$ with various $\beta$, using \eqref{eq60} and \eqref{eq62}.}
\label{table1}
\end{table*}

For large $n$, the least value of the bound \eqref{eq59} occurs when $M \approx \frac{n}{2}$. With this choice of $M$, the error bound is $\mathcal{O}\left( {n^{-\frac{1}{2}} 2^{-n} } \right)$. This is the best accuracy we can achieve using the expansion \eqref{eq60}.

We observe that, for large $n$, the contribution from the second sum in \eqref{eq60} is exponentially small compared to the first one, hence, by neglecting the second component and extending the first sum to infinity, we obtain the formal (non-convergent) expression
\begin{equation}\label{eq62}
U_n \left( {i\cot \beta } \right) \approx \frac{{\left(-1\right)^n \Gamma \left( n \right)}}{{2\pi \left( {2i\left( {\tan \beta  - \beta } \right)} \right)^n }}\sum\limits_{m = 0}^\infty  {\left( {2i\left( {\tan \beta  - \beta } \right)} \right)^m U_m \left( {i\cot \beta } \right)\frac{{\Gamma \left( {n - m} \right)}}{{\Gamma \left( n \right)}}} .
\end{equation}
This is exactly Dingle's expansion for the late coefficients in the asymptotic series of $J_\nu \left(\nu \sec \beta\right)$ and $Y_\nu \left(\nu \sec \beta\right)$ \cite[p. 170]{Dingle}. The mathematically rigorous form of Dingle's series is therefore the formula \eqref{eq60}.

Numerical examples illustrating the efficacy of the expansions \eqref{eq60} and \eqref{eq62}, truncated optimally, are given in Table \ref{table1}. It is seen from the computations that near $\beta = \frac{\pi}{2}$ the contribution from the second series in \eqref{eq60} becomes essential. Indeed, it can be shown that for large $n$ and for $\beta$ satisfying $\tan \beta  - \beta  \ge \pi$, the order of the main term in the second series is comparable with the last kept term in the first series assuming optimal truncation.

Let us now turn our attention to the coefficients $d_{2n}$. Expressing $H_{it}^{\left(1\right)}\left(it\right)$ by the formula \eqref{eq44} in the second integral in \eqref{eq15}, we deduce
\begin{equation}\label{eq61}
d_{2n}  = \frac{{\left( { - 1} \right)^n }}{{\sqrt 3 \pi \left( {2\pi } \right)^{\frac{{2n}}{3}} }}\left( {\sum\limits_{m = 0}^{M - 1} {\left( { - 1} \right)^m \frac{{2\sqrt 3 }}{3}\left( {2\pi } \right)^{\frac{{2m}}{3}} d_{2m} \sin \left( {\frac{{\left( {2m + 1} \right)\pi }}{3}} \right)\Gamma \left( {\frac{{2m + 1}}{3}} \right)\frac{{\Gamma \left( {\frac{{2\left( {n - m} \right)}}{3}} \right)}}{{\Gamma \left( {\frac{{2n + 1}}{3}} \right)}} + A_M \left( n \right)} } \right),
\end{equation}
for any fixed $0 \le M \le n - 2$, provided that $n\geq 2$. The remainder term $A_M \left( n \right)$ is given by the formula
\[
A_M \left( n \right) = \frac{\sqrt 3 \pi \left( {2\pi } \right)^{\frac{{2n}}{3}} i}{\Gamma \left( {\frac{{2n + 1}}{3}} \right)}\int_0^{ + \infty } {t^{\frac{{2n - 2}}{3}} e^{ - 2\pi t} R_M^{\left( H \right)} \left( {it} \right)dt} .
\]
Bounds for this error term follow from the estimates \eqref{eq51}--\eqref{eq38} as
\begin{equation}\label{eq66}
\left| {A_M \left( n \right)} \right| \le \left( {2\pi } \right)^{\frac{{2M}}{3}} \left| {d_{2M} } \right|\Gamma \left( {\frac{{2M + 1}}{3}} \right)\frac{\Gamma \left( {\frac{{2\left( {n - M} \right)}}{3}} \right)}{\Gamma \left( {\frac{{2n + 1}}{3}} \right)} + \left( {2\pi } \right)^{\frac{{2M + 2}}{3}} \left| {d_{2M + 2} } \right|\Gamma \left( {\frac{{2M + 3}}{3}} \right)\frac{\Gamma \left( {\frac{{2\left( {n - M - 1} \right)}}{3}} \right)}{\Gamma \left( {\frac{{2n + 1}}{3}} \right)},
\end{equation}
\begin{equation}\label{eq64}
\left| {A_M \left( n \right)} \right| \le \left( {2\pi } \right)^{\frac{{2M + 2}}{3}} \left| {d_{2M + 2} } \right|\Gamma \left( {\frac{{2M + 3}}{3}} \right)\frac{\Gamma \left( {\frac{{2\left( {n - M - 1} \right)}}{3}} \right)}{\Gamma \left( {\frac{{2n + 1}}{3}} \right)},
\end{equation}
\begin{equation}\label{eq65}
\left| {A_M \left( n \right)} \right| \le \left( {2\pi } \right)^{\frac{{2M}}{3}} \left| {d_{2M} } \right|\Gamma \left( {\frac{{2M + 1}}{3}} \right)\frac{\Gamma \left( {\frac{{2\left( {n - M} \right)}}{3}} \right)}{\Gamma \left( {\frac{{2n + 1}}{3}} \right)}
\end{equation}
according to whether $M\equiv 0 \mod 3$, $M\equiv 1 \mod 3$ or $M\equiv 2 \mod 3$, respectively.

For large $n$, the least values of these bounds occur when $M \approx \frac{n}{2}$. With this choice of $M$, the error bounds are $\mathcal{O}\left( {n^{ - \frac{1}{3}} 2^{ - \frac{{2n}}{3}} } \right)$. This is the best accuracy we can achieve using the expansion \eqref{eq61}. Numerical examples for various $n$ are provided in Table \ref{table2}.

\begin{table*}[!ht]
\begin{center}
\begin{tabular}
[c]{ l r @{\,}c@{\,} l}\hline
 & \\ [-1ex]
 values of $n$ and $M$ & $n=5$, $M=2$ & & \\ [1ex]
 exact numerical value of $d_{2n}$ & $-0.13204080504096204947934$ & $\times$ & $10^{-2}$ \\ [1ex]
 approximation \eqref{eq61} to $d_{2n}$ & $-0.13533519099105530338724$ & $\times$ & $10^{-2}$  \\ [1ex]
 error & $0.329438595009325390791$ & $\times$ & $10^{-4}$\\ [1ex]
 error bound using \eqref{eq65} & $0.741093493984769184813$ & $\times$ & $10^{-4}$\\ [1ex] \hline
 & \\ [-1ex]
 values of $n$ and $M$ & $n=10$, $M=5$ & & \\ [1ex]
 exact numerical value of $d_{2n}$ & $0.22835077129770834034682$ & $\times$ & $10^{-5}$ \\ [1ex]
 approximation \eqref{eq61} to $d_{2n}$ & $0.22863435953438111985117$ & $\times$ & $10^{-5}$  \\ [1ex]
 error & $-0.28358823667277950434$ & $\times$ & $10^{-8}$\\ [1ex]
 error bound using \eqref{eq65} & $0.82073950282859763109$ & $\times$ & $10^{-8}$\\ [1ex] \hline
 & \\ [-1ex]
 values of $n$ and $M$ & $n=25$, $M=12$ & &\\ [1ex]
 exact numerical value of $d_{2n}$ & $-0.17492836534785902720965$ & $\times$ & $10^{-13}$ \\ [1ex]
 approximation \eqref{eq61} to $d_{2n}$ & $-0.17492810163968211380528$ & $\times$ & $10^{-13}$ \\ [1ex]
 error & $-0.26370817691340437$ & $\times$ & $10^{-19}$\\ [1ex]
 error bound using \eqref{eq66} & $0.77781860849182240$ & $\times$ & $10^{-19}$\\ [1ex] \hline
 & \\ [-1ex]
 values of $n$ and $M$ & $n=50$, $M=25$ & &\\ [1ex]
 exact numerical value of $d_{2n}$ & $0.69006711932598958421355$ & $\times$ & $10^{-27}$ \\ [1ex]
 approximation \eqref{eq61} to $d_{2n}$ & $0.69006711932934022685189$ & $\times$ & $10^{-27}$ \\ [1ex]
 error & $-0.335064263835$ & $\times$ & $10^{-38}$\\ [1ex]
 error bound using \eqref{eq64} & $0.1051878923948$ & $\times$ & $10^{-37}$\\ [-1ex]
 & \\\hline
\end{tabular}
\end{center}
\caption{Approximations for $d_{2n}$ with various $n$, using \eqref{eq61}.}
\label{table2}
\end{table*}

By extending the sum in \eqref{eq61} to infinity, we arrive at the formal expansion of Dingle \cite[p. 171]{Dingle} for the late coefficients of the asymptotic series of $J_\nu \left(\nu\right)$ and $Y_\nu \left(\nu\right)$
\[
\frac{{\Gamma \left( {\frac{{2n + 1}}{3}} \right)}}{{\Gamma \left( {\frac{1}{3}} \right)6^{\frac{1}{3}} }}d_{2n}  \approx \frac{{\left( { - 1} \right)^n }}{{\sqrt 3 \pi \left( {2\pi } \right)^{\frac{{2n}}{3}} }}\left( {\Gamma \left( {\frac{{2n}}{3}} \right) - \frac{{3\beta_0 \left( {2\pi } \right)^{\frac{4}{3}} }}{{35}}\Gamma \left( {\frac{{2n}}{3} - \frac{4}{3}} \right) + \frac{{4\pi ^2 }}{{225}}\Gamma \left( {\frac{{2n}}{3} - 2} \right) +  \cdots } \right) .
\]
The constant $\beta_0$ is frequently used by Dingle in his discussion of the asymptotic expansions of integrals with second order saddle points. Its value is given by
\[
\beta_0  = \frac{{6^{\frac{1}{3}} \Gamma \left( {\frac{5}{3}} \right)}}{{4\Gamma \left( {\frac{1}{3}} \right)}} = \frac{{\Gamma \left( {\frac{5}{6}} \right)}}{{2 \cdot 3^{\frac{2}{3}} \sqrt \pi}} \approx 0.1530827453.
\]

We mention that Meissel \cite{Meissel} gave an asymptotic approximation, as $n \to \infty$, for the sequence
\[
\lambda_n  = \left( { - 1} \right)^n 6^{ - \frac{{2n + 1}}{3}} \frac{d_{2n}}{2n + 1} .
\]
A proof of his result, using Darboux's method, is given by Watson \cite[p. 233]{Watson} (see also Olver \cite[pp. 313--315]{Olver}). We show here that Meissel's approximation follows from \eqref{eq61}. Indeed, taking just the first term in \eqref{eq61} and making use of the approximation $\Gamma \left( {x + a} \right) \sim \Gamma \left( x \right)x^a$ yields
\begin{align*}
\lambda _n  & \sim \frac{{\Gamma \left( {\frac{1}{3}} \right)}}{{\sqrt 3 \pi }}\frac{1}{{\left( {12\pi } \right)^{\frac{{2n}}{3}} }}\frac{{\Gamma \left( {\frac{{2n}}{3}} \right)}}{{\left( {2n + 1} \right)\Gamma \left( {\frac{{2n + 1}}{3}} \right)}} = \frac{2}{{3\Gamma \left( {\frac{2}{3}} \right)\left( {12\pi } \right)^{\frac{{2n}}{3}} }}\frac{{\Gamma \left( {\frac{{2n}}{3}} \right)}}{{\left( {2n + 1} \right)\Gamma \left( {\frac{{2n + 1}}{3}} \right)}}
\\ & \sim \frac{2}{{3\Gamma \left( {\frac{2}{3}} \right)\left( {12\pi } \right)^{\frac{{2n}}{3}} }}\frac{1}{{\left( {2n + 1} \right)\left( {\frac{{2n}}{3}} \right)^{\frac{1}{3}} }} \sim \frac{1}{{18^{\frac{1}{3}} \Gamma \left( {\frac{2}{3}} \right)\left( {n + \frac{1}{3}} \right)^{\frac{4}{3}} \left( {12\pi } \right)^{\frac{{2n}}{3}} }}
\end{align*}
for large $n$. This is just Meissel's asymptotic formula.

More accurate approximations could be derived for the coefficients $U_n\left(i\cot \beta\right)$ and $d_{2n}$ by estimating the remainders $A_M^{\left( 1 \right)} \left( {n,\beta } \right)$, $A_M^{\left( 2 \right)} \left( {n,\beta } \right)$, $A_M\left(n\right)$, rather than bounding them. These estimates could perhaps provide an explanation of the apparent fact that for $A_M^{\left( 1 \right)} \left( {n,\beta } \right)$ and $A_M^{\left( 2 \right)} \left( {n,\beta } \right)$, the numerical error bound is always about twice the actual error as it is suggested by Table \ref{table1}; but we do not discuss the details here.

\section{Exponentially improved asymptotic expansions}\label{section5}

We shall find it convenient to express our exponentially improved expansions in terms of the (scaled) Terminant function, which is defined by
\[
\widehat T_p \left( z \right) = \frac{{e^{\pi ip} z^{1 - p} e^{ - z} }}{{2\pi i}}\int_0^{ + \infty } {\frac{{t^{p - 1} e^{ - t} }}{{z + t}}dt} \; \text{ for } \; p>0 \; \text{ and } \; \left| {\arg z} \right| < \pi ,
\]
and by analytic continuation elsewhere. We will also need an integral formula for the range $- 3\pi  < \arg z <  - \pi$. To find the required representation, we apply the connection formula
\[
\widehat T_p \left( z \right) = e^{2\pi ip} \left( {\widehat T_p \left( {ze^{2\pi i} } \right) - 1} \right),
\]
to obtain
\begin{equation}\label{eq72}
\widehat T_p \left( z \right) = \frac{{e^{\pi ip} z^{1 - p} e^{ - z} }}{{2\pi i}}\int_0^{ + \infty } {\frac{{t^{p - 1} e^{ - t} }}{{z + t}}dt}  - e^{2\pi ip} \; \text{ for } \; p>0 \; \text{ and } \; - 3\pi  < \arg z <  - \pi.
\end{equation}
Olver \cite{Olver4} showed that when $p \sim \left|z\right|$ and $z \to \infty$, we have
\begin{equation}\label{eq71}
ie^{ - \pi ip} \widehat T_p \left( z \right) = \begin{cases} \mathcal{O}\left( {e^{ - z - \left| z \right|} } \right) & \; \text{ if } \; \left| {\arg z} \right| \le \pi \\ \mathcal{O}\left(1\right) & \; \text{ if } \; - 3\pi  < \arg z \le  - \pi. \end{cases}
\end{equation}
Concerning the smooth transition of the Stokes discontinuities, we will use the more precise asymptotic formula
\begin{equation}\label{eq82}
e^{ - 2\pi ip} \widehat T_p \left( z \right) =  - \frac{1}{2} + \frac{1}{2}\mathop{\text{erf}} \left( { - \overline {c\left( { - \varphi } \right)} \sqrt {\frac{1}{2}\left| z \right|} } \right) + \mathcal{O}\left( {\frac{{e^{ - \frac{1}{2}\left| z \right|\overline {c^2 \left( { - \varphi } \right)} } }}{{\left| z \right|^{\frac{1}{2}} }}} \right)
\end{equation}
for $- 3\pi  + \delta  \le \arg z \le \pi  - \delta$, $0 < \delta \le 2\pi$. Here $\varphi = \arg z$ and $\text{erf}$ denotes the Error function. The quantity $c\left( \varphi  \right)$ is defined implicitly by the equation
\[
\frac{1}{2}c^2 \left( \varphi  \right) = 1 + i\left( {\varphi  - \pi } \right) - e^{i\left( {\varphi  - \pi } \right)},
\]
and corresponds to the branch of $c\left( \varphi  \right)$ which has the following expansion in the neighbourhood of $\varphi = \pi$:
\begin{equation}\label{eq83}
c\left( \varphi  \right) = \left( {\varphi  - \pi } \right) + \frac{i}{6}\left( {\varphi  - \pi } \right)^2  - \frac{1}{{36}}\left( {\varphi  - \pi } \right)^3  - \frac{i}{{270}}\left( {\varphi  - \pi } \right)^4  +  \cdots .
\end{equation}
For complete asymptotic expansions, see Olver \cite{Olver5}. We remark that Olver uses the different notation $F_p \left( z \right) = ie^{ - \pi ip} \widehat T_p \left( z \right)$ for the Terminant function and the other branch of the function $c\left( \varphi  \right)$. For further properties of the Terminant function, see, for example, Paris and Kaminski \cite[Chapter 6]{Paris2}.

\subsection{Proof of the exponentially improved expansions for $H_\nu^{\left(1\right)}\left(\nu x\right)$}

\subsubsection{Case (i): $x>1$} First, we suppose that $-\frac{\pi}{2} < \arg \nu < \frac{3\pi}{2}$. Our starting point is the representation \eqref{eq21} with $N=0$
\begin{align*}
H_\nu ^{\left( 1 \right)} \left( {\nu \sec \beta } \right) = \; & \frac{{e^{i\nu \left( {\tan \beta  - \beta } \right) - \frac{\pi }{4}i} }}{{2\pi \nu ^{\frac{1}{2}} }}\int_0^{ + \infty } {\frac{{t^{ - \frac{1}{2}} e^{ - t\left( {\tan \beta  - \beta } \right)} }}{{1 + it/\nu }}iH_{it}^{\left( 1 \right)} \left( {it\sec \beta } \right)dt} \\ & + \frac{{e^{i\nu \left( {\tan \beta  - \beta } \right) - \frac{\pi }{4}i} }}{{2\pi \nu ^{\frac{1}{2}} }}\int_0^{ + \infty } {\frac{{t^{ - \frac{1}{2}} e^{ - t\left( {\tan \beta  - \beta  + 2\pi } \right)} }}{{1 + it/\nu }}iH_{it}^{\left( 1 \right)} \left( {it\sec \beta } \right)dt} .
\end{align*}
Let $M$ and $N$ be integers such that $1 \le N \le M$. An application of the expansion \eqref{eq7} shows that
\[
H_\nu ^{\left( 1 \right)} \left( {\nu \sec \beta } \right) = \frac{{e^{i\nu \left( {\tan \beta  - \beta } \right) - \frac{\pi }{4}i} }}{{\left( {\frac{1}{2}\nu \pi \tan \beta } \right)^{\frac{1}{2}} }}\left(\sum\limits_{n = 0}^{N - 1} {\left( { - 1} \right)^n \frac{{U_n \left( {i\cot \beta } \right)}}{{\nu ^n }}}  + \sum\limits_{m = N}^{M - 1} {\left( { - 1} \right)^m \frac{{\widetilde U_m \left( {i\cot \beta } \right)}}{{\nu ^m }}}  + R_{N,M}^{\left( H \right)} \left( {\nu ,\beta } \right)\right),
\]
where
\begin{gather}\label{eq68}
\begin{split}
R_{N,M}^{\left( H \right)} \left( {\nu ,\beta } \right) = \; & \frac{1}{{2\left( {2\pi \cot \beta } \right)^{\frac{1}{2}} \left( {i\nu } \right)^N }}\int_0^{ + \infty } {\frac{{t^{N - \frac{1}{2}} e^{ - t\left( {\tan \beta  - \beta } \right)} }}{{1 + it/\nu }}iH_{it}^{\left( 1 \right)} \left( {it\sec \beta } \right)dt} \\
& + \frac{1}{{2\left( {2\pi \cot \beta } \right)^{\frac{1}{2}} \left( {i\nu } \right)^M }}\int_0^{ + \infty } {\frac{{t^{M - \frac{1}{2}} e^{ - t\left( {\tan \beta  - \beta  + 2\pi } \right)} }}{{1 + it/\nu }}iH_{it}^{\left( 1 \right)} \left( {it\sec \beta } \right)dt} 
\end{split}
\end{gather}
and
\[
\widetilde U_m \left( {i\cot \beta } \right) = \frac{{i^m }}{{2\left( {2\pi \cot \beta } \right)^{\frac{1}{2}} }}\int_0^{ + \infty } {t^{m - \frac{1}{2}} e^{ - t\left( {\tan \beta  - \beta  + 2\pi } \right)} iH_{it}^{\left( 1 \right)} \left( {it\sec \beta } \right)dt} .
\]
We remark that $R_{N,N}^{\left( H \right)} \left( {\nu ,\beta } \right) = R_N^{\left( H \right)} \left( {\nu ,\beta } \right)$. Suppose that $K$ and $L$ are integers such that $0 \le K < N$, $0 \le L < M$. We use \eqref{eq21} again to expand the function $H_{it}^{\left( 1 \right)} \left( {it\sec \beta } \right)$ under the integrals in \eqref{eq68}, to obtain
\begin{gather}\label{eq67}
\begin{split}
R_{N,M}^{\left( H \right)} \left( {\nu ,\beta } \right) = \; & \sum\limits_{k = 0}^{K - 1} {\frac{{U_k \left( {i\cot \beta } \right)}}{{\nu ^k }}\frac{{\left( {i\nu } \right)^{k - N} }}{{2\pi }}\int_0^{ + \infty } {\frac{{t^{N - k - 1} e^{ - 2t\left( {\tan \beta  - \beta } \right)} }}{{1 + it/\nu }}dt} }  \\ & + \sum\limits_{\ell  = 0}^{L - 1} {\frac{{U_\ell  \left( {i\cot \beta } \right)}}{{\nu ^\ell  }}\frac{{\left( {i\nu } \right)^{\ell  - M} }}{{2\pi }}\int_0^{ + \infty } {\frac{{t^{M - \ell  - 1} e^{ - 2t\left( {\tan \beta  - \beta  + \pi } \right)} }}{{1 + it/\nu }}dt} } \\ & + R_{N,M,K,L}^{\left( H \right)} \left( {\nu ,\beta } \right),
\end{split}
\end{gather}
with
\begin{gather}\label{eq69}
\begin{split}
R_{N,M,K,L}^{\left( H \right)} \left( {\nu ,\beta } \right) = \; & \frac{1}{{2\pi \left( {i\nu } \right)^N }}\int_0^{ + \infty } {\frac{{t^{N - 1} e^{ - 2t\left( {\tan \beta  - \beta } \right)} }}{{1 + it/\nu }}R_K^{\left( H \right)} \left( {it,\beta } \right)dt} \\ & + \frac{1}{{2\pi \left( {i\nu } \right)^M }}\int_0^{ + \infty } {\frac{{t^{M - 1} e^{ - 2t\left( {\tan \beta  - \beta  + \pi } \right)} }}{{1 + it/\nu }}R_L^{\left( H \right)} \left( {it,\beta } \right)dt} .
\end{split}
\end{gather}
The integrals in \eqref{eq67} can be identified in terms of the Terminant function since
\[
\frac{{\left( {i\nu } \right)^{k - N} }}{{2\pi }}\int_0^{ + \infty } {\frac{{t^{N - k - 1} e^{ - 2t\left( {\tan \beta  - \beta } \right)} }}{{1 + it/\nu }}dt}  = ie^{ - 2i\nu\left( {\tan \beta  - \beta } \right) } \widehat T_{N - k} \left( { - 2i\nu\left( {\tan \beta  - \beta } \right) } \right)
\]
and
\[
\frac{{\left( {i\nu } \right)^{\ell  - M} }}{{2\pi }}\int_0^{ + \infty } {\frac{{t^{M - \ell  - 1} e^{ - 2t\left( {\tan \beta  - \beta  + \pi } \right)} }}{{1 + it/\nu }}dt}  = ie^{ - 2i\nu\left( {\tan \beta  - \beta  + \pi } \right) } \widehat T_{M - \ell } \left( { - 2i\nu\left( {\tan \beta  - \beta  + \pi } \right) } \right) .
\]
Therefore, we have the following expansion
\begin{align*}
R_{N,M}^{\left( H \right)} \left( {\nu ,\beta } \right) = \; & ie^{ - 2i\nu\left( {\tan \beta  - \beta } \right) } \sum\limits_{k = 0}^{K - 1} {\frac{{U_k \left( {i\cot \beta } \right)}}{{\nu ^k }}\widehat T_{N - k} \left( { - 2i\nu\left( {\tan \beta  - \beta } \right) } \right)} \\ & + ie^{ - 2i \nu \left( {\tan \beta  - \beta  + \pi } \right) } \sum\limits_{\ell  = 0}^{L - 1} {\frac{{U_\ell  \left( {i\cot \beta } \right)}}{{\nu ^\ell  }}\widehat T_{M - \ell } \left( { - 2i\nu\left( {\tan \beta  - \beta  + \pi } \right)} \right)}  + R_{N,M,K,L}^{\left( H \right)} \left( {\nu ,\beta } \right) .
\end{align*}
Taking $\nu  = re^{i\theta }$, the representation \eqref{eq69} becomes
\begin{gather}\label{eq70}
\begin{split}
R_{N,M,K,L}^{\left( H \right)} \left( {\nu ,\beta } \right) = \; & \frac{1}{{2\pi \left( {ie^{i\theta } } \right)^N }}\int_0^{ + \infty } {\frac{{\tau ^{N - 1} e^{ - 2r\tau \left( {\tan \beta  - \beta } \right)} }}{{1 + i\tau e^{ - i\theta } }}R_K^{\left( H \right)} \left( {ir\tau ,\beta } \right)d\tau } \\ & + \frac{1}{{2\pi \left( {ie^{i\theta } } \right)^M }}\int_0^{ + \infty } {\frac{{\tau ^{M - 1} e^{ - 2r\tau \left( {\tan \beta  - \beta +\pi } \right)} }}{{1 + i\tau e^{ - i\theta } }}R_L^{\left( H \right)} \left( {ir\tau ,\beta } \right)d\tau } .
\end{split}
\end{gather}
Using the integral formula \eqref{eq12}, $R_K^{\left( H \right)} \left( {ir\tau ,\beta } \right)$ can be written in the form
\begin{multline*}
R_K^{\left( H \right)} \left( {ir\tau ,\beta } \right) = \frac{{\left( { - 1} \right)^K }}{{2\left( {2\pi \cot \beta } \right)^{\frac{1}{2}}\left( {r\tau } \right)^K }}\left( {\int_0^{ + \infty } {\frac{{s^{K - \frac{1}{2}} e^{ - s\left( {\tan \beta  - \beta } \right)} }}{{1 + s/r}}\left( {1 + e^{ - 2\pi s} } \right)iH_{is}^{\left( 1 \right)} \left( {is\sec \beta } \right)ds}}\right. \\  + \left.{\left( {\tau  - 1} \right)\int_0^{ + \infty } {\frac{{s^{K - \frac{1}{2}} e^{ - s\left( {\tan \beta  - \beta } \right)} }}{{\left( 1+r\tau /s \right)\left( {1 + s/r} \right)}}\left( {1 + e^{ - 2\pi s} } \right)iH_{is}^{\left( 1 \right)} \left( {is\sec \beta } \right)ds} } \right),
\end{multline*}
and similarly for $R_L^{\left( H \right)} \left( {ir\tau ,\beta } \right)$. Noting that
\[
0 < \frac{1}{{1 + s/r}},\frac{1}{{\left( 1 + r\tau /s \right)\left( {1 + s/r} \right)}} < 1
\]
for positive $r$, $\tau$ and $s$, substitution into \eqref{eq70} yields the upper bound
\begin{align*}
\left| {R_{N,M,K,L}^{\left( H \right)} \left( {\nu ,\beta } \right)} \right| \le \; & \frac{{\left| {U_K \left( {i\cot \beta } \right)} \right|}}{{\left| \nu  \right|^K }}\left| {\frac{1}{{2\pi }}\int_0^{ + \infty } {\frac{{\tau ^{N - K - 1} e^{ - 2r\tau \left( {\tan \beta  - \beta } \right)} }}{{1 + i\tau e^{ - i\theta } }}d\tau } } \right| \\ & + \frac{{\left| {U_K \left( {i\cot \beta } \right)} \right|}}{{\left| \nu  \right|^K }}\frac{1}{{2\pi }}\int_0^{ + \infty } {\tau ^{N - K - 1} e^{ - 2r\tau \left( {\tan \beta  - \beta } \right)} \left| {\frac{{\tau  - 1}}{{\tau  - ie^{i\theta } }}} \right|d\tau } \\ & + \frac{{\left| {U_L \left( {i\cot \beta } \right)} \right|}}{{\left| \nu  \right|^L }}\left| {\frac{1}{{2\pi }}\int_0^{ + \infty } {\frac{{\tau ^{M - L - 1} e^{ - 2r\tau \left( {\tan \beta  - \beta +\pi} \right)} }}{{1 + i\tau e^{ - i\theta } }}d\tau } } \right| \\ &+ \frac{{\left| {U_L \left( {i\cot \beta } \right)} \right|}}{{\left| \nu  \right|^L }}\frac{1}{{2\pi }}\int_0^{ + \infty } {\tau ^{M - L - 1} e^{ - 2r\tau \left( {\tan \beta  - \beta +\pi } \right)} \left| {\frac{{\tau  - 1}}{{\tau  - ie^{i\theta } }}} \right|d\tau } .
\end{align*}
Since $\left| {\left( {\tau  - 1} \right)/\left( {\tau  - ie^{i\theta } } \right)} \right| \le 1$, we find that
\begin{align*}
\left| {R_{N,M,K,L}^{\left( H \right)} \left( {\nu ,\beta } \right)} \right| \le \; & \frac{{\left| {U_K \left( {i\cot \beta } \right)} \right|}}{{\left| \nu  \right|^K }}\left| {ie^{ - 2i\nu\left( {\tan \beta  - \beta } \right) } \widehat T_{N - K} \left( { - 2i\nu\left( {\tan \beta  - \beta } \right)} \right)} \right| \\ &+ \frac{{\left| {U_K \left( {i\cot \beta } \right)} \right|\Gamma \left( {N - K} \right)}}{{2\pi \left( {2\left( {\tan \beta  - \beta } \right)} \right)^{N - K} \left| \nu  \right|^N }}
\\& + \frac{{\left| {U_L \left( {i\cot \beta } \right)} \right|}}{{\left| \nu  \right|^L }}\left| {ie^{ - 2i\nu\left( {\tan \beta  - \beta +\pi } \right) } \widehat T_{M - L} \left( { - 2i\nu\left( {\tan \beta  - \beta +\pi} \right) } \right)} \right| \\ & + \frac{{\left| {U_L \left( {i\cot \beta } \right)} \right|\Gamma \left( {M - L} \right)}}{{2\pi \left( {2\left( {\tan \beta  - \beta +\pi } \right)} \right)^{M - L} \left| \nu  \right|^M }}.
\end{align*}
By continuity, this bound holds in the closed sector $-\frac{\pi}{2} \le \arg \nu  \le \frac{3\pi}{2}$. Suppose that $
N = 2\left| \nu  \right|\left( {\tan \beta  - \beta } \right) + \rho$, $M = 2\left| \nu  \right|\left( {\tan \beta  - \beta  + \pi } \right) + \sigma$ where $\rho$ and $\sigma$ are bounded. Employing Stirling's formula, we find that
\[
\frac{{\left| {U_K \left( {i\cot \beta } \right)} \right|\Gamma \left( {N - K} \right)}}{{2\pi \left( {2\left( {\tan \beta  - \beta } \right)} \right)^{N - K} \left| \nu  \right|^N }} = \mathcal{O}_{K,\rho} \left( {\frac{{e^{ - 2\left| \nu  \right|\left( {\tan \beta  - \beta } \right)} }}{{\left( {2\left| \nu  \right|\left( {\tan \beta  - \beta } \right)} \right)^{\frac{1}{2}} }}\frac{{\left| {U_K \left( {i\cot \beta } \right)} \right|}}{{\left| \nu  \right|^K }}} \right)
\]
and
\[
\frac{{\left| {U_L \left( {i\cot \beta } \right)} \right|\Gamma \left( {M - L} \right)}}{{2\pi \left( {2\left( {\tan \beta  - \beta  + \pi } \right)} \right)^{M - L} \left| \nu  \right|^M }} = \mathcal{O}_{L,\sigma} \left( {\frac{{e^{ - 2\left| \nu  \right|\left( {\tan \beta  - \beta  + \pi } \right)} }}{{\left( {2\left| \nu  \right|\left( {\tan \beta  - \beta  + \pi } \right)} \right)^{\frac{1}{2}} }}\frac{{\left| {U_L \left( {i\cot \beta } \right)} \right|}}{{\left| \nu  \right|^L }}} \right)
\]
as $\nu \to \infty$. Olver's estimation \eqref{eq71} shows that
\[
\left| {ie^{ - 2i\nu\left( {\tan \beta  - \beta } \right) } \widehat T_{N - K} \left( { - 2i\nu\left( {\tan \beta  - \beta } \right)} \right)} \right| = \mathcal{O}_{K,\rho} \left( {e^{ - 2\left|\nu\right|\left( {\tan \beta  - \beta } \right)}} \right)
\]
and
\[
\left| {ie^{ - 2i\nu\left( {\tan \beta  - \beta  + \pi } \right) }\widehat T_{M - L} \left( { - 2i\nu\left( {\tan \beta  - \beta  + \pi } \right)} \right)} \right| = \mathcal{O}_{L,\sigma} \left( {e^{ - 2\left| \nu  \right|\left( {\tan \beta  - \beta  + \pi } \right)} } \right)
\]
for large $\nu$. Therefore, we have that
\[
R_{N,M,K,L}^{\left( H \right)} \left( {\nu ,\beta } \right) = \mathcal{O}_{K,\rho } \left( {e^{ - 2\left| \nu  \right|\left( {\tan \beta  - \beta } \right)} \frac{{\left|U_K \left( {i\cot \beta } \right)\right|}}{{\left| \nu  \right|^K }}} \right) + \mathcal{O}_{L,\sigma } \left( {e^{ - 2\left| \nu  \right|\left( {\tan \beta  - \beta  + \pi } \right)} \frac{{\left|U_L \left( {i\cot \beta } \right)\right|}}{{\left| \nu  \right|^L }}} \right)
\]
as $\nu \to \infty$ in the sector $-\frac{\pi}{2} \le \arg \nu \le \frac{3\pi}{2}$.

Consider now the sector $-\frac{3\pi}{2} \leq \arg\nu \leq -\frac{\pi}{2}$. Rotating the path of integration in \eqref{eq69} and employing the residue theorem yields
\begin{gather}\label{eq57}
\begin{split}
R_{N,M,K,L}^{\left( H \right)} \left( {\nu ,\beta } \right) = & - ie^{ - 2i\nu \left( {\tan \beta  - \beta } \right)} R_K^{\left( H \right)} \left( { \nu e^{\pi i},\beta } \right) + \frac{1}{{2\pi \left( {i\nu } \right)^N }}\int_0^{ + \infty } {\frac{{t^{N - 1} e^{ - 2t\left( {\tan \beta  - \beta } \right)} }}{{1 + it/\nu }}R_K^{\left( H \right)} \left( {it,\beta } \right)dt} \\
& - ie^{ - 2i\nu \left( {\tan \beta  - \beta  + \pi } \right)} R_L^{\left( H \right)} \left( { \nu e^{\pi i},\beta } \right) + \frac{1}{{2\pi \left( {i\nu } \right)^M }}\int_0^{ + \infty } {\frac{{t^{M - 1} e^{ - 2t\left( {\tan \beta  - \beta  + \pi } \right)} }}{{1 + it/\nu }}R_L^{\left( H \right)} \left( {it,\beta } \right)dt} 
\end{split}
\end{gather}
when $-\frac{3\pi}{2} < \arg\nu  < -\frac{\pi}{2}$. Proceeding similarly as before, applying the integral formula \eqref{eq72}, we find
\begin{align*}
\left| {R_{N,M,K,L}^{\left( H \right)} \left( {\nu ,\beta } \right)} \right| \le \; & \left| {e^{ - 2i\nu \left( {\tan \beta  - \beta } \right)} R_K^{\left( H \right)} \left( { \nu e^{\pi i},\beta } \right)} \right| 
\\ & + \frac{{\left| {U_K \left( {i\cot \beta } \right)} \right|}}{{\left| \nu  \right|^K }}\left| {ie^{ - 2i\nu\left( {\tan \beta  - \beta } \right) } \left( {\widehat T_{N - K} \left( { - 2i\nu\left( {\tan \beta  - \beta } \right)} \right) + 1} \right)} \right|\\
& + \frac{{\left| {U_K \left( {i\cot \beta } \right)} \right|\Gamma \left( {N - K} \right)}}{{2\pi \left( {2\left( {\tan \beta  - \beta } \right)} \right)^{N - K} \left| \nu  \right|^N }}
 + \left| {e^{ - 2i\nu \left( {\tan \beta  - \beta  + \pi } \right)} R_L^{\left( H \right)} \left( { \nu e^{\pi i},\beta } \right)} \right|\\
& + \frac{{\left| {U_L \left( {i\cot \beta } \right)} \right|}}{{\left| \nu  \right|^L }}\left| {i e^{ - 2i\nu\left( {\tan \beta  - \beta  + \pi } \right) } \left( {\widehat T_{M - L} \left( { - 2i\nu\left( {\tan \beta  - \beta  + \pi } \right) } \right) + 1} \right)} \right| \\ &
 + \frac{{\left| {U_L \left( {i\cot \beta } \right)} \right|\Gamma \left( {M - L} \right)}}{{2\pi \left( {2\left( {\tan \beta  - \beta  + \pi } \right)} \right)^{M - L} \left| \nu  \right|^M }}.
\end{align*}
By continuity, this bound holds in the closed sector $ - \frac{3\pi}{2} \leq \arg \nu  \leq  - \frac{\pi}{2}$. As before, we take $N = 2\left| \nu  \right|\left( {\tan \beta  - \beta } \right) + \rho$, $M = 2\left| \nu  \right|\left( {\tan \beta  - \beta  + \pi } \right) + \sigma$ with $\rho$ and $\sigma$ being bounded. Olver's estimate \eqref{eq71} gives
\[
 \left| {i e^{ - 2i\nu\left( {\tan \beta  - \beta } \right) } \left( {\widehat T_{N - K} \left( { - 2i\nu\left( {\tan \beta  - \beta } \right)} \right) + 1} \right)} \right| = \mathcal{O}_{K,\rho } \left( {e^{2\Im \left( \nu  \right)\left( {\tan \beta  - \beta } \right)} } \right)
\]
and
\[
 \left| {i e^{ - 2i\nu\left( {\tan \beta  - \beta  + \pi } \right) } \left( {\widehat T_{M - L} \left( { - 2i\nu\left( {\tan \beta  - \beta  + \pi } \right) } \right) + 1} \right)} \right| = \mathcal{O}_{L,\sigma } \left( {e^{2\Im \left( \nu  \right)\left( {\tan \beta  - \beta  + \pi } \right)} } \right)
\]
for $\nu$ large and $-\frac{3\pi}{2} \leq \arg\nu \leq -\frac{\pi}{2}$. From \eqref{eq73} and \eqref{eq27}, we deduce that
\[
\left| {e^{ - 2i\nu \left( {\tan \beta  - \beta } \right)} R_K^{\left( H \right)} \left( { \nu e^{\pi i},\beta } \right)} \right| = \mathcal{O}_K  \left( {e^{2\Im \left( \nu  \right)\left( {\tan \beta  - \beta } \right)} \frac{{\left| {U_K \left( {i\cot \beta } \right)} \right|}}{{\left| \nu  \right|^K }}} \right)
\]
and
\[
\left| {e^{ - 2i\nu \left( {\tan \beta  - \beta  + \pi } \right)} R_L^{\left( H \right)} \left( { \nu e^{\pi i},\beta } \right)} \right| = \mathcal{O}_L  \left( {e^{2\Im \left( \nu  \right)\left( {\tan \beta  - \beta  + \pi } \right)} \frac{{\left| {U_L \left( {i\cot \beta } \right)} \right|}}{{\left| \nu  \right|^L }}} \right)
\]
for large $\nu$ with $- \frac{{3\pi }}{2} \leq \arg \nu  \leq  - \frac{\pi }{2}$. Therefore, we have
\[
R_{N,M,K,L}^{\left( H \right)} \left( {\nu ,\beta } \right) = \mathcal{O}_{K,\rho} \left( {e^{2\Im \left( \nu  \right)\left( {\tan \beta  - \beta } \right)} \frac{{\left| {U_K \left( {i\cot \beta } \right)} \right|}}{{\left| \nu  \right|^K }}} \right) + \mathcal{O}_{L,\sigma} \left( {e^{2\Im \left( \nu  \right)\left( {\tan \beta  - \beta  + \pi } \right)} \frac{{\left| {U_L \left( {i\cot \beta } \right)} \right|}}{{\left| \nu  \right|^L }}} \right)
\]
as $\nu \to \infty$ in the sector $- \frac{{3\pi }}{2} \leq \arg \nu  \leq  - \frac{\pi }{2}$.

\subsubsection{Case (ii): $x=1$} First, we suppose that $- \frac{\pi }{2} < \arg \nu < \frac{\pi }{2}$. We write \eqref{eq44} with $N=0$ in the form
\[
H_\nu ^{\left( 1 \right)} \left( \nu  \right) = \frac{{e^{ - \frac{\pi}{3}i} }}{{\sqrt 3 \pi \nu ^{\frac{1}{3}} }}\int_0^{ + \infty } {\frac{{t^{ - \frac{2}{3}} e^{ - 2\pi t} }}{{1 + \left( {t/\nu } \right)^2 }}iH_{it}^{\left( 1 \right)} \left( {it} \right)dt}  - \frac{{e^{\frac{\pi}{3}i} }}{{\sqrt 3 \pi \nu ^{\frac{5}{3}} }}\int_0^{ + \infty } {\frac{{t^{\frac{2}{3}} e^{ - 2\pi t} }}{{1 + \left( {t/\nu } \right)^2 }}iH_{it}^{\left( 1 \right)} \left( {it} \right)dt} .
\]
Let $M$ and $N$ be arbitrary positive integers. Using the expression \eqref{eq7}, we find that
\[
H_\nu ^{\left( 1 \right)} \left( \nu  \right) = \frac{{e^{ - \frac{\pi}{3}i} }}{{\sqrt 3 \pi \nu ^{\frac{1}{3}} }}\sum\limits_{n = 0}^{N - 1} {d_{6n} \frac{{\Gamma \left( {2n+\frac{1}{3}} \right)}}{{\nu ^{2n} }}}  - \frac{{e^{\frac{\pi}{3}i} }}{{\sqrt 3 \pi \nu ^{\frac{5}{3}} }}\sum\limits_{m = 0}^{M - 1} {d_{6m + 4} \frac{{\Gamma \left( {2m+\frac{5}{3}} \right)}}{{\nu ^{2m} }}}  + R_{N,M}^{\left( H \right)} \left( \nu  \right),
\]
where
\begin{gather}\label{eq75}
\begin{split}
R_{N,M}^{\left( H \right)} \left( \nu  \right) = \; & \left( { - 1} \right)^N \frac{{e^{ - \frac{\pi}{3}i} }}{{\sqrt 3 \pi \nu ^{2N + \frac{1}{3}} }}\int_0^{ + \infty } {\frac{{t^{2N - \frac{2}{3}} e^{ - 2\pi t} }}{{1 + \left( {t/\nu } \right)^2 }}iH_{it}^{\left( 1 \right)} \left( {it} \right)dt} \\ & + \left( { - 1} \right)^{M - 1} \frac{{e^{\frac{\pi}{3}i} }}{{\sqrt 3 \pi \nu ^{2M + \frac{5}{3}} }}\int_0^{ + \infty } {\frac{{t^{2M + \frac{2}{3}} e^{ - 2\pi t} }}{{1 + \left( {t/\nu } \right)^2 }}iH_{it}^{\left( 1 \right)} \left( {it} \right)dt} .
\end{split}
\end{gather}
We remark that $R_{N,N}^{\left( H \right)} \left( \nu  \right) = R_{3N}^{\left( H \right)} \left( \nu  \right)$. Assume that $K$ and $L$ are integers such that $0 \leq K < 3N$, $0 \leq L < 3M + 2$ and $K,L \equiv 0 \mod 3$. We apply \eqref{eq44} again to expand the function $H_{it}^{\left( 1 \right)} \left( {it} \right)$ under the integral in \eqref{eq75}, to obtain
\begin{gather}\label{eq76}
\begin{split}
R_{N,M}^{\left( H \right)}  = \; & \frac{{e^{ - \frac{\pi }{3}i} }}{{\sqrt 3 }}\frac{2}{{3\pi }}\sum\limits_{k = 0}^{K - 1} {d_{2k} \sin \left( {\frac{{\left( {2k + 1} \right)\pi }}{3}} \right)\frac{{\Gamma \left( {\frac{{2k + 1}}{3}} \right)}}{{\nu ^{\frac{{2k + 1}}{3}} }}\left( { - 1} \right)^{N + k} \frac{{\nu ^{\frac{{2k}}{3} - 2N} }}{\pi }\int_0^{ + \infty } {\frac{{t^{2N - \frac{{2k}}{3} - 1} e^{ - 2\pi t} }}{{1 + \left( {t/\nu } \right)^2 }}dt} } \\
& - \frac{{e^{\frac{\pi }{3}i} }}{{\sqrt 3 }}\frac{2}{{3\pi }}\sum\limits_{\ell  = 0}^{L - 1} {d_{2\ell } \sin \left( {\frac{{\left( {2\ell  + 1} \right)\pi }}{3}} \right)\frac{{\Gamma \left( {\frac{{2\ell  + 1}}{3}} \right)}}{{\nu ^{\frac{{2\ell  + 1}}{3}} }}\left( { - 1} \right)^{M + \ell } \frac{{\nu ^{ \frac{2\ell  - 4}{3}- 2M} }}{\pi }\int_0^{ + \infty } {\frac{{t^{2M - \frac{{2\ell  - 4}}{3} - 1} e^{ - 2\pi t} }}{{1 + \left( {t/\nu } \right)^2 }}dt} } \\
& + R_{N,M,K,L}^{\left( H \right)} \left( \nu  \right),
\end{split}
\end{gather}
with
\begin{gather}\label{eq77}
\begin{split}
R_{N,M,K,L}^{\left( H \right)} \left( \nu  \right) = & \left( { - 1} \right)^N \frac{{e^{ - \frac{\pi }{3}i} }}{{\sqrt 3 \pi \nu ^{2N + \frac{1}{3}} }}\int_0^{ + \infty } {\frac{{t^{2N - \frac{2}{3}} e^{ - 2\pi t} }}{{1 + \left( {t/\nu } \right)^2 }}iR_K^{\left( H \right)} \left( {it} \right)dt} \\ & + \left( { - 1} \right)^{M - 1} \frac{{e^{\frac{\pi }{3}i} }}{{\sqrt 3 \pi \nu ^{2M + \frac{5}{3}} }}\int_0^{ + \infty } {\frac{{t^{2M + \frac{2}{3}} e^{ - 2\pi t} }}{{1 + \left( {t/\nu } \right)^2 }}iR_L^{\left( H \right)} \left( {it} \right)dt} .
\end{split}
\end{gather}
The integrals in \eqref{eq76} can be identified in terms of the Terminant function since
\[
\left( { - 1} \right)^{N + k} \frac{{\nu ^{\frac{{2k}}{3}- 2N} }}{\pi }\int_0^{ + \infty } {\frac{{t^{2N - \frac{{2k}}{3} - 1} e^{ - 2\pi t} }}{{1 + \left( {t/\nu } \right)^2 }}dt}  = ie^{ - 2\pi i\nu } \widehat T_{2N - \frac{{2k}}{3}} \left( { - 2\pi i\nu } \right) - ie^{\frac{\pi }{3}i} e^{2\pi i\nu } e^{\frac{{2\left( {2k + 1} \right)\pi i}}{3}} \widehat T_{2N - \frac{{2k}}{3}} \left( {2\pi i\nu } \right)
\]
and
\begin{align*}
\left( { - 1} \right)^{M + \ell } \frac{{\nu ^{ \frac{{2\ell  - 4}}{3}- 2M} }}{\pi }\int_0^{ + \infty } {\frac{{t^{2M - \frac{{2\ell  - 4}}{3} - 1} e^{ - 2\pi t} }}{{1 + \left( {t/\nu } \right)^2 }}dt}  = \; & ie^{ - 2\pi i\nu } \widehat T_{2M - \frac{{2\ell  - 4}}{3}} \left( { - 2\pi i\nu } \right)\\ & - ie^{ - \frac{\pi }{3}i} e^{2\pi i\nu } e^{\frac{{2\left( {2\ell  + 1} \right)\pi i}}{3}} \widehat T_{2M - \frac{{2\ell  - 4}}{3}} \left( {2\pi i\nu } \right) .
\end{align*}
Hence, we have the following expansion
\begin{align*}
R_{N,M}^{\left( H \right)} \left( \nu  \right) = \; & ie^{ - \frac{\pi }{3}i} \frac{{e^{ - 2\pi i\nu } }}{{\sqrt 3 }}\frac{2}{{3\pi }}\sum\limits_{k = 0}^{K - 1} {d_{2k} \sin \left( {\frac{{\left( {2k + 1} \right)\pi }}{3}} \right)\frac{{\Gamma \left( {\frac{{2k + 1}}{3}} \right)}}{{\nu ^{\frac{{2k + 1}}{3}} }}\widehat T_{2N - \frac{{2k}}{3}} \left( { - 2\pi i\nu } \right)} \\
& - i\frac{{e^{2\pi i\nu } }}{{\sqrt 3 }}\frac{2}{{3\pi }}\sum\limits_{k = 0}^{K - 1} {d_{2k} e^{\frac{{2\left( {2k + 1} \right)\pi i}}{3}} \sin \left( {\frac{{\left( {2k + 1} \right)\pi }}{3}} \right)\frac{{\Gamma \left( {\frac{{2k + 1}}{3}} \right)}}{{\nu ^{\frac{{2k + 1}}{3}} }}\widehat T_{2N - \frac{{2k}}{3}} \left( {2\pi i\nu } \right)} \\
& - ie^{\frac{\pi }{3}i} \frac{{e^{ - 2\pi i\nu } }}{{\sqrt 3 }}\frac{2}{{3\pi }}\sum\limits_{\ell  = 0}^{L - 1} {d_{2\ell } \sin \left( {\frac{{\left( {2\ell  + 1} \right)\pi }}{3}} \right)\frac{{\Gamma \left( {\frac{{2\ell  + 1}}{3}} \right)}}{{\nu ^{\frac{{2\ell  + 1}}{3}} }}\widehat T_{2M - \frac{{2\ell  - 4}}{3}} \left( { - 2\pi i\nu } \right)} \\
& + i\frac{{e^{2\pi i\nu } }}{{\sqrt 3 }}\frac{2}{{3\pi }}\sum\limits_{\ell  = 0}^{L - 1} {d_{2\ell } e^{\frac{{2\left( {2\ell  + 1} \right)\pi i}}{3}} \sin \left( {\frac{{\left( {2\ell  + 1} \right)\pi }}{3}} \right)\frac{{\Gamma \left( {\frac{{2\ell  + 1}}{3}} \right)}}{{\nu ^{\frac{{2\ell  + 1}}{3}} }} \widehat T_{2M - \frac{{2\ell  - 4}}{3}} \left( {2\pi i\nu } \right)} \\
& + R_{N,M,K,L}^{\left( H \right)} \left( \nu  \right) .
\end{align*}
Taking $\nu = r e^{i \theta}$, the representation \eqref{eq77} becomes
\begin{gather}\label{eq78}
\begin{split}
R_{N,M,K,L}^{\left( H \right)} \left( \nu  \right) = \; & \left( { - 1} \right)^N \frac{{e^{ - \frac{\pi }{3}i} }}{{2\sqrt 3 \pi \left( {e^{i\theta } } \right)^{2N + \frac{1}{3}} }}\int_0^{ + \infty } {\frac{{\tau ^{2N - \frac{2}{3}} e^{ - 2\pi r\tau } }}{{1 + i\tau e^{ - i\theta } }}iR_K^{\left( H \right)} \left( {ir\tau } \right)d\tau } \\
& + \left( { - 1} \right)^N \frac{{e^{ - \frac{\pi }{3}i} }}{{2\sqrt 3 \pi \left( {e^{i\theta } } \right)^{2N + \frac{1}{3}} }}\int_0^{ + \infty } {\frac{{\tau ^{2N - \frac{2}{3}} e^{ - 2\pi r\tau } }}{{1 - i\tau e^{ - i\theta } }}iR_K^{\left( H \right)} \left( {ir\tau } \right)d\tau } \\
& + \left( { - 1} \right)^{M - 1} \frac{{e^{\frac{\pi }{3}i} }}{{2\sqrt 3 \pi \left( {e^{i\theta } } \right)^{2M + \frac{5}{3}} }}\int_0^{ + \infty } {\frac{{\tau ^{2M + \frac{2}{3}} e^{ - 2\pi r\tau } }}{{1 + i\tau e^{ - i\theta } }}iR_L^{\left( H \right)} \left( {ir\tau } \right)d\tau } \\
& + \left( { - 1} \right)^{M - 1} \frac{{e^{\frac{\pi }{3}i} }}{{2\sqrt 3 \pi \left( {e^{i\theta } } \right)^{2M + \frac{5}{3}} }}\int_0^{ + \infty } {\frac{{\tau ^{2M + \frac{2}{3}} e^{ - 2\pi r\tau } }}{{1 - i\tau e^{ - i\theta } }}iR_L^{\left( H \right)} \left( {ir\tau } \right)d\tau } .
\end{split}
\end{gather}
In Appendix \ref{appendixb} it is shown that
\begin{equation}\label{eq90}
\frac{{1 - \left( {s/r\tau } \right)^{\frac{4}{3}} }}{{1 - \left( {s/r\tau } \right)^2 }} = \frac{{1 - \left( {s/r} \right)^{\frac{4}{3}} }}{{1 - \left( {s/r} \right)^2 }} + \left( {\tau  - 1} \right)f\left( {r,\tau ,s} \right)
\end{equation}
for positive $r$, $\tau$ and $s$, with some $f\left(r,\tau ,s\right)$ satisfying $\left|f\left(r,\tau ,s\right)\right| \leq 2$. Using the integral formula \eqref{eq14}, $R_K^{\left( H \right)} \left( {ir\tau } \right)$ can be written as
\begin{align*}
R_K^{\left( H \right)} \left( {ir\tau } \right) = \; & \frac{1}{{\sqrt 3 \pi \left( {r\tau } \right)^{\frac{{2K + 1}}{3}} }}\int_0^{ + \infty } {s^{\frac{{2K - 2}}{3}} e^{ - 2\pi s} \frac{{1 - \left( {s/r\tau } \right)^{\frac{4}{3}} }}{{1 - \left( {s/r\tau } \right)^2 }} H_{is}^{\left( 1 \right)} \left( {is} \right)ds} \\
= \; & \frac{1}{{\sqrt 3 \pi \left( {r\tau } \right)^{\frac{{2K + 1}}{3}} }}\int_0^{ + \infty } {s^{\frac{{2K - 2}}{3}} e^{ - 2\pi s} \frac{{1 - \left( {s/r} \right)^{\frac{4}{3}} }}{{1 - \left( {s/r} \right)^2 }} H_{is}^{\left( 1 \right)} \left( {is} \right)ds} \\ & + \frac{{\tau  - 1}}{{\sqrt 3 \pi \left( {r\tau } \right)^{\frac{{2K + 1}}{3}} }}\int_0^{ + \infty } {s^{\frac{{2K - 2}}{3}} e^{ - 2\pi s} f\left( {r,\tau ,s} \right) H_{is}^{\left( 1 \right)} \left( {is} \right)ds} ,
\end{align*}
and similarly for $R_L^{\left( H \right)} \left( {ir\tau } \right)$. Noting that
\[
0< \frac{{1 - \left( {s/r} \right)^{\frac{4}{3}} }}{{1 - \left( {s/r} \right)^2 }} < 1
\]
for any positive $r$ and $s$, substitution into \eqref{eq78} yields the upper bound
\begin{align*}
\left| {R_{N,M,K,L}^{\left( H \right)} \left( \nu  \right)} \right| \le \; & \frac{1}{{3\pi }}\left| {d_{2K} } \right|\frac{{\Gamma \left( {\frac{{2K + 1}}{3}} \right)}}{{\left| \nu  \right|^{\frac{{2K + 1}}{3}} }}\left| {\frac{1}{{2\pi }}\int_0^{ + \infty } {\frac{{\tau ^{2N - \frac{{2K}}{3} - 1} e^{ - 2\pi r\tau } }}{{1 + i\tau e^{ - i\theta } }}d\tau } } \right| \\ & + \frac{1}{{3\pi ^2 }}\left| {d_{2K} } \right|\frac{{\Gamma \left( {\frac{{2K + 1}}{3}} \right)}}{{\left| \nu  \right|^{\frac{{2K + 1}}{3}} }}\int_0^{ + \infty } {\tau ^{2N - \frac{{2K}}{3} - 1} e^{ - 2\pi r\tau } \left| {\frac{{\tau  - 1}}{{\tau  - ie^{i\theta } }}} \right|d\tau } \\
& + \frac{1}{{3\pi }}\left| {d_{2K} } \right|\frac{{\Gamma \left( {\frac{{2K + 1}}{3}} \right)}}{{\left| \nu  \right|^{\frac{{2K + 1}}{3}} }}\left| {\frac{1}{{2\pi }}\int_0^{ + \infty } {\frac{{\tau ^{2N - \frac{{2K}}{3} - 1} e^{ - 2\pi r\tau } }}{{1 - i\tau e^{ - i\theta } }}d\tau } } \right| \\ & + \frac{1}{{3\pi ^2 }}\left| {d_{2K} } \right|\frac{{\Gamma \left( {\frac{{2K + 1}}{3}} \right)}}{{\left| \nu  \right|^{\frac{{2K + 1}}{3}} }}\int_0^{ + \infty } {\tau ^{2N - \frac{{2K}}{3} - 1} e^{ - 2\pi r\tau } \left| {\frac{{\tau  - 1}}{{\tau  + ie^{i\theta } }}} \right|d\tau } \\
& + \frac{1}{{3\pi }}\left| {d_{2L} } \right|\frac{{\Gamma \left( {\frac{{2L + 1}}{3}} \right)}}{{\left| \nu  \right|^{\frac{{2L + 1}}{3}} }}\left| {\frac{1}{{2\pi }}\int_0^{ + \infty } {\frac{{\tau ^{2M - \frac{{2L - 4}}{3} - 1} e^{ - 2\pi r\tau } }}{{1 + i\tau e^{ - i\theta } }}d\tau } } \right| \\ & + \frac{1}{{3\pi ^2 }}\left| {d_{2L} } \right|\frac{{\Gamma \left( {\frac{{2L + 1}}{3}} \right)}}{{\left| \nu  \right|^{\frac{{2L + 1}}{3}} }}\int_0^{ + \infty } {\tau ^{2M - \frac{{2L - 4}}{3} - 1} e^{ - 2\pi r\tau } \left| {\frac{{\tau  - 1}}{{\tau  - ie^{i\theta } }}} \right|d\tau } \\
& + \frac{1}{{3\pi }}\left| {d_{2L} } \right|\frac{{\Gamma \left( {\frac{{2L + 1}}{3}} \right)}}{{\left| \nu  \right|^{\frac{{2L + 1}}{3}} }}\left| {\frac{1}{{2\pi }}\int_0^{ + \infty } {\frac{{\tau ^{2M - \frac{{2L - 4}}{3} - 1} e^{ - 2\pi r\tau } }}{{1 - i\tau e^{ - i\theta } }}d\tau } } \right| \\ & + \frac{1}{{3\pi ^2 }}\left| {d_{2L} } \right|\frac{{\Gamma \left( {\frac{{2L + 1}}{3}} \right)}}{{\left| \nu  \right|^{\frac{{2L + 1}}{3}} }}\int_0^{ + \infty } {\tau ^{2M - \frac{{2L - 4}}{3} - 1} e^{ - 2\pi r\tau } \left| {\frac{{\tau  - 1}}{{\tau  + ie^{i\theta } }}} \right|d\tau } .
\end{align*}
As $\left| \left(\tau  - 1\right)/\left(\tau \pm i e^{i\theta}\right)  \right| \le 1$, we find that
\begin{align*}
\left| {R_{N,M,K,L}^{\left( H \right)} \left( \nu  \right)} \right| \le \; & \frac{1}{{3\pi }}\left| {d_{2K} } \right|\frac{{\Gamma \left( {\frac{{2K + 1}}{3}} \right)}}{{\left| \nu  \right|^{\frac{{2K + 1}}{3}} }}\left| {e^{ - 2\pi i\nu } \widehat T_{2N - \frac{{2K}}{3}} \left( { - 2\pi i\nu } \right)} \right| + \frac{1}{{3\pi }}\left| {d_{2K} } \right|\frac{{\Gamma \left( {\frac{{2K + 1}}{3}} \right)}}{{\left| \nu  \right|^{\frac{{2K + 1}}{3}} }}\left| {e^{2\pi i\nu } \widehat T_{2N - \frac{{2K}}{3}} \left( {2\pi i\nu } \right)} \right|\\
& + \frac{2}{{3\pi ^2 }}\left| {d_{2K} } \right|\frac{{\Gamma \left( {\frac{{2K + 1}}{3}} \right)\Gamma \left( {2N - \frac{{2K}}{3}} \right)}}{{\left( {2\pi } \right)^{2N - \frac{{2K}}{3}} \left| \nu  \right|^{2N+\frac{1}{3}} }} + \frac{1}{{3\pi }}\left| {d_{2L} } \right|\frac{{\Gamma \left( {\frac{{2L + 1}}{3}} \right)}}{{\left| \nu  \right|^{\frac{{2L + 1}}{3}} }}\left| {e^{ - 2\pi i\nu } \widehat T_{2M - \frac{{2L - 4}}{3}} \left( { - 2\pi i\nu } \right)} \right|\\
& + \frac{1}{{3\pi }}\left| {d_{2L} } \right|\frac{{\Gamma \left( {\frac{{2L + 1}}{3}} \right)}}{{\left| \nu  \right|^{\frac{{2L + 1}}{3}} }}\left| {e^{2\pi i\nu } \widehat T_{2M - \frac{{2L - 4}}{3}} \left( {2\pi i\nu } \right)} \right| + \frac{2}{{3\pi ^2 }}\left| {d_{2L} } \right|\frac{{\Gamma \left( {\frac{{2L + 1}}{3}} \right)\Gamma \left( {2M - \frac{{2L - 4}}{3}} \right)}}{{\left( {2\pi } \right)^{2M - \frac{{2L - 4}}{3}} \left| \nu  \right|^{2M+\frac{5}{3}} }}.
\end{align*}
By continuity, this bound holds in the closed sector $- \frac{\pi }{2} \le \arg \nu \le \frac{\pi }{2}$. Suppose that $
N = \pi \left| \nu  \right| + \rho$, $M = \pi \left| \nu  \right| + \sigma$ where $\rho$ and $\sigma$ are bounded. An application of Stirling's formula shows that
\[
\frac{2}{{3\pi ^2 }}\left| {d_{2K} } \right|\frac{{\Gamma \left( {\frac{{2K + 1}}{3}} \right)\Gamma \left( {2N - \frac{{2K}}{3}} \right)}}{{\left( {2\pi } \right)^{2N - \frac{{2K}}{3}} \left| \nu  \right|^{2N + \frac{1}{3}} }} = \mathcal{O}_{K,\rho } \left( {\frac{{e^{ - 2\pi \left| \nu  \right|} }}{{\left| \nu  \right|^{\frac{1}{2}} }}\left| {d_{2K} } \right|\frac{{\Gamma \left( {\frac{{2K + 1}}{3}} \right)}}{{\left| \nu  \right|^{\frac{{2K + 1}}{3}} }}} \right)
\]
and
\[
\frac{2}{{3\pi ^2 }}\left| {d_{2L} } \right|\frac{{\Gamma \left( {\frac{{2L + 1}}{3}} \right)\Gamma \left( {2M - \frac{{2L - 4}}{3}} \right)}}{{\left( {2\pi } \right)^{2M - \frac{{2L - 4}}{3}} \left| \nu  \right|^{2M + \frac{5}{3}} }} = \mathcal{O}_{L,\sigma } \left( {\frac{{e^{ - 2\pi \left| \nu  \right|} }}{{\left| \nu  \right|^{\frac{1}{2}} }}\left| {d_{2L} } \right|\frac{{\Gamma \left( {\frac{{2L + 1}}{3}} \right)}}{{\left| \nu  \right|^{\frac{{2L + 1}}{3}} }}} \right)
\]
as $\nu \to \infty$. Using Olver's estimation \eqref{eq71} we find
\[
\left| {e^{ \pm 2\pi i\nu } \widehat T_{2N - \frac{{2K}}{3}} \left( { \pm 2\pi i\nu } \right)} \right| = \mathcal{O}_{K,\rho } \left( {e^{ - 2\pi \left| \nu  \right|} } \right)
\]
and
\[
\left| {e^{ \pm 2\pi i\nu } \widehat T_{2M - \frac{{2L-4}}{3}} \left( { \pm 2\pi i\nu } \right)} \right| = \mathcal{O}_{L,\sigma } \left( {e^{ - 2\pi \left| \nu  \right|} } \right)
\]
for large $\nu$. Therefore, we have
\begin{equation}\label{eq79}
R_{N,M,K,L}^{\left( H \right)} \left( \nu  \right) = \mathcal{O}_{K,\rho } \left( {e^{ - 2\pi \left| \nu  \right|} \left| {d_{2K} } \right|\frac{{\Gamma \left( {\frac{{2K + 1}}{3}} \right)}}{{\left| \nu  \right|^{\frac{{2K + 1}}{3}} }}} \right) + \mathcal{O}_{L,\sigma } \left( {e^{ - 2\pi \left| \nu  \right|} \left| {d_{2L} } \right|\frac{{\Gamma \left( {\frac{{2L + 1}}{3}} \right)}}{{\left| \nu  \right|^{\frac{{2L + 1}}{3}} }}} \right)
\end{equation}
as $\nu \to \infty$ in the sector $- \frac{\pi }{2} \le \arg \nu \le \frac{\pi }{2}$.

Next, we consider the sector $\frac{\pi}{2} \leq \arg \nu \le \frac{3\pi}{2}$. Rotating the path of integration in \eqref{eq77} and applying the residue theorem gives
\begin{align*}
R_{N,M,K,L}^{\left( H \right)} \left( \nu  \right) = \; & i\frac{{e^{2\pi i\nu } }}{{\sqrt 3 }}R_K^{\left( H \right)} \left( \nu  \right) + \left( { - 1} \right)^N \frac{{e^{ - \frac{\pi }{3}i} }}{{\sqrt 3 \pi \nu ^{2N + \frac{1}{3}} }}\int_0^{ + \infty } {\frac{{t^{2N - \frac{2}{3}} e^{ - 2\pi t} }}{{1 + \left( {t/\nu } \right)^2 }}iR_K^{\left( H \right)} \left( {it} \right)dt} \\
& - i\frac{{e^{2\pi i\nu } }}{{\sqrt 3 }}R_L^{\left( H \right)} \left( \nu  \right) + \left( { - 1} \right)^{M - 1} \frac{{e^{\frac{\pi }{3}i} }}{{\sqrt 3 \pi \nu ^{2M + \frac{5}{3}} }}\int_0^{ + \infty } {\frac{{t^{2M + \frac{2}{3}} e^{ - 2\pi t} }}{{1 + \left( {t/\nu } \right)^2 }}iR_L^{\left( H \right)} \left( {it} \right)dt} \\
 = & - \overline {R_{N,M,K,L}^{\left( H \right)} \left( {\nu e^{ - \pi i} } \right)}  + i\frac{{e^{2\pi i\nu } }}{{\sqrt 3 }}R_K^{\left( H \right)} \left( \nu  \right) - i\frac{{e^{2\pi i\nu } }}{{\sqrt 3 }}R_L^{\left( H \right)} \left( \nu  \right)
\end{align*}
for $\frac{\pi}{2} < \arg \nu  < \frac{3\pi}{2}$. It follows that when $K=L$, the bound \eqref{eq79} remains valid in the wider sector $- \frac{\pi }{2} \le \arg \nu \le \frac{3\pi }{2}$. Otherwise, we have
\[
\left| {R_{N,M,K,L}^{\left( H \right)} \left( \nu  \right)} \right| \le \left| {R_{N,M,K,L}^{\left( H \right)} \left( {\nu e^{ - \pi i} } \right)} \right| + \frac{{e^{ - 2\pi \Im \left( \nu  \right)} }}{{\sqrt 3 }}\left| {R_K^{\left( H \right)} \left( \nu  \right)} \right| + \frac{{e^{ - 2\pi \Im \left( \nu  \right)} }}{{\sqrt 3 }}\left| {R_L^{\left( H \right)} \left( \nu  \right)} \right|
\]
and therefore, by the bounds from Subsection \ref{bounds} and \eqref{eq79},
\[
R_{N,M,K,L}^{\left( H \right)} \left( \nu  \right) = \mathcal{O}_{K,\rho} \left( {e^{ - 2\pi \Im \left( \nu  \right)} \left| {d_{2K} } \right|\frac{{\Gamma \left( {\frac{{2K + 1}}{3}} \right)}}{{\left| \nu  \right|^{\frac{{2K + 1}}{3}} }}} \right) + \mathcal{O}_{L,\sigma} \left( {e^{ - 2\pi \Im \left( \nu  \right)} \left| {d_{2L} } \right|\frac{{\Gamma \left( {\frac{{2L + 1}}{3}} \right)}}{{\left| \nu  \right|^{\frac{{2L + 1}}{3}} }}} \right)
\]
as $\nu \to \infty$ in the sector $\frac{\pi}{2} \le \arg \nu  \le  \frac{3\pi}{2}$.

Finally, we consider the case $-\frac{3\pi}{2} \le \arg \nu \le -\frac{\pi}{2}$. Rotating the path of integration in \eqref{eq77} and applying the residue theorem gives
\begin{gather}\label{eq84}
\begin{split}
R_{N,M,K,L}^{\left( H \right)} \left( \nu  \right) = & - ie^{\frac{\pi }{3}i} \frac{{e^{ - 2\pi i\nu } }}{{\sqrt 3 }}R_K^{\left( H \right)} \left( {\nu e^{\pi i} } \right) + \left( { - 1} \right)^N \frac{{e^{ - \frac{\pi }{3}i} }}{{\sqrt 3 \pi \nu ^{2N + \frac{1}{3}} }}\int_0^{ + \infty } {\frac{{t^{2N - \frac{2}{3}} e^{ - 2\pi t} }}{{1 + \left( {t/\nu } \right)^2 }}iR_K^{\left( H \right)} \left( {it} \right)dt} \\
& + ie^{ - \frac{\pi }{3}i} \frac{{e^{ - 2\pi i\nu } }}{{\sqrt 3 }}R_L^{\left( H \right)} \left( { \nu e^{\pi i} } \right) + \left( { - 1} \right)^{M - 1} \frac{{e^{\frac{\pi }{3}i} }}{{\sqrt 3 \pi \nu ^{2M + \frac{5}{3}} }}\int_0^{ + \infty } {\frac{{t^{2M + \frac{2}{3}} e^{ - 2\pi t} }}{{1 + \left( {t/\nu } \right)^2 }}iR_L^{\left( H \right)} \left( {it} \right)dt} \\
= & - ie^{\frac{\pi }{3}i} \frac{{e^{ - 2\pi i\nu } }}{{\sqrt 3 }}R_K^{\left( H \right)} \left( { \nu e^{\pi i} } \right) + \left( { - 1} \right)^N \frac{1}{{\sqrt 3 \pi \left( {\nu e^{\pi i} } \right)^{2N + \frac{1}{3}} }}\int_0^{ + \infty } {\frac{{t^{2N - \frac{2}{3}} e^{ - 2\pi t} }}{{1 + \left( {t/\nu e^{\pi i} } \right)^2 }}iR_K^{\left( H \right)} \left( {it} \right)dt} \\
& + ie^{ - \frac{\pi }{3}i} \frac{{e^{ - 2\pi i\nu } }}{{\sqrt 3 }}R_L^{\left( H \right)} \left( { \nu e^{\pi i} } \right) + \left( { - 1} \right)^{M - 1} \frac{1}{{\sqrt 3 \pi \left( {\nu e^{\pi i} } \right)^{2M + \frac{5}{3}} }}\int_0^{ + \infty } {\frac{{t^{2M + \frac{2}{3}} e^{ - 2\pi t} }}{{1 + \left( {t/\nu e^{\pi i} } \right)^2 }}iR_L^{\left( H \right)} \left( {it} \right)dt} 
\end{split}
\end{gather}
for $-\frac{3\pi}{2} < \arg \nu < -\frac{\pi}{2}$. It is easy to see that the sum of the two integrals has the order of magnitude given in the right-hand side of \eqref{eq79}. Since
\[
\left| { - ie^{\frac{\pi }{3}i} \frac{{e^{ - 2\pi i\nu } }}{{\sqrt 3 }}R_K^{\left( H \right)} \left( {\nu e^{\pi i} } \right) + ie^{ - \frac{\pi }{3}i} \frac{{e^{ - 2\pi i\nu } }}{{\sqrt 3 }}R_L^{\left( H \right)} \left( { \nu e^{\pi i} } \right)} \right| \le \frac{{e^{2\pi \Im \left( \nu  \right)} }}{{\sqrt 3 }}\left| {R_K^{\left( H \right)} \left( { \nu e^{\pi i} } \right)} \right| + \frac{{e^{2\pi \Im \left( \nu  \right)} }}{{\sqrt 3 }}\left| {R_L^{\left( H \right)} \left( { \nu e^{\pi i} } \right)} \right|,
\]
applying the error bounds given in Subsection \ref{bounds}, produces
\[
R_{N,M,K,L}^{\left( H \right)} \left( \nu  \right) = \mathcal{O}_{K,\rho} \left( {e^{ 2\pi \Im \left( \nu  \right)} \left| {d_{2K} } \right|\frac{{\Gamma \left( {\frac{{2K + 1}}{3}} \right)}}{{\left| \nu  \right|^{\frac{{2K + 1}}{3}} }}} \right) + \mathcal{O}_{L,\sigma} \left( {e^{ 2\pi \Im \left( \nu  \right)} \left| {d_{2L} } \right|\frac{{\Gamma \left( {\frac{{2L + 1}}{3}} \right)}}{{\left| \nu  \right|^{\frac{{2L + 1}}{3}} }}} \right)
\]
as $\nu \to \infty$ in the sector $-\frac{3\pi}{2} \le \arg \nu \le -\frac{\pi}{2}$.

\subsection{Stokes phenomenon and Berry's transition}

\subsubsection{Case (i): $x>1$} We study the Stokes phenomenon related to Debye's expansion for $H_\nu ^{\left( 1 \right)} \left( {\nu \sec \beta } \right)$ occurring when $\arg \nu$ passes through the value $-\frac{\pi}{2}$. From \eqref{eq57} we have
\begin{align*}
R_N^{\left( H \right)} \left( {\nu ,\beta } \right) = R_{N,N,0,0}^{\left( H \right)} \left( {\nu ,\beta } \right) = & - ie^{ - 2i\nu \left( {\tan \beta  - \beta } \right)} R_0^{\left( H \right)} \left( { \nu e^{\pi i},\beta } \right) - ie^{ - 2i\nu \left( {\tan \beta  - \beta  + \pi } \right)} R_0^{\left( H \right)} \left( { \nu e^{\pi i},\beta } \right)
\\ & + \frac{1}{{2\left( {2\pi \cot \beta } \right)^{\frac{1}{2}} \left( {i\nu } \right)^N }}\int_0^{ + \infty } {\frac{{t^{N - \frac{1}{2}} e^{ - t\left( {\tan \beta  - \beta } \right)} }}{{1 + it/\nu }}\left( {1 + e^{ - 2\pi t} } \right)iH_{it}^{\left( 1 \right)} \left( {it\sec \beta } \right)dt} 
\end{align*}
when $-\frac{3\pi}{2} < \arg \nu < -\frac{\pi}{2}$. Then, by Debye's expansions, we infer that
\[
R_N^{\left( H \right)} \left( {\nu ,\beta } \right) \sim  - ie^{ - 2i\nu \left( {\tan \beta  - \beta } \right)} \sum\limits_{k = 0}^\infty  {\frac{{U_k \left( {i\cot \beta } \right)}}{{\nu ^k }}}  - ie^{ - 2i\nu \left( {\tan \beta  - \beta  + \pi } \right)} \sum\limits_{\ell  = 0}^\infty  {\frac{{U_\ell  \left( {i\cot \beta } \right)}}{{\nu ^\ell  }}}  + \sum\limits_{n = N}^\infty  {\left( { - 1} \right)^n \frac{{U_n \left( {i\cot \beta } \right)}}{{\nu ^n }}} 
\]
as $\nu \to \infty$ in the sector $-\frac{3\pi}{2} < \arg \nu < -\frac{\pi}{2}$. Therefore, as the line $\arg \nu = -\frac{\pi}{2}$ is crossed, the two additional series
\begin{equation}\label{eq81}
- ie^{ - 2i\nu \left( {\tan \beta  - \beta } \right)} \sum\limits_{k = 0}^\infty  {\frac{U_k \left( {i\cot \beta } \right)}{\nu ^k}}  - ie^{ - 2i\nu \left( {\tan \beta  - \beta  + \pi } \right)} \sum\limits_{\ell  = 0}^\infty  {\frac{U_\ell  \left( {i\cot \beta } \right)}{\nu ^\ell}}
\end{equation}
appear in the asymptotic expansion of $H_\nu ^{\left( 1 \right)} \left( {\nu \sec \beta } \right)$ beside the original series $\sum\nolimits_{n = 0}^\infty  {\left( { - 1} \right)^n U_n \left( {i\cot \beta } \right)/\nu ^n }$. We have encountered a Stokes phenomenon with Stokes line $\arg \nu = -\frac{\pi}{2}$.

In the important papers \cite{Berry3, Berry2}, Berry provided a new interpretation of the Stokes phenomenon; he found that assuming optimal truncation, the transition between compound asymptotic expansions is of Error function type, thus yielding a smooth, although very rapid, transition as a Stokes line is crossed. Similar results for general linear or nonlinear ordinary differential equations were proved by Costin and Kruskal \cite{Costin}.

Using the exponentially improved expansion given in Theorem \ref{thm3}, we show that Debye's expansion exhibits the Berry transition between the two asymptotic series across the Stokes line $\arg \nu = -\frac{\pi}{2}$. More precisely, we shall find that the first few terms of the two series in \eqref{eq81} ``emerge" in a rapid and smooth way as $\arg \nu$ decreases through $-\frac{\pi}{2}$.

From Theorem \ref{thm3}, we conclude that if $N \approx 2\left| \nu  \right|\left( {\tan \beta  - \beta } \right)$, $M \approx 2\left| \nu  \right|\left( {\tan \beta  - \beta  + \pi } \right)$ then for large $\nu$, $ - \pi < \arg \nu \leq \frac{3\pi}{2}$, we have
\begin{multline*}
H_\nu ^{\left( 1 \right)} \left( {\nu \sec \beta } \right) \approx \frac{{e^{i\nu \left( {\tan \beta  - \beta } \right) - \frac{\pi }{4}i} }}{{\left( {\frac{1}{2}\nu \pi \tan \beta } \right)^{\frac{1}{2}} }}\left( {\sum\limits_{n = 0}^{N - 1} {\left( { - 1} \right)^n \frac{{U_n \left( {i\cot \beta } \right)}}{{\nu ^n }}}  + \sum\limits_{m = N}^{M - 1} {\left( { - 1} \right)^m \frac{{\widetilde U_m \left( {i\cot \beta } \right)}}{{\nu ^m }}}}\right. \\  + ie^{ - 2i\nu \left( {\tan \beta  - \beta } \right)} \sum\limits_{k = 0} {\frac{{U_k \left( {i\cot \beta } \right)}}{{\nu ^k }}\widehat T_{N - k} \left( { - 2i\nu \left( {\tan \beta  - \beta } \right)} \right)} \\ \left.{ + ie^{ - 2i\nu \left( {\tan \beta  - \beta  + \pi } \right)} \sum\limits_{\ell  = 0} {\frac{{U_\ell  \left( {i\cot \beta } \right)}}{{\nu ^\ell  }}\widehat T_{M - \ell } \left( { - 2i\nu \left( {\tan \beta  - \beta  + \pi } \right)} \right)} } \right),
\end{multline*}
where $\sum\nolimits_{k = 0}$ and $\sum\nolimits_{\ell = 0}$ mean that the sums are restricted to the first few terms of the series.

Under the above assumptions on $N$ and $M$, from \eqref{eq82} and \eqref{eq83}, the Terminant functions have the asymptotic behaviour
\begin{align*}
\widehat T_{N - k} \left( { - 2i\nu \left( {\tan \beta  - \beta } \right)} \right) & \sim  - \frac{1}{2} + \frac{1}{2}\mathop{\text{erf}} \left( { - \overline {c\left( {\frac{\pi }{2} - \theta } \right)} \sqrt {\left| \nu  \right|\left( {\tan \beta  - \beta } \right)} } \right)\\ & \sim  - \frac{1}{2} + \frac{1}{2}\mathop{\text{erf}} \left( {\left( {\theta  + \frac{\pi }{2}} \right)\sqrt {\left| \nu  \right|\left( {\tan \beta  - \beta } \right)} } \right),
\end{align*}
\begin{align*}
\widehat T_{M - \ell } \left( { - 2i\nu \left( {\tan \beta  - \beta  + \pi } \right)} \right) & \sim  - \frac{1}{2} + \frac{1}{2}\mathop{\text{erf}} \left( { - \overline {c\left( {\frac{\pi }{2} - \theta } \right)} \sqrt {\left| \nu  \right|\left( {\tan \beta  - \beta  + \pi } \right)} } \right)\\ & \sim  - \frac{1}{2} + \frac{1}{2}\mathop{\text{erf}} \left( {\left( {\theta  + \frac{\pi }{2}} \right)\sqrt {\left| \nu  \right|\left( {\tan \beta  - \beta  + \pi } \right)} } \right)
\end{align*}
provided that $\arg \nu = \theta$ is close to $-\frac{\pi}{2}$, $\nu$ is large and $k,\ell$ are small in comparison with $N$ and $M$. Therefore, when $\theta  >  - \frac{\pi }{2}$, the Terminant functions are exponentially small; for $\theta  =  -\frac{\pi }{2}$, they are asymptotically $-\frac{1}{2}$ up to an exponentially small error; and when $\theta  <  - \frac{\pi}{2}$, the Terminant functions are asymptotic to $-1$ with an exponentially small error. Thus, the transition across the Stokes line $\arg \nu = -\frac{\pi}{2}$ is effected rapidly and smoothly.

\subsubsection{Case (ii): $x=1$} The analysis of the Stokes phenomenon for the asymptotic expansion of $H_\nu ^{\left( 1 \right)} \left( \nu \right)$ is similar to case $x>1$. From \eqref{eq84} we infer that
\begin{align*}
& R_{3N}^{\left( H \right)} \left( \nu  \right) = R_{N,N,0,0}^{\left( H \right)} \left( \nu  \right) =  ie^{\frac{\pi }{3}i} \frac{{e^{ - 2\pi i\nu } }}{{\sqrt 3 }}H_\nu ^{\left( 2 \right)} \left( \nu \right) - ie^{ - \frac{\pi }{3}i} \frac{{e^{ - 2\pi i\nu } }}{{\sqrt 3 }}H_\nu ^{\left( 2 \right)} \left( \nu  \right) \\
& + \frac{{\left( { - 1} \right)^{N + 1} }}{{\sqrt 3 \pi \nu ^{2N + \frac{1}{3}} }}\int_0^{ + \infty } {t^{2N - \frac{2}{3}} e^{ - 2\pi t} e^{\frac{{2\left( {6N + 1} \right)\pi i}}{3}} \frac{{1 + \left( {t/\nu } \right)^{\frac{2}{3}} e^{\frac{\pi }{3}i} }}{{\left( {1 + \left( {t/\nu } \right)^{\frac{2}{3}} e^{\frac{{2\pi i}}{3}} } \right)\left( {1 + \left( {t/\nu } \right)^{\frac{2}{3}} } \right)}}iH_{it}^{\left( 1 \right)} \left( {it} \right)dt} 
\end{align*}
when $-\frac{3\pi}{2} < \arg \nu < -\frac{\pi}{2}$. Then, from Debye's expansions, we obtain the compound asymptotic series
\begin{align*}
R_{3N}^{\left( H \right)} \left( \nu  \right) \sim & - ie^{\frac{\pi }{3}i} \frac{{e^{ - 2\pi i\nu } }}{{\sqrt 3 }}\frac{2}{{3\pi }}\sum\limits_{k = 0}^\infty  {d_{2k} e^{ - \frac{{2\left( {2k + 1} \right)\pi i}}{3}} \sin \left( {\frac{{\left( {2k + 1} \right)\pi }}{3}} \right)\frac{{\Gamma \left( {\frac{{2k + 1}}{3}} \right)}}{{\nu ^{\frac{{2k + 1}}{3}} }}} \\
& + ie^{ - \frac{\pi }{3}i} \frac{{e^{ - 2\pi i\nu } }}{{\sqrt 3 }}\frac{2}{{3\pi }}\sum\limits_{\ell  = 0}^\infty  {d_{2\ell } e^{ - \frac{{2\left( {2\ell  + 1} \right)\pi i}}{3}} \sin \left( {\frac{{\left( {2\ell  + 1} \right)\pi }}{3}} \right)\frac{{\Gamma \left( {\frac{{2\ell  + 1}}{3}} \right)}}{{\nu ^{\frac{{2\ell  + 1}}{3}} }}} \\
& - \frac{2}{{3\pi }}\sum\limits_{n = 3N}^\infty  {d_{2n} e^{\frac{{2\left( {2n + 1} \right)\pi i}}{3}} \sin \left( {\frac{{\left( {2n + 1} \right)\pi }}{3}} \right)\frac{{\Gamma \left( {\frac{{2n + 1}}{3}} \right)}}{{\nu ^{\frac{{2n + 1}}{3}} }}} 
\end{align*}
as $\nu \to \infty$ in the sector $-\frac{3\pi}{2} < \arg \nu < -\frac{\pi}{2}$. Hence, as $\arg \nu$ decreases through the value $-\frac{\pi}{2}$, the two additional series
\begin{gather}\label{eq85}
\begin{split}
& - ie^{\frac{\pi }{3}i} \frac{{e^{ - 2\pi i\nu } }}{{\sqrt 3 }}\frac{2}{{3\pi }}\sum\limits_{k = 0}^\infty  {d_{2k} e^{ - \frac{{2\left( {2k + 1} \right)\pi i}}{3}} \sin \left( {\frac{{\left( {2k + 1} \right)\pi }}{3}} \right)\frac{{\Gamma \left( {\frac{{2k + 1}}{3}} \right)}}{{\nu ^{\frac{{2k + 1}}{3}} }}} \\
& + ie^{ - \frac{\pi }{3}i} \frac{{e^{ - 2\pi i\nu } }}{{\sqrt 3 }}\frac{2}{{3\pi }}\sum\limits_{\ell  = 0}^\infty  {d_{2\ell } e^{ - \frac{{2\left( {2\ell  + 1} \right)\pi i}}{3}} \sin \left( {\frac{{\left( {2\ell  + 1} \right)\pi }}{3}} \right)\frac{{\Gamma \left( {\frac{{2\ell  + 1}}{3}} \right)}}{{\nu ^{\frac{{2\ell  + 1}}{3}} }}}
\end{split}
\end{gather}
appear in the asymptotic expansion of $H_\nu ^{\left( 1 \right)} \left( \nu \right)$. With the aid of the exponentially improved expansion given in Theorem \ref{thm4}, we shall find that the asymptotic series of $H_\nu ^{\left( 1 \right)} \left( \nu \right)$ shows the Berry transition property: the two series in \eqref{eq85} ``emerge" in a rapid and smooth way as the Stokes line $\arg \nu = -\frac{\pi}{2}$ is crossed.

Let us assume that $M, N \approx \pi \left| \nu \right|$. Under these conditions, Olver's estimation \eqref{eq71} gives that
\[
e^{2\pi i\nu } \widehat T_{2N - \frac{{2k}}{3}} \left( {2\pi i\nu } \right) = \mathcal{O}\left( {e^{ - 2\pi \left| \nu  \right|} } \right)
\]
and
\[
e^{2\pi i\nu } \widehat T_{2M - \frac{{2\ell  - 4}}{3}} \left( {2\pi i\nu } \right) = \mathcal{O}\left( {e^{ - 2\pi \left| \nu  \right|} } \right)
\]
as $\nu \to \infty$ and $-\frac{3\pi}{2} < \arg \nu < \frac{\pi}{2}$. Therefore, from Theorem \ref{thm4}, we deduce that for large $\nu$, $ - \pi < \arg \nu < \frac{\pi}{2}$, we have
\begin{align*}
H_\nu ^{\left( 1 \right)} \left( \nu  \right) \approx  \; & \frac{{e^{ - \frac{\pi}{3}i} }}{{\sqrt 3 \pi \nu ^{\frac{1}{3}} }}\sum\limits_{n = 0}^{N - 1} {d_{6n} \frac{{\Gamma \left( {\frac{{6n + 1}}{3}} \right)}}{{\nu ^{2n} }}}  - \frac{{e^{\frac{\pi}{3}i} }}{{\sqrt 3 \pi \nu ^{\frac{5}{3}} }}\sum\limits_{m = 0}^{M - 1} {d_{6m + 4} \frac{{\Gamma \left( {\frac{{6m + 5}}{3}} \right)}}{{\nu ^{2m} }}} \\
& + ie^{\frac{\pi }{3}i} \frac{{e^{ - 2\pi i\nu } }}{{\sqrt 3 }}\frac{2}{{3\pi }}\sum\limits_{k = 0} {d_{2k} e^{ - \frac{{2\left( {2k + 1} \right)\pi i}}{3}} \sin \left( {\frac{{\left( {2k + 1} \right)\pi }}{3}} \right)\frac{{\Gamma \left( {\frac{{2k + 1}}{3}} \right)}}{{\nu ^{\frac{{2k + 1}}{3}} }}e^{2\pi i \frac{{2k}}{3}} \widehat T_{2N - \frac{{2k}}{3}} \left( { - 2\pi i\nu } \right)} \\
& - ie^{ - \frac{\pi }{3}i} \frac{{e^{ - 2\pi i\nu } }}{{\sqrt 3 }}\frac{2}{{3\pi }}\sum\limits_{\ell  = 0} {d_{2\ell } e^{ - \frac{{2\left( {2\ell  + 1} \right)\pi i}}{3}} \sin \left( {\frac{{\left( {2\ell  + 1} \right)\pi }}{3}} \right)\frac{{\Gamma \left( {\frac{{2\ell  + 1}}{3}} \right)}}{{\nu ^{\frac{{2\ell  + 1}}{3}} }}e^{2\pi i \frac{{2\ell  - 4}}{3}} \widehat T_{2M - \frac{{2\ell  - 4}}{3}} \left( { - 2\pi i\nu } \right)},
\end{align*}
where, as before, $\sum\nolimits_{k = 0}$ and $\sum\nolimits_{\ell = 0}$ mean that the sums are restricted to the first few terms of the series.

Since $M, N \approx \pi \left| \nu \right|$, from \eqref{eq82} and \eqref{eq83}, the normalised Terminant functions have the asymptotic behaviour
\[
e^{2\pi i\frac{{2k}}{3}} \widehat T_{2N - \frac{{2k}}{3}} \left( { - 2\pi i\nu } \right) \sim - \frac{1}{2} + \frac{1}{2}\mathop{\text{erf}}\left( { - \overline {c\left( {\frac{\pi }{2} - \theta } \right)} \sqrt {\pi \left| \nu  \right|} } \right) \sim - \frac{1}{2} + \frac{1}{2}\mathop{\text{erf}}\left( {\left( {\theta  + \frac{\pi }{2}} \right)\sqrt {\pi \left| \nu  \right|} } \right),
\]
\[
e^{2\pi i\frac{{2\ell  - 4}}{3}} \widehat T_{2M - \frac{{2\ell  - 4}}{3}} \left( { - 2\pi i\nu } \right) \sim  - \frac{1}{2} + \frac{1}{2}\mathop{\text{erf}}\left( { - \overline {c\left( {\frac{\pi }{2} - \theta } \right)} \sqrt {\pi \left| \nu  \right|} } \right) \sim  - \frac{1}{2} + \frac{1}{2}\mathop{\text{erf}}\left( {\left( {\theta  + \frac{\pi }{2}} \right)\sqrt {\pi \left| \nu  \right|} } \right)
\]
under the conditions that $\arg \nu = \theta$ is close to $-\frac{\pi}{2}$, $\nu$ is large and $k,\ell$ are small compared to $N$ and $M$. Thus, when $\theta  >  - \frac{\pi }{2}$, the normalised Terminant functions are exponentially small; for $\theta  =  -\frac{\pi }{2}$, they are asymptotic to $-\frac{1}{2}$ with an exponentially small error; and when $\theta  <  - \frac{\pi}{2}$, the normalised Terminant functions are asymptotically $-1$ up to an exponentially small error. Thus, the transition through the Stokes line $\arg \nu = -\frac{\pi}{2}$ is carried out rapidly and smoothly.

\section{Discussion}\label{section6}

In this paper, we have discussed in detail the large order and argument asymptotics of the Hankel functions $H_\nu^{\left(1\right)}\left(\nu x\right)$, $H_\nu^{\left(2\right)}\left(\nu x\right)$ and the Bessel functions $J_\nu \left(\nu x\right)$, $Y_\nu \left(\nu x\right)$ when $x\geq 1$. As for the case $0 < x < 1$, we ran into trouble with the adjacent saddles when conducting the research. It turned out that there are infinitely many adjacent saddle points and the critical phase of $\nu$ corresponding to an adjacent saddle depends on $x$ and the adjacent saddle point, making the analysis hopeless to carry out. This is in agreement with the observation of Uchiyama \cite{Uchiyama} that he made about the structure of the adjacent saddles during the resurgence analysis of the Bessel function of the third kind $K_\nu\left(\nu z\right)$ when $z$ lies in a certain bounded domain of the complex plane, especially when $0 < \Im\left(z\right) < 1$.

Nevertheless, there are some results in the literature concerning approximations for the late coefficients, exponential improvement and error bounds. With the notation $x= \mathop{\text{sech}}\alpha$, $\alpha>0$, Debye's expansions for the Bessel functions can be written as
\begin{equation}\label{eq86}
J_\nu  \left( {\nu \mathop{\text{sech}}\alpha } \right) \sim \frac{{e^{\nu \left( {\tanh \alpha  - \alpha } \right)} }}{{\left( {2\pi \nu \tanh \alpha } \right)^{\frac{1}{2}} }}\sum\limits_{n = 0}^\infty  {\frac{{U_n \left( {\coth \alpha } \right)}}{{\nu ^n }}} ,
\end{equation}
\begin{equation}\label{eq87}
Y_\nu  \left( {\nu \mathop{\text{sech}}\alpha } \right) \sim  - \frac{{e^{ \nu \left( {\alpha - \tanh \alpha } \right)} }}{{\left( {\frac{1}{2}\pi \nu \tanh \alpha } \right)^{\frac{1}{2}} }}\sum\limits_{n = 0}^\infty  {\left( { - 1} \right)^n \frac{{U_n \left( {\coth \alpha } \right)}}{{\nu ^n }}} 
\end{equation}
when $\nu \to \infty$, provided that $\left|\arg \nu\right| < \pi$ (see \cite[p. 231]{NIST}). Based on Darboux's method, Dingle \cite[p. 168]{Dingle} gave a formal asymptotic expansion for the coefficients $U_n \left( {\coth \alpha } \right)$ when $n$ is large. His result, in our notation, may be written as
\begin{equation}\label{eq88}
U_n \left( {\coth \alpha } \right) \approx \frac{{\left( { - 1} \right)^n \Gamma \left( n \right)}}{{2\pi \left( {2\left( {\alpha  - \tanh \alpha } \right)} \right)^n }}\sum\limits_{m = 0}^\infty  {\left( {2\left( {\alpha  - \tanh \alpha } \right)} \right)^m U_m \left( {\coth \alpha } \right)\frac{{\Gamma \left( {n - m} \right)}}{{\Gamma \left( n \right)}}} .
\end{equation}
Numerical calculations indicate that this approximation is correct if it is truncated to the first few terms. As an application of this asymptotic expansion, using Borel summation, Dingle \cite[p. 482]{Dingle} derived the resurgence formulas
\[
U_n \left( {\coth \alpha } \right) = \frac{{\left( { - 1} \right)^n }}{{\left( {2\pi \coth \alpha } \right)^{\frac{1}{2}} }}\int_0^{ + \infty } {t^{n - \frac{1}{2}} e^{ - t\left( {\alpha  - \tanh \alpha } \right)} J_t \left( t \mathop{\text{sech}} \alpha  \right)dt} 
\]
and
\[
J_\nu  \left( {\nu \mathop{\text{sech}}\alpha } \right) = \frac{{e^{\nu \left( {\tanh \alpha  - \alpha } \right)} }}{{2\pi \nu ^{\frac{1}{2}} }}\int_0^{ + \infty } {\frac{{t^{ - \frac{1}{2}} e^{ - t\left( {\alpha  - \tanh \alpha } \right)} }}{{1 + t/\nu }}J_t \left( {t\mathop{\text{sech}} \alpha } \right)dt} .
\]
Unfortunately, numerical computations prove these formulas to be incorrect. Using his formal theory of terminants, he gave exponentially improved versions of \eqref{eq86} and \eqref{eq87} \cite[p. 467--468]{Dingle}. The expansion \eqref{eq86} can be rearranged to
\[
J_\nu  \left( {\nu \mathop{\text{sech}}\alpha } \right) \sim \frac{{e^{\nu \left( {\tanh \alpha  - \alpha } \right)} }}{{\left( {2\pi \nu \tanh \alpha } \right)^{\frac{1}{2}} }}\sum\limits_{n = 0}^\infty  {\frac{{C_n \left( \nu  \right)}}{{\left( {\nu \tanh \alpha } \right)^n }}} 
\]
where $C_n \left( \nu  \right)$ is a polynomial in $\nu^2$ of degree $\left\lfloor {n/3} \right\rfloor$. For $\nu, \alpha >0$ and $N\geq 1$, define the remainder $M_N \left( {\nu ,\alpha } \right)$ by
\[
J_\nu  \left( {\nu \mathop{\text{sech}}\alpha } \right) = \frac{{e^{\nu \left( {\tanh \alpha  - \alpha } \right)} }}{{\left( {2\pi \nu \tanh \alpha } \right)^{\frac{1}{2}} }}\left( {\sum\limits_{n = 0}^{N - 1} {\frac{{C_n \left( \nu  \right)}}{{\left( {\nu \tanh \alpha } \right)^n }}}  + M_N \left( {\nu ,\alpha } \right)} \right) .
\]
It was showed by van Veen \cite{van Veen} that
\begin{equation}\label{eq89}
\left| M_N \left( {\nu ,\alpha } \right) \right| \le \frac{{\left| {C_N \left( {i\nu } \right)} \right|}}{{\left( {\nu \tanh \alpha } \right)^N }} \begin{cases} 2^N \left( {2N + 2} \right) & \; \text{ for all } \; \nu >0 \\ 2^{N + 1} \left( {1 - \frac{{2N - 1}}{{4\nu \tanh \alpha }}} \right)^{ - \frac{1}{2}}  & \; \text{ if } \; 2N - 1 < 4\nu \tanh \alpha, \end{cases}
\end{equation}
provided that $N \equiv 0 \mod 3$. We remark that his notations slightly differ from ours.

A possible direction for further research is to find a rigorous proof and form of Dingle's late coefficient formula \eqref{eq88}. Another interesting problem is the extension of van Veen's error bounds \eqref{eq89} to every $N\geq 1$ and to complex values of $\nu$. The resurgence analysis of the large $\nu$ asymptotics of the functions $H_\nu ^{\left( 1 \right) \prime}  \left( {\nu x} \right)$, $H_\nu ^{\left( 2 \right) \prime}  \left( {\nu x} \right)$ and $J_\nu' \left(\nu x\right)$, $Y_\nu' \left(\nu x\right)$ is another possible direction of future research.

\section*{Acknowledgement} I would like to thank the two anonymous referees for their constructive and helpful comments and suggestions on the manuscript.

\appendix

\section{Computation of the coefficients $U_n\left(i \cot \beta\right)$ and $d_{2n}$}\label{appendixa}

In this appendix we collect some formulas for the computation of the coefficients that appear in Debye's expansions.

\subsection{The coefficients $U_n\left(i \cot \beta\right)$}

It is known that $U_n \left( {i\cot \beta } \right) = \left[U_n \left( x \right)\right]_{x = i\cot \beta }$ where $U_n\left(x\right)$ is a polynomial in $x$ of degree $3n$. We consider these polynomials. They can be generated by the following recurrence
\begin{equation}\label{eq74}
U_n \left( x \right) = \frac{1}{2}x^2 \left( {1 - x^2 } \right)U'_{n-1} \left( x \right) + \frac{1}{8}\int_0^x {\left( {1 - 5t^2 } \right)U_{n-1} \left( t \right)dt} 
\end{equation}
for $n \ge 1$ with $U_0\left(x\right) = 1$ (see, e.g., \cite[p. 376]{Olver}, \cite[p. 256]{NIST}). For $n=1,2,3$, we have
\[
U_1 \left( x \right) =  - \frac{5}{{24}}x^3  + \frac{1}{8}x,
\]
\[
U_2 \left( x \right) = \frac{{385}}{{1152}}x^6  - \frac{{77}}{{192}}x^4  + \frac{9}{{128}}x^2 ,
\]
\[
U_3 \left( x \right) =  - \frac{{85085}}{{82944}}x^9  + \frac{{17017}}{{9216}}x^7  - \frac{{4563}}{{5120}}x^5  + \frac{{75}}{{1024}}x^3 .
\]
For $U_4 \left( x \right)$, $U_5 \left( x \right)$ and $U_6 \left( x \right)$, see Bickley et al. \cite[p. xxxv]{Bickley}. If we write
\[
U_n \left( x \right) = \frac{\left( 2n \right)!}{2^{2n} n!}x^n\sum\limits_{k = 0}^n { u_{n,k} x^{2k} } ,
\]
then substitution into \eqref{eq74} gives
\[
u_{n,k}  = \frac{{2n + 4k - 1}}{{4\left( {2n - 1} \right)\left( {n + 2k} \right)}}\left[ {\left( {2n + 4k - 1} \right)u_{n - 1,k}  - \left( {2n + 4k - 5} \right)u_{n - 1,k - 1} } \right]
\]
for $n \ge 1$, $0 \le k \le n$ with $u_{0,0}=1$ and $u_{n,-1}=u_{n,n+1}=0$ for $n\geq 0$. This is equivalent to a result of Meijer \cite{Meijer}.

An analytic continuation argument in \eqref{eq13} yields the formula
\begin{equation}\label{eq91}
U_n \left( x \right) = \left( { - 1} \right)^n \frac{{x^n }}{{2^n n!}}\left[ {\frac{{d^{2n} }}{{dt^{2n} }}\left( {\frac{1}{2}\frac{{t^2 }}{{x\left( {t - \sinh t} \right) + \cosh t - 1}}} \right)^{n + \frac{1}{2}} } \right]_{t = 0} .
\end{equation}
There has been a recent interest in finding explicit formulas for the coefficients in asymptotic expansions of Laplace-type integrals (see \cite{Lopez}, \cite{Nemes}, \cite{Wojdylo1} and \cite{Wojdylo2}). There are two general formulas for these coefficients, one containing Potential polynomials and one containing Bell polynomials. We derive them here for the special case of the coefficients $U_n \left( x \right)$. Let
\[
x\left( {t - \sinh t} \right) + \cosh t - 1 = \sum\limits_{k = 0}^\infty  {a_k t^{k + 2} } ,
\]
so that
\[
a_{2k}  = \frac{1}{{\left( {2k + 2} \right)!}},\; a_{2k + 1}  =  - \frac{x}{{\left( {2k + 3} \right)!}} \; \text{ for } \; k\geq 0.
\]
Let $0 \leq j \leq k$ be integers and $\rho$ be a complex number. We define the Potential polynomials
\[
\mathsf{A}_{\rho ,k}  = \mathsf{A}_{\rho ,k} \left( {\frac{{a_1 }}{{a_0 }},\frac{{a_2 }}{{a_0 }}, \ldots ,\frac{{a_k }}{{a_0 }}} \right)
\]
and the Bell polynomials
\[
\mathsf{B}_{k,j}  = \mathsf{B}_{k,j} \left( {a_1 ,a_2 , \ldots ,a_{k - j + 1} } \right)
\]
via the expansions
\begin{equation}\label{eq92}
\left( {1 + \sum\limits_{k = 1}^\infty  {\frac{{a_k }}{{a_0 }}t^k } } \right)^\rho   = \sum\limits_{k = 0}^\infty  {\mathsf{A}_{\rho ,k} t^k } \; \text{ and } \; \mathsf{A}_{\rho ,k}  = \sum\limits_{j = 0}^k {\binom{\rho}{j}\frac{1}{a_0^j}\mathsf{B}_{k,j}} .
\end{equation}
Naturally, these polynomials can be defined for arbitrary power series with $a_0\neq 0$. It is possible to express the Potential polynomials with complex parameter in terms of Potential polynomials with integer parameter using the following formula of Comtet \cite[p. 142]{Comtet}
\begin{equation}\label{eq93}
\mathsf{A}_{\rho ,k}  = \frac{\Gamma \left( { - \rho  + k + 1} \right)}{k!\Gamma \left( { - \rho } \right)}\sum\limits_{j = 0}^k {\frac{\left( { - 1} \right)^j}{- \rho  + j}\binom{k}{j}\mathsf{A}_{j,k} } .
\end{equation}
With these notations we can write \eqref{eq91} as
\[
U_n \left( x \right) = \left( { - 1} \right)^n \frac{{x^n }}{{2^n n!}}\left[ {\frac{{d^{2n} }}{{dt^{2n} }}\left( {1 + \sum\limits_{k = 1}^\infty  {\frac{{a_k }}{{a_0 }}t^k } } \right)^{ - n - \frac{1}{2}} } \right]_{t = 0}  = \left( { - 1} \right)^n \frac{{\left( {2n} \right)!}}{{2^n n!}}x^n \mathsf{A}_{ - n - \frac{1}{2},2n} .
\]
Using \eqref{eq92} and \eqref{eq93} we find
\[
U_n \left( x \right) = \left( { - 1} \right)^n \left( {2x} \right)^n \sum\limits_{k = 0}^{2n} {\left( { - 1} \right)^k \frac{2^k \Gamma \left( {n + k + \frac{1}{2}} \right)}{\sqrt \pi  k!}\mathsf{B}_{2n,k} } 
\]
and
\[
U_n \left( x \right) = \left( { - 1} \right)^n \frac{{\left( {2x} \right)^n }}{{\left( {2n} \right)!}}\frac{{2\Gamma \left( {3n + \frac{3}{2}} \right)}}{{\sqrt \pi  }}\sum\limits_{k = 0}^{2n} {\frac{\left( { - 1} \right)^k}{2n + 2k + 1}\binom{2n}{k} \mathsf{A}_{k,2n} } .
\]
The quantities $\mathsf{B}_{n,k}$ and $\mathsf{A}_{k,n}$ appearing in these formulas may be computed from the recurrence relations
\[
\mathsf{B}_{n,k}  = \sum\limits_{j = 1}^{n - k + 1} {a_j \mathsf{B}_{n - j,k - 1} } \; \text{ and } \; \mathsf{A}_{k,n}  = \sum\limits_{j = 0}^{n} {\frac{a_j}{a_0}\mathsf{A}_{k - 1,n - j} }
\]
with $\mathsf{B}_{0,0} = \mathsf{A}_{0,0} = 1$, $\mathsf{B}_{j,0} = \mathsf{A}_{0,j} = 0$ $\left(j \ge 1\right)$, $\mathsf{B}_{j,1} = a_0 \mathsf{A}_{1,j} = a_j$ (see Nemes \cite{Nemes}). We remark that since
\[
x\left( {t - \sinh t} \right) + \cosh t - 1 = \frac{1-x}{2}\left( {e^t  - t - 1} \right) + \frac{1+x}{2}\left( {e^{ - t}  + t - 1} \right),
\]
it is possible to derive an expression for the Potential polynomials $\mathsf{A}_{k,2n}$ in terms of the $1$-associated Stirling numbers of the second kind (see, e.g., Howard \cite{Howard}), but we do not discuss the details here.

For the general theory of Potential polynomials and Bell polynomials, see Comtet \cite[pp. 133--153]{Comtet}. Some other formulas for the coefficients $U_n\left( x \right)$ involving Bell polynomials are given by L\'{o}pez and Pagola \cite{Lopez}. For another approach, using rearrangement of the asymptotic series, see Luke \cite{Luke}.

\subsection{The coefficients $d_{2n}$}

In 1952, Lauwerier \cite{Lauwerier} showed that the coefficients in asymptotic expansions of Laplace-type integrals can be calculated by means of linear recurrence relations. As an illustration of his method, he considered, inter alia, the coefficients $d_{2n}$. Define the sequence $P_0 \left( x \right), P_1 \left( x \right), P_2 \left( x \right),\ldots$ of polynomials via the recurrence
\[
P_n \left( x \right) =  - \sum\limits_{k = 1}^n {\frac{1}{{\left( {2k + 3} \right)!}}\int_0^x {P_{n - k} \left( t \right)dt} }
\]
with $P_0 \left( x \right) = 1$. Then the coefficients $d_{2n}$ can be recovered from the formula
\[
d_{2n}  = \frac{1}{{\Gamma \left( {\frac{2n + 1}{3}} \right)}}\int_0^{+\infty} {t^{\frac{{2n - 2}}{3}} e^{ - \frac{t}{6}} P_n \left( t \right)dt} .
\]
The first few are given explicitly by
\[
d_0  = 6^{\frac{1}{3}} ,\; d_2  =  - \frac{3}{10},\; d_4  = \frac{3}{140}6^{\frac{2}{3}},\; d_6  =  - \frac{6^{\frac{1}{3}} }{100},\; d_8  = \frac{10449}{2156000},\; d_{10}  =  - \frac{3639}{9100000}6^{\frac{2}{3}}.
\]
For higher coefficients, see Jentschura and L\"{o}tstedt \cite{Jentschura}.

It is possible to derive representations for the $d_{2n}$'s in terms of the Potential polynomials and Bell polynomials like for the polynomials $U_n \left( x \right)$. Let
\[
\sinh t - t = \sum\limits_{k = 0}^\infty  {b_k t^{k + 3} } ,
\]
that is
\[
b_{2k}  = \frac{1}{\left( {2k + 3} \right)!},\; b_{2k + 1}  = 0 \; \text{ for } \; k\geq 0.
\]
Employing this notation in \eqref{eq15} gives
\[
d_{2n}  = \frac{6^{\frac{2n + 1}{3}}}{\left( {2n} \right)!}\left[ {\frac{{d^{2n} }}{{dt^{2n} }}\left( {1 + \sum\limits_{k = 1}^\infty  {\frac{b_k}{b_0}t^k } } \right)^{ - \frac{2n + 1}{3}} } \right]_{t = 0}  = 6^{\frac{2n + 1}{3}} \mathsf{A}_{ - \frac{2n + 1}{3},2n} .
\]
Applying \eqref{eq92} and \eqref{eq93} (with $b_k$ in place of $a_k$) we obtain the formulas
\[
d_{2n}  = \frac{6^{\frac{2n + 1}{3}}}{\Gamma \left( {\frac{2n + 1}{3}} \right)}\sum\limits_{k = 0}^{2n} {\left( { - 1} \right)^k \frac{{6^k \Gamma \left( {\frac{2n + 1}{3} + k} \right)}}{{k!}} \mathsf{B}_{2n,k} } 
\]
and
\begin{equation}\label{eq94}
d_{2n}  = \frac{6^{\frac{2n + 1}{3}}}{\left(2n\right)!}\frac{3\Gamma \left( {4\frac{{2n + 1}}{3}} \right)}{\Gamma \left( {\frac{{2n + 1}}{3}} \right)}\sum\limits_{k = 0}^{2n} {\frac{\left( { - 1} \right)^k }{2n + 3k + 1}\binom{2n}{k}\mathsf{A}_{k,2n} } .
\end{equation}
The quantities $\mathsf{B}_{n,k}$ and $\mathsf{A}_{k,n}$ can be generated via the recurrence relations
\[
\mathsf{B}_{n,k}  = \sum\limits_{j = 1}^{n - k + 1} {b_j \mathsf{B}_{n - j,k - 1} } \; \text{ and } \; \mathsf{A}_{k,n}  = \sum\limits_{j = 0}^{n} {\frac{b_j}{b_0}\mathsf{A}_{k - 1,n - j} }
\]
with $\mathsf{B}_{0,0} = \mathsf{A}_{0,0} = 1$, $\mathsf{B}_{j,0} = \mathsf{A}_{0,j} = 0$ $\left(j \ge 1\right)$, $\mathsf{B}_{j,1} = b_0 \mathsf{A}_{1,j} = b_j$. Finally, we show that the Potential polynomials $\mathsf{A}_{k,2n}$ in \eqref{eq94} can be written in terms of the generalised Bernoulli polynomials $B_n^{\left( \kappa  \right)} \left(\lambda\right)$ which are defined by the exponential generating function
\[
\left( \frac{z}{e^z  - 1} \right)^\kappa  e^{\lambda z}  = \sum\limits_{n = 0}^\infty  {B_n^{\left( \kappa  \right)} \left(\lambda\right)\frac{z^n}{n!}} \; \text{ for } \; \left|z\right| < 2\pi.
\]
For basic properties of these polynomials, see Milne-Thomson \cite{Milne-Thomson} or N\"{o}rlund \cite{Norlund}. A straightforward computation gives
\begin{align*}
\mathsf{A}_{k,2n} & = \frac{1}{{\left( {2n} \right)!}}\left[ {\frac{{d^{2n} }}{{dt^{2n} }}\left( {6\frac{{\sinh t - t}}{{t^3 }}} \right)^k } \right]_{t = 0} = \frac{1}{{2\pi i}}\oint_{\left( {0^ +  } \right)} {\left( {6\frac{{\sinh z - z}}{{z^3 }}} \right)^k \frac{{dz}}{{z^{2n + 1} }}} \\
& = \frac{{6^k }}{{2\pi i}}\oint_{\left( {0^ +  } \right)} {\left( {\sum\limits_{j = 0}^k {\left( { - 1} \right)^{k - j} \binom{k}{j}z^{k - j} \sinh ^j z} } \right)\frac{{dz}}{{z^{2n + 3k + 1} }}}  = \sum\limits_{j = 0}^k {\left( { - 1} \right)^{k - j} \binom{k}{j}\frac{{6^k }}{{2\pi i}}\oint_{\left( {0^ +  } \right)} {\left( {\frac{{\sinh z}}{z}} \right)^j \frac{{dz}}{{z^{2n + 2k + 1} }}} } \\
& = \sum\limits_{j = 0}^k {\left( { - 1} \right)^{k - j} \binom{k}{j}\frac{{6^k }}{{2\pi i}}\oint_{\left( {0^ +  } \right)} {\left( {\frac{{2z}}{{e^{2z}  - 1}}} \right)^{ - j} e^{ - jz} \frac{{dz}}{{z^{2n + 2k + 1} }}} }  = \sum\limits_{j = 0}^k {\left( { - 1} \right)^{k - j} \binom{k}{j}\frac{{2^{2n + 2k} 6^k }}{{\left( {2n + 2k} \right)!}}B_{2n + 2k}^{\left( { - j} \right)} \left( { - \frac{j}{2}} \right)} ,
\end{align*}
and substitution into \eqref{eq94} provides a representation of the coefficients $d_{2n}$ in terms of the generalised Bernoulli polynomials.

An other formula for the coefficients $d_{2n}$ involving Bell polynomials is derived by L\'{o}pez and Pagola \cite{Lopez}.

\section{Auxiliary inequalities}\label{appendixb}

\subsection{Proof of the inequality \eqref{eq47}} It is enough to show that
\[
\left| {1 + re^{ - \frac{2}{3}\theta i} } \right|^2 \left| {1 + re^{\frac{2}{3}\left( {\pi  - \theta } \right)i} } \right|^2  \geq \begin{cases} \cos^2 \theta & \; \text{ if } \; { - \frac{\pi }{2} < \theta  < 0 \; \text{ or } \; \pi  < \theta  < \frac{{3\pi }}{2}} \\ 1 & \; \text{ if } \; {0 \le \theta  \le \pi }, \end{cases}
\]
for any $r>0$. First, we observe that
\[
\cos \left( {\frac{2}{3}\left( {\frac{\pi}{2}  - \theta } \right)} \right) \ge 0,\; \frac{1}{2} + \cos \left( {\frac{4}{3}\left({ \frac{\pi}{2}  - \theta } \right)} \right) \ge 0
\]
when $0 \leq \theta \leq \pi$, therefore
\[
\left| {1 + re^{ - \frac{2}{3}\theta i} } \right|^2 \left| {1 + re^{\frac{2}{3}\left( {\pi  - \theta } \right)i} } \right|^2  = 1 + 2r\left( {r^2  + 1} \right)\cos \left( {\frac{2}{3}\left( {\frac{\pi }{2} - \theta } \right)} \right) + 2r^2 \left( {\frac{1}{2} + \cos \left( {\frac{4}{3}\left( {\frac{\pi }{2} - \theta } \right)} \right)} \right) + r^4 \geq 1
\]
for any $0 \leq \theta \leq \pi$. To prove the other inequality, we note that
\begin{align*}
\left| {1 + re^{ - \frac{2}{3}\theta i} } \right|^2 \left| {1 + re^{\frac{2}{3}\left( {\pi  - \theta } \right)i} } \right|^2 & = \left( {1 + 2r\cos \left( {\frac{{2\theta }}{3}} \right) + r^2 } \right)\left( {1 + 2r\cos \left( {\frac{2}{3}\left( {\pi  - \theta } \right)} \right) + r^2 } \right)\\
& = \frac{3}{4}\left( {1 - r^2 } \right)^2  + \left( {\frac{{1 + r^2 }}{2} + 2r\cos \left( {\frac{{\pi  - 2\theta }}{3}} \right)} \right)^2 .
\end{align*}
On the other hand,
\[
\cos ^2 \theta  = \frac{{1 + \cos \left( {2\theta } \right)}}{2} = \frac{{1 - \cos \left( {3\frac{{\pi  - 2\theta }}{3}} \right)}}{2} = \frac{1}{2}\left( {1 + 3\cos \left( {\frac{{\pi  - 2\theta }}{3}} \right) - 4\cos ^3 \left( {\frac{{\pi  - 2\theta }}{3}} \right)} \right) .
\]
Hence, if we denote $q = \cos \left( {\frac{{\pi  - 2\theta }}{3}} \right)$, all we need to show is that for $r>0$ and $-\frac{1}{2} < q < \frac{1}{2}$ one has
\[
\frac{3}{4}\left( {1 - r^2 } \right)^2  + \left( {\frac{{1 + r^2 }}{2} + 2rq} \right)^2  \ge \frac{1}{2}\left( {1 + 3q - 4q^3 } \right)
\]
or, equivalently
\begin{equation}\label{eq46}
P\left( {q,r} \right) = 4q^3  + 8r^2 q^2  + \left( {4r^3  + 4r - 3} \right)q + \left( {2r^4  - 2r^2  + 1} \right) \ge 0.
\end{equation}
Now since it can be easily shown that
\[
P\left( { - \frac{1}{2},r} \right) = 2\left( {r - 1} \right)^2 \left( {r^2  + r + 1} \right) \ge 0,
\]
\[
P\left( {\frac{1}{2},r} \right) = 2r\left( {r^3  + r^2  + 1} \right) > 0,
\]
the only possibility for \eqref{eq46} to be violated is for the function $q \mapsto P\left( {q,r} \right)$ to take a minimum value at some point $-\frac{1}{2} < q^\ast < \frac{1}{2}$. Differentiating $P$ with respect to $q$ and solving the resulting quadratic equation for $q$, one can show that such a point can exist only if $0 < r \leq \frac{3}{4}$ or $\left(\frac{3}{4}\right)^{\frac{1}{3}} \leq r < 1$ and is explicitly given by
\[
q^\ast = q^\ast \left( r \right) = \frac{1}{6}\left( {\sqrt {\left( {4r - 3} \right)\left( {4r^3  - 3} \right)}  - 4r^2 } \right).
\]
Finally, it can be shown by standard single-variable methods that $P\left( {q^\ast,r} \right) \ge 0$ for the values of $q^\ast$ given above and any $r$ satisfying $0 < r \leq \frac{3}{4}$ or $\left(\frac{3}{4}\right)^{\frac{1}{3}} \leq r <1$.

\subsection{Proof of the inequality \eqref{eq48}} It is enough to show that
\[
\left| {1 - r^{\frac{2}{3}} e^{ - \frac{2}{3}\theta i} } \right|^2  \le \begin{cases} \left| {1 + r^2 e^{ - 2\theta i} } \right|^2 \csc^2 \left(2\theta\right) & \; \text{ if } \; \frac{\pi}{4} < \left|\theta\right| < \frac{\pi}{2} \\ \left| {1 + r^2 e^{ - 2\theta i} } \right|^2 & \; \text{ if } \; \left|\theta\right| \le \frac{\pi}{4}, \end{cases} 
\]
for any $r>0$. Suppose that $\left|\theta\right| \le \frac{\pi}{4}$, then
\[
\left| {1 + r^2 e^{ - 2\theta i} } \right|^2  - \left| {1 - r^{\frac{2}{3}} e^{ - \frac{2}{3}\theta i} } \right|^2  \ge 2r^{\frac{2}{3}} \cos \left( {\frac{{2\theta }}{3}} \right) + 2r^2 \cos \left( {2\theta } \right) + r^4  - r^{\frac{4}{3}}  \ge \sqrt 3 r^{\frac{2}{3}}  + r^4  - r^{\frac{4}{3}}  \ge 0,
\]
for $r>0$. Now we consider the case $\frac{\pi}{4} < \left|\theta\right| < \frac{\pi}{2}$. The inequality in question is equivalent to the assertion
\[
\left( {r^2  - 1} \right)^2  + 2r^2 \left( {1 + \cos \left( {2\theta } \right)} \right) + \left( {\left( {2\cos \left( {\frac{{2\theta }}{3}} \right) - r^{\frac{2}{3}} } \right)r^{\frac{2}{3}}  - 1} \right)\left(1-\cos^2 \left( {2\theta } \right)\right) \ge 0
\]
for any $r>0$. We can assume that $\frac{\pi}{4} < \theta < \frac{\pi}{2}$. By monotonicity
\[
\left( {2\cos \left( {\frac{{2\theta }}{3}} \right) - r^{\frac{2}{3}} } \right)r^{\frac{2}{3}}  - 1 \ge \left( {2\cos \left( {\frac{2}{3}\frac{\pi }{2}} \right) - r^{\frac{2}{3}} } \right)r^{\frac{2}{3}}  - 1 =  - r^{\frac{4}{3}}  + r^{\frac{2}{3}}  - 1 .
\]
Hence, if we denote $q = \cos \left( 2\theta \right)$, all we need to show is that for $r>0$ and $-1 < q < 0$ one has
\[
\left( {r^2  - 1} \right)^2  + 2r^2 \left( {1 + q} \right) + \left( { - r^{\frac{4}{3}}  + r^{\frac{2}{3}}  - 1} \right)\left( {1 - q^2 } \right) \ge 0
\]
or, equivalently
\[
Q\left(q,r\right) = \left( {r^{\frac{4}{3}}  - r^{\frac{2}{3}}  + 1} \right)q^2 + 2r^2 q + \left( {r^4  - r^{\frac{4}{3}}  + r^{\frac{2}{3}} } \right) \ge 0 .
\]
We have
\[
Q\left( { - 1,r} \right) = \left( {r^2  - 1} \right)^2  \ge 0,
\]
\[
Q\left( {0,r} \right) = r^4  - r^{\frac{4}{3}}  + r^{\frac{2}{3}}  \ge 0,
\]
If $0<r < 1$, the function $q \mapsto Q\left( {q,r} \right)$ has a minimum at
\[
q^{\ast} = q^{\ast} \left( r \right) =  - \frac{{r^2 }}{{r^{\frac{4}{3}}  - r^{\frac{2}{3}}  + 1}}
\]
and $Q\left( {q^{\ast},r} \right) \geq 0$. This completes the proof of the inequality \eqref{eq48}.

\subsection{Proof of formula \eqref{eq90}} Let $a>0$ be fixed. It is enough to prove that if $f$ and $g$ are defined by
\[
g\left( {x,a} \right) = \frac{{1 - \left( {a/x} \right)^{\frac{4}{3}} }}{{1 - \left( {a/x} \right)^2 }} = \frac{{1 - a^{\frac{4}{3}} }}{{1 - a^2 }} + \left( {x - 1} \right)f\left( {x,a} \right) \; \text{ for } \; x>0,
\]
then $\left|f\left( {x,a} \right)\right| < 2$. It is easy to show that $g'\left( {x,a} \right) \ge 0$, $\lim _{x \to 0^+} g'\left( {x,a} \right) =  + \infty$, $g'\left( {\frac{1}{2},a} \right) < 2$ and $g''\left( {x,a} \right) \le 0$. From these it follows that there is an $0< x^{\ast}\left(a\right) < \frac{1}{2}$ such that $g'\left( {x,a} \right) \ge 2$ for all $0 < x \leq x^{\ast}\left(a\right)$, and $g'\left( {x,a} \right) < 2$ for all $x > x^{\ast}\left(a\right)$. Hence, if $x>x^{\ast}\left(a\right)$ then by the Taylor formula with Lagrange remainder it follows that $\left| {f\left( {x,a} \right)} \right|  < 2$. Now, we restrict our attention to the interval $0< x \leq x^{\ast}\left(a\right) \left(<\frac{1}{2}\right)$. We have
\[
f'\left( {x,a} \right) = \left( {\frac{{g\left( {x,a} \right) - g\left( {1,a} \right)}}{{x - 1}}} \right)^\prime   = \frac{{g'\left( {x,a} \right)\left( {x - 1} \right) - \left( {g\left( {x,a} \right) - g\left( {1,a} \right)} \right)}}{{\left( {x - 1} \right)^2 }} .
\]
By assumption,
\[
g'\left( {x,a} \right)\left( {x - 1} \right) <  - \frac{1}{2}g'\left( {x,a} \right) \le  - 1
\]
and
\[
 - \left( {g\left( {x,a} \right) - g\left( {1,a} \right)} \right) = g\left( {1,a} \right) - g\left( {x,a} \right) \le g\left( {1,a} \right) - g\left( {0,a} \right) = g\left( {1,a} \right) \le 1 .
\]
Therefore $f'\left( {x,a} \right) \le 0$ for all $0< x \leq x^{\ast}\left(a\right)$, i.e., $f$ is non-increasing for $0< x \leq x^{\ast}\left(a\right)$. It follows that
\[
f\left( {x,a} \right) \le f\left( {0,a} \right) = g\left( {1,a} \right) \le 1 <2
\]
for $0< x \leq x^{\ast}\left(a\right)$. Finally,
\[
f\left( {x,a} \right) = \frac{{g\left( {1,a} \right) - g\left( {x,a} \right)}}{{1 - x}} > g\left( {1,a} \right) - g\left( {x,a} \right) \ge g\left( {1,a} \right) - g\left( {1,a} \right) = 0 >  - 2
\]
in the range $0< x \leq x^{\ast}\left(a\right)$.

\end{document}